\documentclass[12pt]{amsart}
\usepackage{amssymb, amscd, latexsym, color, amscd, float}
\usepackage{url}
\usepackage[final]{graphicx}

\usepackage{ulem,ifthen}
\newcommand{\ifemptythenelse}[3]{%
  \begingroup
    \def\dummy{#1}%
    \def\empty{}%
    \ifx\dummy\empty{#2}\else{#3}\fi
  \endgroup
  }
\DeclareRobustCommand{\change}[2][]{%
  \ifemptythenelse{#1}{%
    \ifemptythenelse{#2}{}{\begin{color}{blue}{{#2}}\end{color}}
  }{%
    \ifemptythenelse{#2}{%
      \sout{\begin{color}{red}{{#1}}\end{color}}%
    }{%
      \ifmmode
        \begin{color}{green}{{#2}}\end{color}%
      \else
      \begin{color}{blue}{{#2}}\end{color}%
      \footnote{was: \begin{color}{red}{{#1}}\end{color}}%
      \fi
    }%
  }%
}

\newcommand{\washere}[1]{}

\newtheorem{thm}{Theorem}[section]
\newtheorem{prop}[thm]{Proposition}
\newtheorem{lem}[thm]{Lemma}
\newtheorem{cor}[thm]{Corollary}

\newtheorem{defn}[thm]{Definition}
\newtheorem{rmk}[thm]{Remark}

\begin{document}



\title
[]
{\normalsize Equations in acylindrically hyperbolic groups
and verbal closedness}

\bigskip

\author{Oleg Bogopolski}
\address{{Sobolev Institute of Mathematics of Siberian Branch of Russian Academy
of Sciences, Novosibirsk, Russia}\newline {and D\"{u}sseldorf University, Germany}}
\email{Oleg$\_$Bogopolski@yahoo.com}

\begin{abstract}
We describe solutions of the equation $x^ny^m=a^nb^m$ in acylindrically hyperbolic groups (AH-groups), where $a,b$ are non-commensurable special loxodromic elements and $n,m$ are integers with sufficiently large common divisor.
Using this description and certain test words in AH-groups, we study the verbal closedness of AH-subgroups
in groups.

A subgroup $H$ of a group $G$ is called {\it verbally closed} if for any word $w(x_1,\dots, x_n)$
in variables $x_1,\dots,x_n$ and any element $h\in H$, the equation $w(x_1,\dots, x_n)=h$ has a solution in $G$ if and only if it has a solution in $H$.


{\bf Main Theorem:} Suppose that $G$ is a finitely presented group and $H$ is a finitely generated acylindrically 
hyperbolic subgroup of $G$ such that $H$ does not have nontrivial finite normal subgroups.
Then $H$ is verbally closed in $G$ if and only if $H$ is a retract of $G$.

The condition that $G$ is finitely presented and $H$ is finitely generated can be replaced by
the condition that $G$ is finitely generated over $H$ and $H$ is equationally noetherian.

As a corollary, we solve Problem 5.2 from the paper~\cite{MR} of Myasnikov and Roman'kov:
Verbally closed subgroups of torsion-free hyperbolic groups are retracts.

\end{abstract}

\maketitle

\setcounter{tocdepth}{1}
\tableofcontents

\bigskip





\section{Introduction}





In~1943, Neumann~\cite{Neumann_0} considered systems of equations over arbitrary groups and,
taking in mind field theory, introduced for groups such notions as adjoining of solutions,
algebraic and transcendent extensions.
Inspired by this paper, Scott~\cite{Scott} introduced the notion of algebraically closed groups in the class of all groups.
Since than the theory of equations over groups developed in two directions.

In the first direction one studies which types of equations are solvable over groups from certain classes
(e.g. over finite, residually finite, locally indicable, or torsion-free groups).
The branch which studies properties of algebraic sets in groups is called algebraic geometry over groups, see~\cite{BMR,Myas_Remesl}.
An exten\-sive list of problems and results in this area can be found in the survey
of Roman'kov~\cite{Roman'kov} of 2012 and in the recent paper of Klyachko and Thom~\cite{KT}.

In the second direction one studies properties of {\it algebraically, existentially,
and verbally closed groups} in certain overgroups or classes of groups (see Definition~\ref{alg,verb,retr} below and a general definition of $\frak{S}$-closedness suggested by Neumann in~\cite{Neumann_2}).
For problems and results in this area see the surveys of Leinen~\cite{Leinen} and Roman'kov~~\cite{Roman'kov}, and the papers~\cite{BL,KM,KMM,Macintyre,MR, Neumann_1,Neumann_2,RK1,RK2,RKK,Scott}. Note that this branch of group theory is closely related to logic in the form
of model theory and recursive functions,
see the book of Higman and Scott~\cite{Higman_Scott}, appendix A.4 in the book of Hodges~\cite{Hodges},
and the paper~\cite{JA}.

An important class of groups, where these two directions have a good chance for developing, is the class
of {\it acylindrically hyperbolic groups}.
These groups were implicitly studied in~\cite{BF,Bowditch,Ham,DOG,Sisto}
before they were formally defined by Osin in~\cite{Osin_1}.
In~\cite[Theorem 1.2]{Osin_1}, Osin proved that all definitions used in the above mentioned papers
are equivalent to his definition $({\rm AH}_1)$, see Section~3 below.

\newpage

The class of acylindrically hyperbolic groups is large. It includes non-(virtually cyclic) groups
that are hyperbolic relative to proper subgroups,
many 3-manifold groups, groups of deficiency at least~2, many groups acting on trees,
non-(virtually cyclic) groups acting properly on proper CAT(0)-spaces and containing rank-one elements, non-cyclic directly indecomposable right-angled Artin groups,
all but finitely many mapping class groups,
$Out(F_n)$ for $n\geqslant 2$,
and many other interesting groups; see the survey of Osin~\cite{Osin_3}.
Though this class is large, it can be universally studied by using methods of geometric group theory.


\medskip

In this paper, we first describe solutions of equations of type $x^ny^m=a^nb^m$ in acylindrically hyperbolic groups. Using this, we characterize (under some mild conditions)
verbally closed acylindrically hyperbolic subgroups of groups (see Theorems~\ref{1_prop 4.1} and~\ref{1_Theorem_2}).
In Corollary~\ref{solution_Myasnikov}, we solve Problem 5.2 from the paper~\cite{MR} of Myasnikov and Roman'kov:
{\it Verbally closed subgroups of torsion-free hyperbolic groups are retracts.}


The methods we use are mostly geometric, but we also use test words, whose construction
is combinatorial.


It is interesting to compare this corollary with known Theorems B and C from the appendix.
In particular, Theorem C says that a subgroup $H$ of a torsion-free hyperbolic group $G$
is existentially closed if and only if for any nontrivial element $g\in G$
there exists a retraction $\varphi:G\rightarrow H$ such that $\varphi(g)\neq 1$.

\medskip

{\bf Acknowledgements.} I am grateful to Denis Osin for useful discussions, in particular, for
pointing out results in~\cite{HO} on the extension of quasi-morphisms. I am also grateful to
David Bradley-Williams for helpful discussions on model theory.

\section{Main results}
\subsection{Algebraic closedness, verbal closedness, and retracts}\hfill

\noindent
Let $X=\{x_1,x_2,\dots \}$ be a countably infinite set of variables, and let $H$ be a group.
An {\it equation} with variables $x_1,\dots,x_n$ and constants from $H$ is an arbitrary
expression $W(x_1,\dots ,x_n;H)=1$, where $W(x_1,\dots ,x_n;H)$ is a word in the alphabet $\{x_1,\dots,x_n\}^{\pm}\cup H$. In other words $W(x_1,\dots ,x_n;H)$ lies in the free product
$F(X)\ast H$, where $F(X)$ is the free group  with basis $X$.\break If the left side of the equation
does not contain letters from $H$, we will omit~$H$ from this expression.
Let $G$ be an overgroup of $H$. A tuple $(g_1,\dots ,g_n)$ with components from $G$
is called a {\it solution of the equation $W(x_1,\dots ,x_n;H)=1$ in} $G$ if $W(g_1,\dots,g_n;H)=1$ in $G$.

We recall the definitions of algebraically (verbally) closed subgroups and retracts.

\begin{defn}\label{alg,verb,retr} {\rm Let $H$ be a subgroup of a group $G$.

\begin{enumerate}
\item[(a)] (see~\cite{Neumann_2,MR})
The subgroup $H$ is called {\it algebraically closed} in $G$ if
for any finite system of equations
$$
S=\{W_i(x_1,\dots,x_n;H)=1\,|\, i=1,\dots,m\}
$$
with constants from $H$ the following holds: if $S$ has a solution in $G$, then it has a solution in $H$.

\item[(b)] (see~\cite[Definition 1.1]{MR}) The subgroup $H$ is called {\it verbally closed} in $G$ if
for any word $W\in F(X)$ and any element $h\in H$ the following holds: if the equation $W(x_1,\dots,x_n)=h$
has a solution in $G$, then it has a solution in $H$.

\item[(c)] The subgroup $H$ is called a {\it retract} of $G$ if there is a homomorphism
$\varphi:G\rightarrow H$ such that $\varphi|_{H}={\rm id}$. The homomorphism $\varphi$ is called a {\it retraction}.
\end{enumerate}
}
\end{defn}

Obviously, if $H$ is a retract of $G$, then $H$ is algebraically closed in $G$.
The converse is true if $H$ is finitely generated and $G$ is finitely presented
(see~\cite[Proposition 2.2]{MR}).

Algebraic closedness implies verbal closedness, but the converse implication is not valid in general.
An example of a verbally closed but not algebraically closed finitely generated subgroup of a finitely generated virtually free group is given
in Remark~\ref{example}.



In~\cite{MR}, Myasnikov and Roman'kov proved that if $F$ is a free group of finite rank, then
every verbally closed subgroup $H$ of $F$ is a retract of $F$, and hence algebraically closed in $F$.
They write that not much is known
in general about verbally closed subgroups of a given group $G$ and raise two problems.

\medskip

\noindent
{\bf Problem 5.1 in~\cite{MR}.} What are the verbally closed subgroups of a free nilpotent group
of finite rank?

\medskip

\noindent
{\bf Problem 5.2 in~\cite{MR}.} Prove that verbally closed subgroups of a torsion-free hyperbolic group are retracts.

\medskip

Problem 5.1 was solved by Roman'kov and Khisamiev in~\cite{RK1}.
They proved the following. Let $\mathcal{N}_c$ be the variety of all nilpotent groups of class at most~$c$ and $N_{r,c}$
a free nilpotent group of finite rank $r$ and nilpotency class~$c$. A subgroup
$H$ of $N_{r,c}$ is verbally closed in $N_{r,c}$ if and only if $H$ is a free factor of $N_{r,c}$ in the variety $\mathcal{N}_c$ (equivalently, an algebraically closed subgroup, or a retract of $N_{r,c}$).
Some other results on verbal closedness can be found in~\cite{KM,KMM, Mazhuga_1,Mazhuga_2,Mazhuga_3}.

\medskip

Problem 5.2 is solved in this paper, see Corollary~\ref{solution_Myasnikov}.
Our main results are Theorems~\ref{1_prop 4.1} and~\ref{1_Theorem_2} about the verbal closedness
of acylindrically hyperbolic subgroups in groups.







\medskip









\begin{thm}\label{1_prop 4.1} {\rm (Theorem~\ref{prop 4.1})}
Suppose that $G$ is a finitely presented group and $H$ is a finitely generated acylindrically hyperbolic subgroup of $G$ such that $H$ does not have nontrivial finite normal subgroups.
Then $H$ is verbally closed in $G$ if and only if $H$ is a retract of $G$.
\end{thm}

\begin{rmk}\label{first_remark}\hspace*{-4.5mm}.
{\rm
\begin{enumerate}
\item[1)] The assumption that $H$ does not have nontrivial finite normal subgroups cannot be omitted (see example in Remark~\ref{example}).

\medskip

\item[2)] Any non-virtually-cyclic group that is hyperbolic relative to a (possibly infinite) collection of proper subgroups is acylindrically hyperbolic~\cite{Osin_3}. Thus, Theorem~\ref{1_prop 4.1} is applicable in this case;
see Corollary~\ref{corol_rel_hyp_1}.

\medskip

\item[3)]  Theorem~\ref{1_prop 4.1} implies a positive solution to Problem 5.2 in the case of finitely generated subgroups, see Corollary~\ref{corollary_hyp_finitely_generated_1}.
 To solve this problem in the general case we need a ``noetherian version'' of this theorem, see~Theorem~\ref{1_Theorem_2}.
\end{enumerate}
}
\end{rmk}









\noindent
Recall that a group $H$ is called {\it equationally noetherian} if
every system of equations with constants from $H$ and a finite number of variables
is equivalent to a finite subsystem, see~\cite{BMR}.
The equational noetherian property is important in the study of equations over groups, model theory of groups,
and other questions; see~\cite{BMR,{BMRom},Groves_2,JS,KM_0,Myas_Remesl,AOH,Sela_1,Sela_2,Sela_3}.

\begin{enumerate}
\item[(1)]  The Hilbert Basis Theorem implies that all linear groups over a commutative noetherian unitary ring, e.g., a field, are equationally noetherian, see the paper of
    Baumslag, Myasnikov and Remeslennikov~\cite[Theorem~B1]{BMR}. In particular, any free group is equationally noetherian (this was originally proved by Guba in~\cite{Guba}).


\item[(2)] Bryant proved in~\cite{Bryant} that every finitely generated abelian-by-nilpotent group is
equationally noetherian.

\item[(3)] Rigid solvable groups are equationally noetherian, see the paper of Romanovskii~\cite{Rom_1}. The case where the rigid solvable group is finitely generated was considered by Gupta and Romanovskii in~\cite{Gupta_Rom}.
Recall that a group $G$ is called {\it rigid} if it possesses a normal series of the form
$G=G_1>G_2>\dots >G_m>G_{m+1}=1$, where the quotients $G_i/G_{i+1}$ are abelian and torsion-free as right $\mathbb{Z}[G/G_i]$-modules.
In particular, free solvable groups are rigid and hence equationally noetherian.

\item[(4)]If $A$ and $B$ are equationally noetherian groups, then their free product $A\ast B$ is also equationally noetherian, see the preprint of Sela~\cite[Theorem 9.1]{Sela_3}.

\item[(5)] Hyperbolic groups are equationally noetherian. The torsion-free case was considered by Sela~\cite[Theorem 1.22]{Sela_2} and the general case by Reinfeldt and Weidmann~\cite[Corolary 6.13]{RW} (they expanded Sela's methods).


\item[(6)] In~\cite[Theorem D]{Groves_1}, Groves and Hull proved the following:\\
  Suppose that $G$ is a relatively hyperbolic group with respect to a finite collection of subgroups $\{H_1,\dots ,H_n\}$. Then $G$ is equationally noetherian if and only if each $H_i$ is equationally noetherian.


\end{enumerate}

Groves has announced in~\cite[Section 3]{Groves_1} that the mapping class group of a surface of finite type is equationally noetherian.

We say that a group $G$ is {\it finitely generated over a subgroup} $H$ if
there exists a finite subset $A\subset G$ such that $G=\langle A,H\rangle$.


\begin{thm}\label{1_Theorem_2}{\rm (Theorem~\ref{Theorem_2})}
Let $G$ be a group and let $H$ be a subgroup~of~$G$ such that $G$ is finitely generated over $H$. Suppose that
$H$ is equationally noetherian, acylindrically hyperbolic,
and does not have nontrivial finite normal subgroups.
Then $H$ is verbally closed in $G$ if and only if $H$ is a retract of~$G$.
\end{thm}


We deduce three corollaries from this theorem. For terminology concerning relatively hyperbolic groups see
the paper of Osin~\cite{Osin_0}.

\begin{cor}\label{1_noeth_rel_hyp}
{\rm (Corollary~\ref{noeth_rel_hyp_a})}
Let $G$ be a group and let $H$ be a subgroup~of~$G$ such that $G$ is finitely generated over $H$. Suppose that
$H$ is hyperbolic relative to a finite collection of equationally noetherian proper subgroups and does not have nontrivial finite normal subgroups. Then $H$ is verbally closed in $G$ if and only if $H$ is a retract of $G$.
\end{cor}

The special case of this corollary, where $H$ is a free group
(in this case the assumptions concerning $H$ are fulfilled automatically),
was considered in~\cite{KM}.

\begin{cor}\label{general_solution_Myasnikov}{\rm (Corollary~\ref{general_solution_Myasnikov_1})}
Let $G$ be a relatively hyperbolic group with respect to a finite collection of finitely generated
equationally noetherian subgroups. Suppose that $H$ is a non-parabolic subgroup
of $G$ such that $H$ does not have nontrivial finite normal subgroups.
Then $H$ is verbally closed in $G$ if and only if $H$ is a retract of $G$.
\end{cor}

The next corollary follows directly from the previous one. Indeed,
every hyperbolic group is relatively hyperbolic with respect to the trivial subgroup.
One can also deduce Corollary~\ref{solution_Myasnikov} directly from Theorem~\ref{1_Theorem_2}, see Remark~\ref{easier_proof}.

\begin{cor}\label{solution_Myasnikov} {\rm (Corollary~\ref{solution_Myasnikov_1})
(Solution to Problem 5.2 in~\cite{MR})}\\
Let $G$ be a hyperbolic group and $H$ be a subgroup of $G$.
Suppose that $H$ does not have nontrivial finite normal subgroups.
Then $H$ is verbally closed in $G$ if and only if $H$ is a retract of $G$.
\end{cor}

Note that the assumption that $H$ does not normalise a nontrivial finite subgroup of $H$ cannot be omitted (see example in Remark~\ref{example}).

\subsection{Solutions of certain equations in acylindrically hyperbolic groups}
\noindent
In the course of the proof of Theorem~\ref{1_prop 4.1}, we obtain a description of solutions
of the equation $x^ny^m=a^nb^m$ in acylindrically hyperbolic groups for non-commensurable special loxodromic elements $a,b$
and numbers $n,m$ with sufficiently large common divisor.

Suppose that $G$ is an acylindrically hyperbolic group {\it with respect to a generating set} $Y$, see
Definition~\ref{Definition_of_Osin}.
Then any loxodromic (with respect to~$Y$) element $g\in G $ is contained in a unique maximal virtually cyclic subgroup $E_G(g)$ of~$G$ (see~\cite[Lemma 6.5]{DOG}).
This subgroup is called the {\it elementary subgroup associated with} $g$.
An element $g\in G$ is called {\it special} with respect to $Y$ if it is loxodromic with respect to $Y$
and $E_G(g)=\langle g\rangle$ (see precise definitions in Section~3).
Two elements $a,b\in G$ of infinite order are called {\it commensurable}
if there exist $g\in G$ and $s,t\in \mathbb{Z}\setminus \{0\}$ such that $a^s=g^{-1}b^tg$.

\begin{prop}\label{1_lem 0.1}{\rm (Proposition~\ref{lem 0.1})}
Let $G$ be an acylindrically hyperbolic group with respect to a generating set $Z$.
Suppose that $a$ and $b$ are two non-commensurable special (with respect to $Z$)
elements of $G$.
Then there exists a generating set $Y$
containing $\mathcal{E}=\langle a\rangle\cup \langle b\rangle$ and there exists a number $N\in \mathbb{N}$ such that for all $n,m>N$ the following holds:

If $(c,d)$ is a solution of the equation $x^ny^m=a^nb^m$, then one of the following holds:

\begin{enumerate}
   \item[1)] $c$ and $d$ are loxodromic with respect to $Y$, and
$E_G(d)=E_G(c)$;

\item[2)] $c$ is loxodromic with respect to $Y$ and $d$ is elliptic, and $d^m\in E_G(c)$;

\item[3)] $d$ is loxodromic with respect to $Y$ and $c$ is elliptic, and
$c^n\in E_G(d)$;

\item[4)] $c$ and $d$ are elliptic with respect to $Y$, and one of the following holds:

\begin{enumerate}
\item[(a)]
$c$ is conjugate to $a$ and $d$ is conjugate to $b$;

\item[(b)]
$c$ is conjugate to $b$ and $d$ is conjugate to $a$, and $|n-m|\leqslant N$.
\end{enumerate}
\end{enumerate}
\end{prop}

Case 4) in this proposition does not give any information about the conjugators.
However, in the following special case we prove that the conjugators can be made equal and we give a simple description of solutions.

\begin{cor}\label{1_prop} {\rm (Corollary~\ref{prop})}
Let $G$ be an acylindrically hyperbolic group with respect to a generating set $S$.
Suppose that $a,b\in G$ are two non-commensurable special elements (with respect to $S$).
Then there exists a number $\ell=\ell(a,b)\in \mathbb{N}$ such that for all $n,m\in \mathcal{\ell}\mathbb{N}$, $n\neq m$, the equation $x^ny^m=a^nb^m$ is perfect, i.e. any solution of this equation in
$G$ is conjugate to $(a,b)$ by a power of $a^nb^m$.
\end{cor}

\newpage

\subsection{Uniform divergence of quasi-geodesics determined by loxodromic elements
in acylindrically hyperbolic groups}\hfill

Proposition~\ref{1_lem 0.1} is proved with the help
of the following two lemmas, which seem to be interesting for their own sake.
The first one says that the quasi-geodesics determined by two loxodromic elements
in acylindrically hyperbolic groups diverge uniformly.



\begin{lem}\label{lem 3.4_1} {\rm (Lemma~\ref{lem 3.4})}
Let $G$ be a group and let $X$ be a generating set of $G$.
Suppose that the Cayley graph $\Gamma(G,X)$ is hyperbolic and acylindrical.
Then there exists a constant $N_0>0$ such that
for any loxodromic (with respect to $X$) elements $c,d\in G$
with $E_G(c)\neq E_G(d)$ and for any $n,m\in \mathbb{N}$ we have that
$$|c^nd^m|_X>\frac{\min\{n,m\}}{N_0}.$$


\end{lem}



\medskip

\begin{lem}\label{lem 3.6_1} {\rm (Lemma~\ref{lem 3.6})}
Let $G$ be a group and let $X$ be a generating set of $G$.
Suppose that the Cayley graph $\Gamma(G,X)$ is hyperbolic and acylindrical.
Then there exists a constant $N_1>0$ such that
for any loxodromic (with respect to $X$) element $c\in G$, for
any elliptic element $e\in G\setminus E_G(c)$, and for any $n\in \mathbb{N}$ we have that
$$|c^ne|_X>\frac{n}{N_1}.$$


\end{lem}

We prove these lemmas with the help of the {\it periodicity theorem} for acylindrically hyperbolic groups,
see~\cite[Theorem~1.4]{Bog_1}. This theorem and relevant notions are reproduced in subsection~5.1 of the present paper. A special case of this theorem, where $G$ is a free group and $r=0$, can be found in the book of Adian~\cite{A} devoted to a solution of the Burnside problem (see statement 2.3 in Chapter I there).

\medskip

In a forthcoming paper we will prove the following theorem.

\medskip

\noindent
{\bf Theorem.} {\it Let $H$ be an acylindrically hyperbolic group without nontrivial finite normal subgroups.
Then for any group $G$ containing $H$ as a subgroup the following holds:
If $H$ is verbally closed in $G$, then $H$ is algebraically closed in~$G$.}

\bigskip

All actions of groups on metric spaces are assumed to be isometric in this paper.

\newpage

\section{Acylindrically hyperbolic groups}

We introduce general notation and recall some relevant definitions and statements
from the papers~\cite{Bog_1,DOG,Osin_1}.

\subsection{General notation} All generating sets considered in this paper are assumed to be symmetric, i.e. closed under taking inverse elements.
Let $G$ be a group generated by a subset $X$. For $g\in G$ let $|g|_X$ be the length of a shortest word in $X$ representing $g$. The corresponding metric
on $G$ is denoted by ${d}_X$ (or by ${d}$ if $X$ is clear from the context); thus ${d}_X(a,b)=|a^{-1}b|_X$. The right Cayley graph of $G$ with respect to $X$ is denoted by $\Gamma(G,X)$.
By a path $p$ in the Cayley graph we mean a combinatorial path; the initial and the terminal vertices of $p$ are denoted by $p_{-}$ and $p_{+}$, respectively.
The path inverse to $p$ is denoted by $\overline{p}$.
The label of $p$ is denoted by ${\bold {Lab}}(p)$;
we stress that the label is a formal word in the alphabet $X$.
The length of $p$ is denoted by $\ell(p)$.\\
Given a real number $K\geqslant 0$, two paths $p$ and $q$ in $\Gamma(G,X)$ are called {\it $K$-similar} if
${d}(p_{-},q_{-})\leqslant K$ and ${d}(p_{+},q_{+})\leqslant K$.

Recall that a path $p$ in $\Gamma(G,X)$ is called ($\varkappa,\varepsilon)$-{\it quasi-geodesic} for some $\varkappa\geqslant 1$,
$\varepsilon\geqslant 0$, if ${d}(q_{-},q_{+})\geqslant \frac{1}{\varkappa}\ell(q)-\varepsilon$ for any subpath $q$ of $p$.

The following remark is important. Suppose that $(X_{\lambda})_{\lambda \in \Lambda}$ is a collection of subsets of a group $G$ such that
$\underset{\lambda\in \Lambda}\cup X_{\lambda}$ generates $G$. The alphabet $\mathcal{X}=\underset{\lambda\in \Lambda}\sqcup X_{\lambda}$ determines the Cayley graph $\Gamma(G,\mathcal{X})$, where
two vertices may be connected by many edges. This happens if some element $x\in G$
belongs to subsets $X_{\lambda}$ and $X_{\mu}$ of $G$
for different $\lambda,\mu\in \Lambda$. In this case $\Gamma(G,\mathcal{X})$ contains two edges from $g$ to $gx$ for any
vertex $g$. The labels of these edges are different since they belong to disjoint subsets of the alphabet~$\mathcal{X}$; however these labels represent the same element $x$ in $G$.

The following notation will shorten the forthcoming proofs.
For $a,b,c\in \mathbb{R}$, we write $a\approx_{c}b$ if $|a-b|\leqslant c$.
Note that $a\approx_{c}b$ and $b\approx_{c_1}d$ imply $a\approx_{c+c_1}d$.

For a group $G$ and an element $a\in G$, we define a homomorphism $\widehat{a}:G\rightarrow G$
by the rule $\widehat{a}(g)=a^{-1}ga$. We also write $g^a$ for $a^{-1}ga$.

\subsection{Hyperbolic spaces}

Let $A,B,C$ be three points in a metric space $\frak{X}$.
Recall that the {\it Gromov product} of $A,B$ with respect to $C$ is the number
$$
(A,B)_C:=\frac{d(C,A)+d(C,B)-d(A,B)}{2}.
$$

We use the following definition of a $\delta$-hyperbolic space (see~\cite[Chapter III.H, Definition~1.16 and Proposition~1.17]{BH}).

For $\delta\geqslant 0$,
we say that a geodesic triangle $ABC$ in $\frak{X}$ is {\it $\delta$-thin} at the vertex $C$
if for any two points $A_1$ and $B_1$ on the sides $[C,A]$ and $[C,B]$ with $d(C,A_1)=d(C,B_1)\leqslant (A,B)_C$,
we have $d(A_1,B_1)\leqslant \delta$.

We say that a metric space $\frak{X}$ is {\it $\delta$-hyperbolic} if it is geodesic and every geodesic triangle in $\frak{X}$ is $\delta$-thin at each of its vertices.

\subsection{Definitions of acylindrically hyperbolic groups}

\begin{defn} {\rm (see~\cite{Bowditch} and Introduction in~\cite{Osin_1})
An action of a group $G$ on a metric space $S$ is called
{\it acylindrical}
if for every $\varepsilon>0$ there exist $R,N>0$ such that for every two points $x,y$ with $d(x,y)\geqslant R$,
there are at most $N$ elements $g\in G$ satisfying
$$
d(x,gx)\leqslant \varepsilon\hspace*{2mm}{\text{\rm and}}\hspace*{2mm} d(y,gy)\leqslant \varepsilon.
$$
}
\end{defn}

Given a generating set $X$ of a group $G$, we say that the Cayley graph $\Gamma(G,X)$ is
{\it acylindrical} if the left action of $G$ on $\Gamma(G,X)$ is acylindrical.
For Cayley graphs, the acylindricity condition can be rewritten as follows:
for every $\varepsilon>0$ there exist $R,N>0$ such that for any $g\in G$ of length $|g|_X\geqslant R$
we have
$$
\bigl|\{f\in G\,|\, |f|_X\leqslant \varepsilon,\hspace*{2mm} |g^{-1}fg|_X\leqslant \varepsilon \}\bigr|\leqslant N.
$$

Recall that an action of a group $G$ on a hyperbolic space $S$ is called {\it elementary} if the limit set
of $G$ on the Gromov boundary $\partial S$ contains at most 2 points.

\begin{defn}\label{Definition_of_Osin} {\rm (see~\cite[Definition 1.3]{Osin_1})
A group $G$ is called {\it acylindrically hyperbolic} if it satisfies one of the following equivalent
conditions:

\begin{enumerate}
\item[(${\rm AH}_1$)] There exists a generating set $X$ of $G$ such that the corresponding Cayley graph $\Gamma(G,X)$
is hyperbolic, $|\partial \Gamma (G,X)|>2$, and the natural action of $G$ on $\Gamma(G,X)$ is acylindrical.

\medskip

\item[(${\rm AH}_2$)] $G$ admits a non-elementary acylindrical action on a hyperbolic space.
\end{enumerate}
}
\end{defn}

In the case (AH$_1$), we also write that $G$ is {\it acylindrically hyperbolic with respect to $X$}.

\medskip

Recall the following useful lemma.

\begin{lem}\label{sup} {\rm (\cite[Lemma 5.1]{Osin_1})} For any group $G$ and any generating sets $X$ and $Y$
of $G$ such that
$$
\underset{x\in X}{\sup}|x|_Y<\infty\hspace*{3mm}{\rm and}\hspace*{3mm} \underset{y\in Y}{\sup}|y|_X<\infty,
$$
the following hold.

\begin{enumerate}
\item[(a)] $\Gamma(G,X)$ is hyperbolic if and only if $\Gamma(G,Y)$ is hyperbolic.

\item[(b)] $\Gamma(G,X)$ is acylindrical if and only if $\Gamma(G,Y)$ is acylindrical.
\end{enumerate}

\end{lem}

\subsection{Elliptic and loxodromic elements in acylindrically hyperbolic groups}

Let $G$ be a group acting on a metric space $S$.
Recall that the {\it stable norm} of an element $g\in G$ for this action is defined as $$||g||=\underset{n\rightarrow \infty}{\lim}\frac{1}{n}\,d(x,g^nx),$$ where $x$ is an arbitrary point in $S$,
see~\cite{CDP}. One verifies that this number is well-defined, independently of $x$, that it is a conjugacy invariant, and that $||g^k||=|k|\cdot||g||$ for all $k\in \mathbb{Z}$.
The following definition is standard.

\begin{defn}
{\rm
Given a group $G$ acting on a metric space $S$, an element $g\in G$ is called {\it elliptic}
if some (equivalently, any) orbit of $g$ is bounded, and {\it loxodromic} if the map
$\mathbb{Z}\rightarrow S$ defined by
$n\mapsto g^nx$ is a quasi-isometric embedding for some (equivalently, any) $x\in S$. That is,
for $x\in S$, there exist $\varkappa\geqslant 1$ and $\varepsilon\geqslant 0$ such that for any $n,m\in \mathbb{Z}$ we have
$$
d(g^nx,g^mx)\geqslant \frac{1}{\varkappa} |n-m|-\varepsilon.
$$

Let $X$ be a generating set of $G$.
We say that $g\in G$ is {\it elliptic (respectively loxodromic) with respect to $X$} if $g$ is elliptic (respectively loxodromic) for the canonical left action of $G$ on the Cayley graph $\Gamma(G,X)$.
If $X$ is clear from a context, we omit the words ``with respect to $X$''.

The set of all elliptic
(respectively loxodromic) elements of $G$ with respect to $X$ is denoted by ${\rm Ell}(G,X)$ (respectively by ${\rm Lox}(G,X))$.
}
\end{defn}

In the case of groups acting on hyperbolic spaces, there may be other types of actions
(see~\cite[Section 8.2]{Gromov} and~\cite[Section 3]{Osin_1}).


Bowditch~\cite[Lemma 2.2]{Bowditch} proved that every element of a group acting acylindrically on a hyperbolic space is either elliptic or loxodromic (see a more general statement in~\cite[Theorem 1.1]{Osin_1}).
Moreover, he proved there
that the infimum of the set of stable norms of all loxodromic elements for such an action is larger than zero
(we assume that $\inf \emptyset=+\infty$).

From this fundamental result, we deduced in~\cite[Corollary 2.12]{Bog_1} that, under certain assumptions, the quasi-geodesics associated with loxodromic elements (see Definition~\ref{Def_L}) have universal quasi-geodesic constants (see Corollary~\ref{qg}).


\begin{defn}\label{Def_L}
{\rm Let $G$ be a group and $X$ be a generating set of $G$.
For any two elements $u,v\in G$, we choose a path $[u,v]$ in $\Gamma(G,X)$ from $u$ to $v$ so that
$g[u,v]=[gu,gv]$ for any $g\in G$.
With any element $x\in G$ and any loxodromic element $g\in G$, we associate the bi-infinite quasi-geodesic
$$
L(x,g)=\overset{\infty}{\underset{i=-\infty}{\cup}}x[g^i,g^{i+1}].
$$
We have $L(x,g)=x\, L(1,g)$. The path $L(1,g)$ is called the {\it quasi-geodesic associated with $g$}.
}
\end{defn}

\begin{cor}\label{qg} {\rm (\cite[Corollary 2.12]{Bog_1})} Let $G$ be a group and $X$ be a generating set of $G$. Suppose that the Cayley graph $\Gamma(G,X)$ is hyperbolic and acylindrical.
Then there exist $\varkappa\geqslant 1$ and $\varepsilon\geqslant 0$ such that the following holds:

If an element $g\in G$ is loxodromic and shortest in its conjugacy class, then the quasi-geodesic $L(1,g)$
associated with $g$ is a $(\varkappa,\varepsilon)$-quasi-geodesic.
\end{cor}

Recall that any loxodromic element $g$ in an acylindrically hyperbolic group $G$ is contained in a
unique maximal virtually cyclic subgroup~\cite[Lemma 6.5]{DOG}. This subgroup, denoted by $E_G(g)$, is called the {\it elementary subgroup associated with $g$}; it can be described as follows (see equivalent definitions in~\cite[Corollary~6.6]{DOG}):
$$
\hspace*{12.5mm}\begin{array}{ll}
E_G(g)\! \!\! & =\{f\in G\,|\, \exists  n\in \mathbb{N}:  f^{-1}g^nf=g^{\pm n}\}\vspace*{3mm}\\
\! \!\! & =\{f\in G\,|\, \exists  k,m\in \mathbb{Z}\setminus \{0\}:  f^{-1}g^kf=g^{m}\}.
\end{array}\eqno{(3.1)}
$$

\begin{lem}\label{elem_index} {\rm (see~\cite[Lemma 6.8]{Osin_1})}
Suppose that a group $G$ acts acylindrically on a hyperbolic space $S$. Then there exists $L\in \mathbb{N}$
such that for every loxodromic element $g\in G$, $E_G(g)$ contains a normal infinite cyclic subgroup
of index~$L$.
\end{lem}

\begin{defn}\label{exact_defn_special}
{\rm  Suppose that $G$ is an acylindrically hyperbolic group.

\begin{enumerate}
\item[(a)] An element $g\in G$ is called {\it special} if there exists a generating set $X$ of $G$
such that

- $G$ is acylindrically hyperbolic with respect to $X$,

- $g$ is loxodromic with respect to $X$, and

- $E_G(g)=\langle g\rangle$.

\noindent
In this case $g$ is called {\it special with respect to~$X$}.

\medskip

\item[(b)] Elements $g_1,\dots ,g_k\in G$ are called {\it jointly special} if there exists a generating set $X$ of $G$
such that each $g_i$ is special with respect to $X$.
\end{enumerate}
}

\end{defn}

\medskip

Note that point (a) of this definition was already used in the case of relatively hyperbolic groups
(see comments in~\cite[Section 3]{OT}).



\medskip

The following theorem helps to verify, whether an acylindrical action of a group on a hyperbolic space
is elementary or not.

\medskip

\begin{thm}\label{check_elementary} {\rm (see~\cite[Theorem 1.1]{Osin_1})}
Let $G$ be a group acting acylindrically on a hyperbolic space. Then $G$
satisfies exactly one of the following conditions:
\begin{enumerate}
\item[(a)] $G$ has bounded orbits.

\item[(b)] $G$ is virtually cyclic and contains a loxodromic element.

\item[(c)] $G$ contains infinitely many loxodromic elements, whose limit sets are pairwise disjoint. In this case the action of $G$ is non-elementary and $G$ is acylindrically hyperbolic.
\end{enumerate}
\end{thm}


\subsection{Hyperbolically embedded subgroups}

Let $G$ be a group, $\{H_{\lambda}\}_{\lambda\in \Lambda}$ a collection of subgroups of $G$.
A subset $X$ of $G$ is called a {\it relative generating set of $G$ with respect to}
$\{H_{\lambda}\}_{\lambda\in \Lambda}$ if $G$ is generated by $X$ together with the union of all $H_{\lambda}$.
All relative generating sets are assumed to be symmetric.
We define
$$
\mathcal{H}=\bigsqcup_{\lambda\in\Lambda}H_{\lambda}.
$$


For the following two definitions, we assume that $X$ is a relative generating set of $G$ with respect to $\{H_{\lambda}\}_{\lambda\in \Lambda}$.

\begin{defn} {\rm (see~\cite[Definition 4.1]{DOG})
The group $G$ is called {\it weakly hyperbolic} relative to $X$ and $\{H_{\lambda}\}_{\lambda\in \Lambda}$ if the Cayley graph $\Gamma(G, X\sqcup \mathcal{H})$ is hyperbolic.
}
\end{defn}

\noindent
We consider the Cayley graph
$\Gamma(H_{\lambda},H_{\lambda})$ as a complete subgraph of $\Gamma(G,X\sqcup~\mathcal{H})$.

\begin{defn}
{\rm (see~\cite[Definition 4.2]{DOG})
For every $\lambda\in \Lambda$, we introduce a {\it relative metric}
$\widehat{d}_{\lambda}:H_{\lambda}\times H_{\lambda}\rightarrow [0,+\infty]$ as follows:

Let $a,b\in H_{\lambda}$. A path
in $\Gamma(G,X\sqcup \mathcal{H})$ from $a$ to $b$ is called {\it $H_{\lambda}$-admissible} if it has no edges in the subgraph $\Gamma(H_{\lambda},H_{\lambda})$.

The distance $\widehat{d}_{\lambda}(a,b)$ is defined to be the length of a shortest
{\it $H_{\lambda}$-admissible} path connecting $a$ to $b$ if such exists.
If no such path exists, we set $\widehat{d}_{\lambda}(a,b)=\!
\infty$.


}
\end{defn}

\begin{defn}\label{def_hyperb_embedd} {\rm (see~\cite[Definition 4.25]{DOG})
Let $G$ be a group, $X$ a symmetric subset of $G$. A collection of subgroups $\{H_{\lambda}\}_{\lambda\in \Lambda}$
of $G$ is called {\it hyperbolically embedded in $G$ with respect to $X$}
(we write $\{H_{\lambda}\}_{\lambda\in \Lambda}\hookrightarrow_h (G,X)$) if the following hold.

\begin{enumerate}
\item[(a)] The group $G$ is generated by $X$ together with the union of all $H_{\lambda}$ and the Cayley graph
$\Gamma(G,X\sqcup \mathcal{H})$ is hyperbolic.

\item[(b)] For every $\lambda\in \Lambda$, the metric space $(H_{\lambda},\widehat{d}_{\lambda})$ is
proper. That is, any ball of finite radius in $H_{\lambda}$ contains finitely many elements.
\end{enumerate}

\medskip

Further, we say that $\{H_{\lambda}\}_{\lambda\in \Lambda}$ is {\it hyperbolically embedded}
in $G$ and write $\{H_{\lambda}\}_{\lambda\in \Lambda}\hookrightarrow_h G$ if
$\{H_{\lambda}\}_{\lambda\in \Lambda}\hookrightarrow_h (G,X)$ for some $X\subseteq G$.
}
\end{defn}

It was proved in~\cite[Theorem 1.2]{Osin_1} that a group $G$ is acylindrically hyperbolic if and only if
it contains a proper infinite hyperbolically embedded subgroup.

\begin{lem}\label{hyp_embed_1} {\rm (see~\cite[Corollary~4.27]{DOG})} Let $G$ be a group, $\{H_{\lambda}\}_{\lambda\in \Lambda}$ a collection of subgroups of $G$, and $X,Y$ relative generating sets of $G$ with respect to $\{H_{\lambda}\}_{\lambda\in \Lambda}$.
Suppose that $|X\Delta Y|<\infty$.
Then
$\{H_{\lambda}\}_{\lambda\in \Lambda}\hookrightarrow_h (G,X)$
if and only if $\{H_{\lambda}\}_{\lambda\in \Lambda}\hookrightarrow_h (G,Y)$.
\end{lem}

There are examples which show that the condition $|X\Delta Y|<\infty$ cannot be replaced by the condition using supremum as in Lemma~\ref{sup}.


\begin{lem}\label{periph_subgr} {\rm  (see \cite[Proposition 4.33]{DOG})}
Suppose that $\{H_{\lambda}\}_{\lambda\in \Lambda}\hookrightarrow_h(G,X)$. Then, for each $\lambda\in \Lambda$,
we have $\{g\in G\,|\, |H_{\lambda}\cap g^{-1}H_{\lambda}g|=\infty\}\subseteq H_{\lambda}$.
\end{lem}

We use the following nontrivial theorem.

\begin{thm}\label{enlarging} {\rm (see~\cite[Theorem 5.4]{Osin_1})}
Let $G$ be a group, $\{H_{\lambda}\}_{\lambda\in \Lambda}$ a finite collection
of subgroups of $G$, $X$ a subset of $G$. Suppose that $\{H_{\lambda}\}_{\lambda\in \Lambda}\hookrightarrow_h (G,X)$.
Then there exists $Y\subseteq G$ such that $X\subseteq Y$ and the following conditions hold.

\begin{enumerate}
\item[(a)] $\{H_{\lambda}\}_{\lambda\in \Lambda}\hookrightarrow_h (G,Y)$. In particular, the Cayley graph
$\Gamma(G, Y\sqcup \mathcal{H})$ is hyperbolic.

\item[(b)] The action of $G$ on $\Gamma(G,Y\sqcup \mathcal{H})$ is acylindrical.
\end{enumerate}

\end{thm}

\medskip




\section{Preliminary statements}

The main aim of this section is to prove Lemmas~\ref{blue} and~\ref{small_conjug}.
These lemmas will be used in Section 5.

\begin{lem}\label{middle_points}
Let $ABC$ be a geodesic triangle in a $\delta$-hyperbolic space and let $A_1$ and $B_1$ be the middle points
of the sides $[A,C]$ and $[B,D]$, respectively. Then $d(A_1,B_1)\leqslant d(A,B)+2\delta$.
\end{lem}

\medskip

{\it Proof.} Denote $k=d(A,B)$.
First suppose that $A_1$ lies in the $\delta$-neighborhood of the side $[B,C]$, i.e. there
exists $D\in [B,C]$ such that $d(A_1,D)\leqslant \delta$. Using triangle inequalities, we deduce
$$
d(A,A_1)\approx_{k+\delta} d(B,D)
\hspace*{2mm}{\rm and}\hspace*{2mm} d(A_1,C)\approx_{\delta} d(D,C).
$$
Since $d(A,A_1)=d(A_1,C)$, we have $d(B,D)\approx_{k+2\delta} d(D,C)$.
Since $B_1$ is the middle point of $[B,C]$, we have
$d(D,B_1)\leqslant \frac{k}{2}+\delta$. Then
$$
d(A_1,B_1)\leqslant d(A_1,D)+d(D,B_1)\leqslant \frac{k}{2}+2\delta.
$$
The same estimation is valid in the case that $B_1$ lies in the $\delta$-neighborhood of $[A,C]$.
We consider the remaining case that $A_1$ and $B_1$ lie in the $\delta$-neighborhood of $[A,B]$.
Then there are points $A_2,B_2\in [A,B]$ such that $d(A_1,A_2)\leqslant \delta$ and
$d(B_1,B_2)\leqslant \delta$. This implies $$d(A_1,B_1)\leqslant d(A_1,A_2)+d(A_2,B_2)+d(B_2,B_1)\leqslant k+2\delta.$$
\hfill $\Box$

\medskip

\begin{lem}\label{Hausdorff} {\rm (see~\cite[Chapter III.H, Theorem 1.7]{BH})} For all $\delta\geqslant 0$, $\varkappa\geqslant 1$, $\epsilon\geqslant 0$, there exists a constant
$\mu=\mu(\delta,\varkappa,\epsilon)> 0$ with the following property:

If $\frak{X}$ is a $\delta$-hyperbolic space, $p$ is a $(\varkappa,\epsilon)$-quasi-geodesic in $\frak{X}$,
and $[x,y]$ is a geodesic segment joining the endpoints of $p$, then the Hausdorff distance between $[x,y]$
and the image of $p$ is at most $\mu$.
\end{lem}

An easy generalization of this statement is the following lemma.

\begin{lem}\label{close_qg} {\rm (see~\cite[Corollary 2.3]{Bog_1})} For any $\delta\geqslant 0$, $\varkappa\geqslant 1$, $\epsilon\geqslant 0$, $r\geqslant 0$ the following holds:

If $\frak{X}$ is a $\delta$-hyperbolic space, $p$ and $q$ are $(\varkappa,\epsilon)$-quasi-geodesics in $\frak{X}$ such that $\max\{d(p_{-},q_{-}),d(p_{+},q_{+})\}\leqslant r$, then
the Hausdorff distance between the images of $p$ and $q$ is at most $\eta(\delta,\varkappa,\epsilon,r)=
r+2\delta+2\mu$, where $\mu=\mu(\delta,\varkappa,\epsilon)$ is the constant from
Lemma~\ref{Hausdorff}.
\end{lem}





\medskip

\noindent
The following lemma can be deduced straightforward from the definition of Gromov product.

\begin{lem}\label{exercise}
Let $ABC$ be a geodesic triangle in a metric space and
let $P$ and $Q$ be points on its sides $[A,B]$ and $[B,C]$. Then
$d(P,Q)\geqslant d(P,B)+d(B,Q)-2(A,C)_B$.
\end{lem}

\begin{lem}\label{triangle}
Let $p,q$ be two $(\varkappa,\varepsilon)$-quasi-geodesics in a $\delta$-hyperbolic space such that $p_{+}=q_{-}$.
Then their concatenation $pq$ is a $(\varkappa,2\alpha+\beta)$-quasi-geodesic, where $\alpha$ is the Gromov product of $p_{-},q_{+}$ with respect to $p_{+}$ and
$\beta$ is a constant depending only on $\delta,\varkappa,\varepsilon$.
\end{lem}

\medskip

{\it Proof.} Let $r$ be a subpath of $pq$. We shall estimate $d(r_{-},r_{+})$ from below by using $\ell(r)$.
Consider only a nontrivial case: $r=p_1q_1$, where $p_1$ is a terminal subpath of $p$ and $q_1$
is an initial subpath of $q$.

Denote $A=p_{-}$, $B=p_{+}$, $C=q_{+}$, $P_1=r_{-}$ and $Q_1=r_{+}$.
By Lemma~\ref{Hausdorff}, there exists a constant $\mu=\mu(\delta,\varkappa,\varepsilon)$ and points $P\in[A,B]$ and $Q\in [B,C]$ such that
$$
d(P_1,P)\leqslant \mu,\hspace*{2mm} d(Q_1,Q)\leqslant \mu.
$$
By Lemma~\ref{exercise}, we have
$$d(P,Q)\geqslant d(P,B)+d(B,Q)-2\alpha.$$
Then
$$
\begin{array}{ll}
d(P_1,Q_1) & \geqslant d(P_1,B)+d(B,Q_1)-2\alpha-4\mu\vspace*{3mm}\\
& \geqslant \bigl(\frac{1}{\varkappa}\ell(p_1)-\varepsilon\bigr) + \bigl(\frac{1}{\varkappa}\ell(q_1)-\varepsilon\bigr)-2\alpha-4\mu
\vspace*{3mm}\\
& =\frac{1}{\varkappa}\ell(r)-(2\alpha+2\varepsilon+4\mu).
\end{array}
$$
We set $\beta=2\varepsilon+4\mu$. \hfill $\Box$

\medskip

The following lemma is a generalisation of Lemma~\ref{triangle}.

\medskip

\begin{lem}\label{concat_quasigeod} For any $\delta\geqslant 0$, $\alpha\geqslant 0$, $\varkappa\geqslant 1$, $\varepsilon\geqslant 0$ and any $m\in \mathbb{N}\cup \{0\}$, there exists $\varepsilon_0\geqslant 0$ such that the following holds.
Let  $\frak{X}$ be a $\delta$-hyperbolic space, $q=q_0q_1\dots q_mq_{m+1}$ a path in $\frak{X}$
such that
$q_0, q_1,\dots,q_m, q_{m+1}$ are $(\varkappa,\varepsilon)$-quasi-geodesic paths satisfying
\begin{enumerate}
\item[(1)] $\bigl((q_i)_{-},(q_{i+1})_{+}\bigr)_{(q_i)_{+}}<  \alpha$ for $i=0,\dots,m$;
\vspace*{1mm}
\item[(2)] $d\bigl((q_i)_{-},(q_i)_{+}\bigr)\geqslant 2\alpha +2m^2\delta$ for $i=1,\dots,m$.
\end{enumerate}
Then $q$ is a $(\varkappa,\varepsilon_0)$-quasi-geodesic path.
\end{lem}

\medskip

{\it Proof.} Induction by $m$. Induction basis $m=0$ was considered in Lemma~\ref{triangle}.
Suppose that $m\geqslant 1$ and proceed the induction step from $m-1$ to $m$. Let $q_0,q_1,\dots q_m,q_{m+1}$ be
$(\varkappa,\varepsilon)$-quasi-geodesic paths in a $\delta$-hyperbolic space $\frak{X}$,
which satisfy conditions (1) and (2) for some $\alpha\geqslant 0$. By Lemma~\ref{triangle}, the path $(q_0q_1)$
is a $(\varkappa,\varepsilon_1)$-quasi-geodesic for some
$\varepsilon_1=\varepsilon_1(\delta,\varkappa,\varepsilon,\alpha)$.
We may assume that $\varepsilon_1\geqslant \varepsilon$.
Then the paths $(q_0q_1),q_2,\dots ,q_m,q_{m+1}$ are $(\varkappa,\varepsilon_1)$-quasi-geodesic.
We show that these paths (after appropriate enumeration) satisfy conditions (1) and~(2)
with $\alpha_1:=\alpha+\delta$ instead of $\alpha$ and $m-1$ instead of $m$. Then we can apply induction.

To check condition (1), it suffices to show that
$\bigl((q_{0}q_1)_{-},(q_2)_{+}\bigr)_{(q_0q_1)_{+}}< \alpha_1$.
Denote $$A=(q_0)_{-},\hspace*{2mm} B=(q_1)_{-},\hspace*{2mm} C=(q_1)_{+},\hspace*{2mm} D=(q_2)_{+}.$$

\vspace*{-30mm}
\hspace*{-2.5mm}
\includegraphics[scale=0.7]{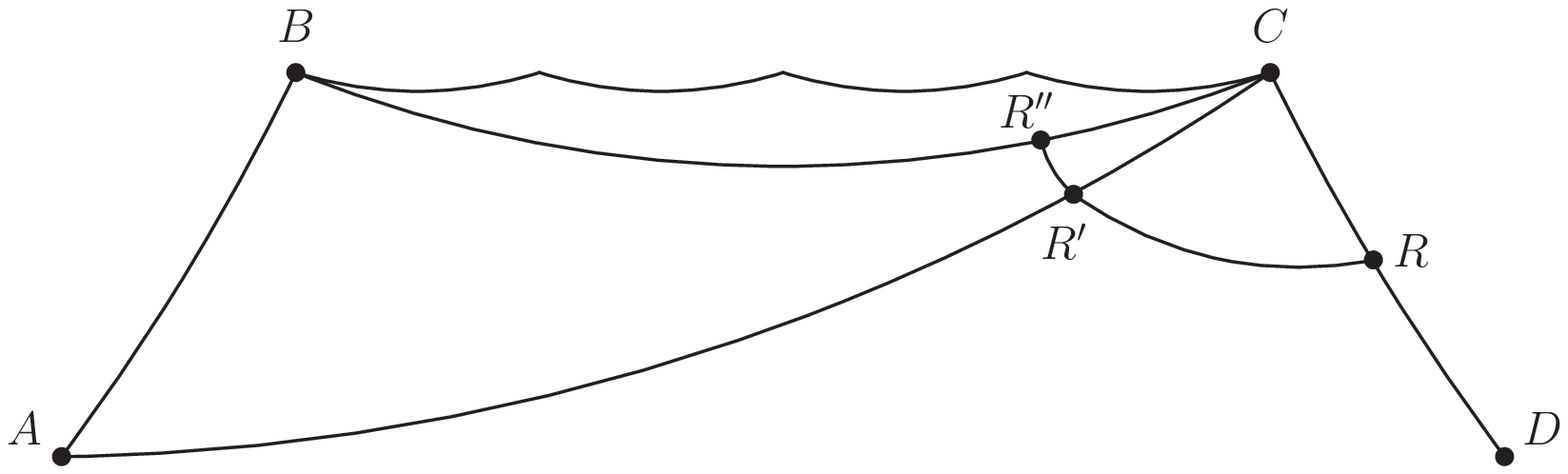}

\vspace*{-13.7cm}

\begin{center}
Fig. 1.
\end{center}

\medskip

Suppose the converse, i.e. $(A,D)_C\geqslant\alpha_1$.
Then there exist two points $R\in [C,D]$ and
$R'\in [C,A]$ such that
$$
d(C,R)=d(C,R')=\alpha_1\hspace*{2mm}{\text{\rm and}}\hspace*{2mm} d(R,R')\leqslant \delta,
$$
see Fig. 1.
We have
$$(A,B)_C=d(B,C)-(A,C)_B> (2\alpha+2m^2\delta)-\alpha
\geqslant \alpha+2\delta\geqslant \alpha_1.$$
Then there exists $R''\in [C,B]$ such that
$$d(C,R'')=\alpha_1 \hspace*{2mm}{\text{\rm and}}\hspace*{2mm} d(R'',R')\leqslant \delta.$$
Thus, $d(R,R'')\leqslant 2\delta$.
On the other hand, using Lemma~\ref{exercise}, we deduce
$$
d(R,R'')\geqslant d(C,R)+d(C,R'')-2(B,D)_C> 2\alpha_1-2\alpha=2\delta.
$$
A contradiction.

Condition (2) is fulfilled automatically since $2\alpha+2m^2\delta\geqslant 2\alpha_1+2(m-1)^2\delta$.
\hfill $\Box$


\begin{lem}\label{blue}
Let $G$ be a group and $X$ be a generating set of $G$. Suppose that the Cayley graph $\Gamma(G,X)$ is
hyperbolic and acylindrical.
Then there exist real numbers $\varkappa\geqslant 1, \varepsilon_0\geqslant 0$
and a number $n_0\in \mathbb{N}$ with the following property.

Suppose that $n\geqslant n_0$ and $c\in G$ is a loxodromic element.
Let $S(c)$ be the set of shortest elements in the conjugacy class
of $c$ and let $g\in G$ be a shortest element for which there exists $c_1\in S(c)$ with $c=g^{-1}c_1g$.
Then any path $p_0p_1\dots p_np_{n+1}$ in $\Gamma(G,X)$, where $p_0, p_1,\dots,p_n,p_{n+1}$
are geodesics with labels representing $g^{-1},c_1,\dots, c_1,g$, is a $(\varkappa, \varepsilon_0)$-quasi-geodesic.
In particular, $$|c^n|_X\geqslant \frac{1}{\varkappa}\bigl(n|c_1|_X+2|g|_X\bigr)-\varepsilon_0 \geqslant \frac{1}{\varkappa}n-\varepsilon_0.$$
\end{lem}

\medskip

{\it Proof.} Let $\delta\geqslant 0$ be a constant such that $\Gamma(G,X)$ is $\delta$-hyperbolic.
Let $n$ be an arbitrary positive integer.
We set $q_0=p_0$, $q_1=p_1p_2\dots p_n$, and $q_2=p_{n+1}$.
Then $q_1$ is a $(\varkappa,\varepsilon)$-quasi-geodesic, where $\varkappa$ and $\varepsilon$
are the constants from Corollary~\ref{qg}.
According to Lemma~\ref{Hausdorff}, the Hausdorff distance between $q_1$ and $[(q_1)_{-},(q_1)_{+}]$
is at most $\mu=\mu(\delta,\varkappa,\varepsilon)$.
We set $$n_0:=\lceil\varkappa(4\delta+2\mu+\varepsilon+2)\rceil.$$

Now we suppose that $n\geqslant n_0$. It suffices to show that the paths $q_0$, $q_1$, $q_2$ satisfy conditions (1) and (2) of Lemma~\ref{concat_quasigeod} for $\alpha=\delta+\mu+1$ and $m=1$.
Denote
$$A=(q_0)_{-},\hspace*{2mm} B=(q_1)_{-},\hspace*{2mm} C=(q_1)_{+},\hspace*{2mm} D=(q_2)_{+}.$$

\vspace*{-28mm}
\hspace*{-2.5mm}
\includegraphics[scale=0.7]{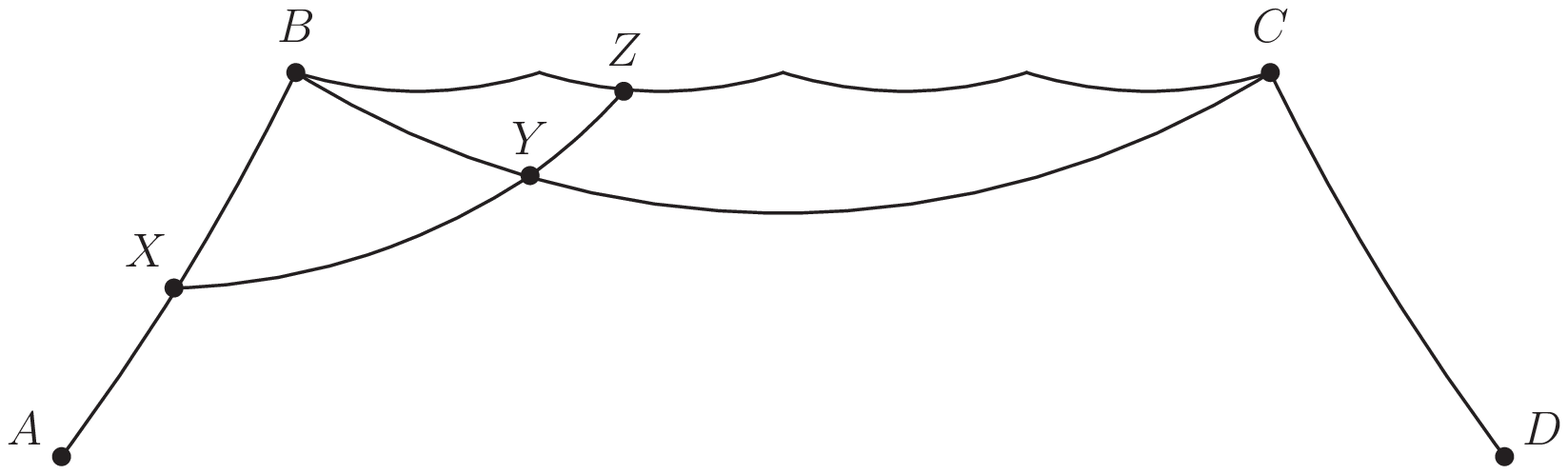}

\vspace*{-13.5cm}

\begin{center}
Fig. 2.
\end{center}

\bigskip

For (1), we shall check that
 $$(A,C)_B< \delta+\mu+1\hspace*{2mm}{\rm and}\hspace*{2mm} (B,D)_C< \delta+\mu+1.$$

Because of symmetry, we check only the first inequality. To the contrary, suppose that
$(A,C)_B\geqslant \delta+\mu+1$. Then there exist vertices $X\in [B,A]$ and $Y\in [B,C]$
such that $d(B,X)=d(B,Y)=\lfloor\delta+\mu+1\rfloor$ and $d(X,Y)\leqslant \delta$.
Since the Hausdorff distance between $q_1$ and $[B,C]$ is at most $\mu$,
there exists a vertex $Z\in q_1$ such that $d(Y,Z)\leqslant \mu$. Then
$$
\begin{array}{ll}
d(A,Z) & \leqslant d(A,X)+d(X,Y)+d(Y,Z)\vspace*{2mm}\\
       & \leqslant d(A,X)+\delta+\mu<d(A,X)+d(B,X)=d(A,B).
\end{array}\eqno{(4.1)}
$$

\medskip

We have $Z\in p_i=[Bc_1^i,Bc_1^{i+1}]$ for some $0\leqslant i\leqslant n-1$. Let $a$ and $b$ be the elements
of $G$ representing the labels of $[Bc_1^i,Z]$ and $[Z,Bc_1^{i+1}]$, respectively. Then $c_1=ab$. We set $c_2=ba$ and
$f=a^{-1}c_1^{-i}g$.
Then $c_2$ is also shortest in the conjugacy class of $c$ and we have $g^{-1}c_1g=f^{-1}c_2f$.
Since $Z=Bc_1^ia$ and $B=Ag^{-1}$, we have $f=Z^{-1}A$. Then  $$|f|_X=d(Z,A)\overset{(4.1)}{<}d(B,A)=|g|_X,$$ that contradicts the choice of $g$.
Condition (2) follows from the above mentioned fact that $q_1$ is a $(\varkappa,\varepsilon)$-quasi-geodesic:
$$d((q_1)_{-},(q_1)_{+})\geqslant \frac{1}{\varkappa}\ell(p_1p_2\dots p_n)-\varepsilon\geqslant \frac{n}{\varkappa}-\varepsilon\geqslant \frac{n_0}{\varkappa}-\varepsilon\geqslant 2\alpha+2\delta.
$$

Thus, by Lemma~\ref{concat_quasigeod}, $q_0q_1q_2=p_0p_1\dots p_np_{n+1}$ is a $(\varkappa, \varepsilon_0)$-quasi-geodesic for some universal constants
$\varkappa$, $\varepsilon_0$.
In particular, $$|c^n|_X\geqslant \frac{1}{\varkappa}\ell(p_0p_1\dots p_np_{n+1})-\varepsilon_0
=\frac{1}{\varkappa}\bigl(n|c_1|_X+2|g|_X\bigr)-\varepsilon_0\geqslant \frac{1}{\varkappa}n-\varepsilon_0.$$
\hfill $\Box$




\medskip

\begin{lem}\label{small_conjug}
For every $\delta\geqslant 0$, there exists $\varepsilon_1=\varepsilon_1(\delta)\geqslant 0$ such that the following holds.
Suppose that the Cayley graph of a group $G$ with respect to a generating set $X$ is $\delta$-hyperbolic
for some integer $\delta\geqslant 0$.
Let $a,b\in G$ be conjugate elements satisfying $|a|_X\geqslant |b|_X+4\delta+2$.
Then there exist $x,y\in G$ with the following properties:

\begin{enumerate}
\item[(a)] $a=x^{-1}yx$;
\vspace*{1mm}
\item[(b)] $|y|_X\in \{|b|_X+4\delta+1,|b|_X+4\delta+2\}$;
\vspace*{1mm}
\item[(c)] any path $q_0q_1q_2$ in $\Gamma(G,X)$, where $q_0,q_1,q_2$ are geodesics with labels representing $x^{-1},y,x$, is
a $(1,\varepsilon_1)$-quasi-geodesic.
\end{enumerate}
\end{lem}

\medskip

{\it Proof.}
In the set $S:=\{(y,x)\in G\times G\,|\, a=x^{-1}yx,\hspace*{2mm} |y|_X\leqslant |b|_X+4\delta +2\}$,
we choose a pair $(y,x)\in S$ with minimal $|x|_X$. Clearly, (a) is valid.
We claim that (b) is valid.

Suppose the contrary, i.e. $|y|_X\leqslant  |b|_X+4\delta$. Since $|a|_X\geqslant |b|_X+4\delta+2$,
we have $y\neq a$, hence $x\neq 1$.
We write $x=x_1x_2\dots x_n$ with $x_i\in X^{\pm}$, $i=1,\dots,n$, and $n=|x|_X$.
Then $(x_1^{-1}yx_1, x_2x_3\dots x_n)\in S$. A contradiction to minimality of $|x|_X$.

Now we verify that (c) is valid for some $\varepsilon_1$ depending only on $\delta$.
Let $q=q_0q_1q_2$ be a path in $\Gamma(G,X)$ such that its subpaths $q_0$, $q_1$, $q_2$ are geodesics with labels
representing $x^{-1}$, $y$, $x$. Without loss of generality, we may assume that $(q_0)_{-}=1$. First we show that conditions (1) and (2) of Lemma~\ref{concat_quasigeod} are satisfied for  $\alpha=\delta+1$ and $m=1$.

To the contrary, suppose that condition (1) of this lemma is not valid, say
$\bigl((q_1)_{-},(q_2)_{+}\bigr)_{(q_1)_{+}}\geqslant  \delta+1$, i.e.
$\Bigl(x^{-1},x^{-1}yx\Bigr)_{x^{-1}y}\geqslant \delta+1$. Because of $G$-invariance of Gromov product,
we have $$\Bigl(1,yx\Bigr)_{y}\geqslant \delta+1.$$
 Then there exist expressions $y=v_1v_2$, $x=u_1u_2$ such that $|y|_X=|v_1|_X+|v_2|_X$, $|x|_X=|u_1|_X+|u_2|_X$,
$|v_2|_X=|u_1|_X=\delta+1$, and $|v_2u_1|_X\leqslant \delta$.
Then the pair $( v_2v_1,v_2u_1u_2)$ lies in $S$ and
$$
|v_2u_1u_2|_X\leqslant |v_2u_1|_X+|u_2|_X\leqslant \delta+|u_2|_X< |u_1|_X+|u_2|_X=|x|_X.
$$
A contradiction to minimality of $|x|_X$. Thus, condition (1) of Lemma~\ref{concat_quasigeod} is valid.
Condition (2) of this lemma is also valid, since
$$
d((q_1)_{-},(q_1)_{+})=|y|_X\overset{\rm (b)}{\geqslant}  |b|_X+4\delta +1\geqslant 4\delta+2=2\alpha+2\delta.
$$
Recall that the paths $q_0,q_1,q_2$ are geodesic and hence $(1,0)$-quasi-geodesic.
Then, by Lemma~\ref{concat_quasigeod},
their concatenation $q_0q_1q_2$ is a $(1,\varepsilon_1)$-quasi-geodesic for some $\varepsilon_1$ depending only on $\delta$.
\hfill $\Box$

\medskip

\section{Uniform divergence of quasi-geodesics associated with loxodromic elements in acylindrically hyperbolic groups}

The main aim of this section is to prove Lemmas~\ref{lem 3.4} and~\ref{lem 3.6}.
We use Theorem~\ref{acylindric} proved by the author in~\cite{Bog_1}, see Theorem~1.4 there.
For convenience, we recall all necessary notions in the following subsection.

\subsection{A periodicity theorem for acylindrically hyperbolic groups}

\begin{defn}
{\rm
Let $Y$ be a generating set of $G$.
Given a loxodromic element $a\in G$ and an element $x\in G$, consider the bi-infinite path $L(x,a)$ in $\Gamma(G,Y)$
obtained by connecting consequent points $\dots, xa^{-1},x,xa,\dots $ by geodesic paths so that,
for all $n\in \mathbb{Z}$, the path $p_n$ connecting $xa^n$ and $xa^{n+1}$ has the same label as the path $p_0$
connecting $x$ and $xa$.
The paths $p_n$ are called {\it $a$-periods} of $L(x,a)$, and the vertices $xa^i$, $i\in \mathbb{Z}$, are called the {\it phase vertices} of $L(x,a)$. For $k\in \mathbb{N}$, we say that a subpath $p\subset L(x,a)$ {\it contains $k$ $a$-periods} if there exists $n\in \mathbb{Z}$ such that $p_np_{n+1}\dots p_{n+k-1}$ is a subpath of $p$.
}
\end{defn}

Let $G$ be a group and $X$ a generating set of $G$. Suppose that the Cayley graph $\Gamma(G,X)$
is hyperbolic and that $G$ acts acylindrically on $\Gamma(G,X)$. In~\cite{Bowditch}
Bowditch proved that the infimum of the set of stable norms (see Section~3 of the present paper) of all loxodromic elements of $G$ with respect to $X$ is a positive number.
We denote this number by ${\bf inj}(G,X)$ and call it the {\it injectivity radius} of $G$ with respect to $X$.

\begin{thm}\label{acylindric} {\rm (\cite[Theorem~1.4]{Bog_1})}
Let $G$ be a group and $X$ a generating set of $G$. Suppose that the Cayley graph $\Gamma(G,X)$
is hyperbolic and that $G$ acts acylindrically on $\Gamma(G,X)$.
Then there exists a constant $\mathcal{C}>0$ such that the following holds.

Let $a,b\in G$ be two loxodromic elements which are shortest in their conjugacy classes
and such that $|a|_X\geqslant |b|_X$.
Let $x,y\in G$ be arbitrary elements and $r$ an arbitrary non-negative integer. We set $f(r)=\frac{2r}{{\bf inj}(G,X)}+\mathcal{C}$.

If a subpath $p\subset L(x,a)$ contains at least $f(r)$ $a$-periods and lies in the $r$-neighborhood of $L(y,b)$,
then there exist $s,t\neq 0$ such that $(x^{-1}y)b^s(y^{-1}x)=a^t$.
In particular, $a$ and $b$ are commensurable.
\end{thm}




\medskip

\vspace*{-30mm}
\hspace*{-5mm}
\includegraphics[scale=0.7]{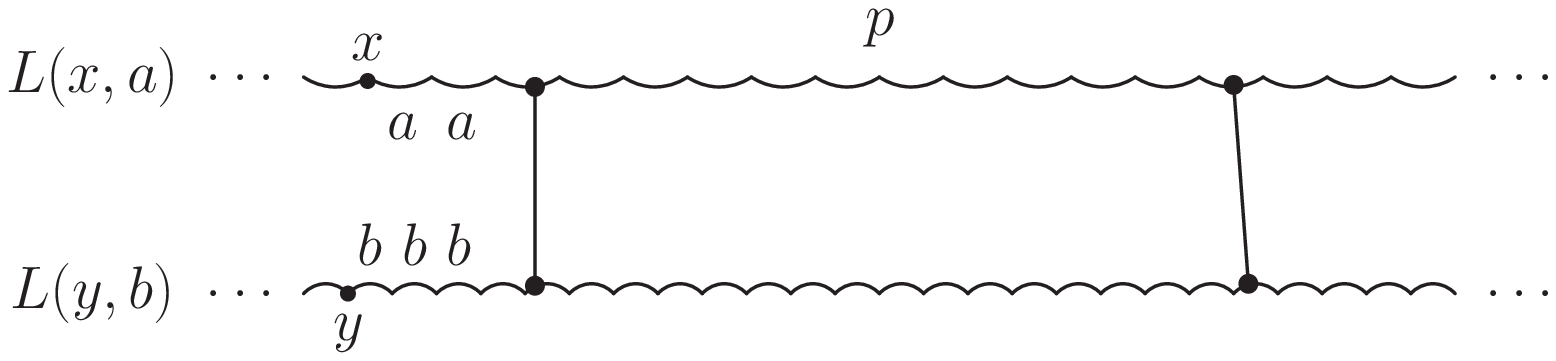}

\vspace*{-14.5cm}

\begin{center}
Fig. 3. Illustration to Theorem~\ref{acylindric}.
\end{center}

\medskip

\begin{rmk}\label{appropriate_conjugator}
{\rm
1) The main issue of Theorem~\ref{acylindric} is that the function $f$ is linear and does not depend on $|a|_X$ and $|b|_X$. Another point is that $X$ is allowed to be infinite.

2) In the conclusion of Theorem~\ref{acylindric}, we can write $z^{-1}b^sz=a^t$, where $z\in G$ is the element corresponding to the label of any path from any phase vertex of $L(y,b)$ to any phase vertex of $L(x,a)$.
}
\end{rmk}

\subsection{Loxodromic-loxodromic case}

\begin{lem}\label{lem 3.4}
Let $G$ be a group and let $X$ be a generating set of $G$.
Suppose that the Cayley graph $\Gamma(G,X)$ is hyperbolic and acylindrical.
Then there exists a constant $N_0>0$ such that
for any loxodromic (with respect to $X$) elements $c,d\in G$
with $E_G(c)\neq E_G(d)$ and for any $n,m\in \mathbb{N}$ we have that
$$|c^nd^m|_X>\frac{\min\{n,m\}}{N_0}.$$
\end{lem}



\medskip

{\it Proof.} Let $\delta\geqslant 0$ be a constant such that $\Gamma(G,X)$
is $\delta$-hyperbolic. We use the following constants:

{\small $\bullet$} $\varkappa\geqslant 1$, $\varepsilon_0\geqslant 0$ and $n_0\in \mathbb{N}$
are the constants from Lemma~\ref{blue};

{\small $\bullet$} $\mu=\mu(\delta,\varkappa,\varepsilon_0)$, see Lemma~\ref{Hausdorff};

{\small $\bullet$} $\mathcal{C}$ is the constant from Theorem~\ref{acylindric}.

We show that the lemma is valid for
$$
N_0=\max\bigl\{n_0, \frac{4000(\mu +\delta+1)}{{\bf inj}(G,X)}+4\mathcal{C}\bigr\}.\eqno{(5.1)}
$$

Suppose to the contrary that there exist two loxodromic elements
$c,d\in G$ satisfying $E_G(c)\neq E_G(d)$ and there exist $n,m\in \mathbb{N}$
such that
$$
\min\{n,m\}\geqslant N_0 k,\hspace*{2mm} {\rm where}\hspace*{2mm} k=|c^nd^m|_X.\eqno{(5.2)}
$$
Clearly, $k\neq 0$.

\medskip

For any $g\in G$, let $S(g)$ be the set of shortest elements in the conjugacy class of $g$.
Let $u\in G$ be a shortest element for which there exists $c_1\in S(c)$ with the property $c=u^{-1}c_1u$.
Let $v\in G$ be a shortest element for which there exists $d_1\in S(d)$ with the property $d=v^{-1}d_1v$.







Without loss of generality, we assume that $$|c_1|_X\geqslant |d_1|_X.\eqno{(5.3)}$$

Denote $w=d^{-m}c^{-n}$. Then $|w|_X=k$ and we have the equation
$$
u^{-1}\underbrace{c_1c_1\dots c_1}_nuv^{-1}\underbrace{d_1\dots d_1d_1}_{m}vw=1.\eqno{(5.4)}
$$

\medskip

\noindent
Consider a geodesic $(n+m+5)$-gon $\mathcal{P}=p_0(p_1p_2\dots p_n)p_{n+1}\bar{q}_{m+1}(\bar{q}_m\dots \bar{q}_2\bar{q}_1)\bar{q}_0h$
in the Cayley graph $\Gamma(G,X)$ such that the labels of its sides correspond
to the elements in the left side of equation (5.4) (see Fig. 4). In particular,
the labels of the paths $q_0, q_1, q_2,\dots ,q_m,q_{m+1}$ correspond to the elements
$v^{-1},d_1^{-1}, d_1^{-1},\dots ,d_1^{-1},v$.

\medskip

\vspace*{-25mm}
\hspace*{-12mm}
\includegraphics[scale=0.7]{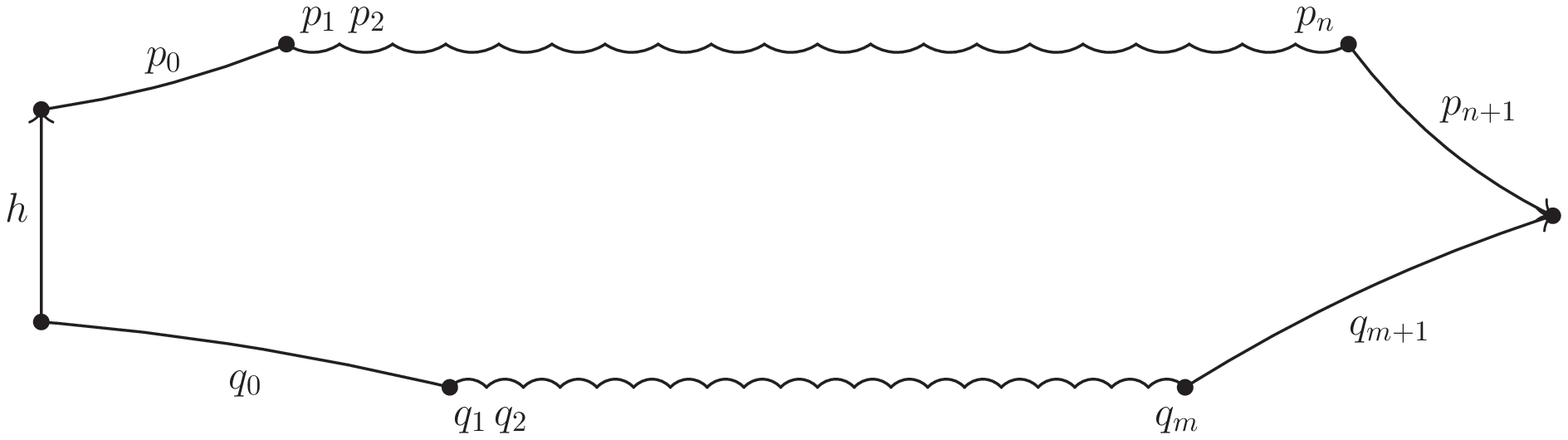}

\vspace*{-12.5cm}

\begin{center}
Fig. 4.
\end{center}

\bigskip

Since $\min\{n,m\}\geqslant N_0k\geqslant N_0\geqslant n_0$, we have by Lemma~\ref{blue} that
the paths $p:=p_0p_1p_2\dots p_np_{n+1}$ and $q=q_0q_1q_2\dots q_mq_{m+1}$ are $(\varkappa,\varepsilon_0)$-quasi-geodesics.
In the following claims we use the following constants.
$$
K=k+36\mu+26\delta,\hspace*{3mm}
K_1=K+4\mu+4\delta,\hspace*{3mm}\\
K_2=10K_1+2\mu+2\delta.
$$


{\bf Claim 1.} The quasi-geodesic $p_1p_2\dots p_n$ contains $n\geqslant 4f(K_2)$\,\, $c_1$-periods.
The quasi-geodesic $q_1q_2\dots q_m$ contains $m\geqslant 4f(K_2)$\,\, $d_1$-periods.

\medskip

{\it Proof.} The claim follows straightforward from the definition of function
$f(r)$ in Theorem~\ref{acylindric} and from (5.1) and (5.2).\hfill $\Box$

\medskip

We have in mind to apply Theorem~\ref{acylindric} to these quasi-geodesics or to their parts. However, we cannot claim that the first quasi-geodesic is contained
in the $K_2$-neighborhood of the second one.

Let $P$ be the middle point of the quasi-geodesic $p_1p_2\dots p_n$ and let
$Q$ be the middle point of the quasi-geodesic $q_1q_2\dots q_m$.

\medskip

{\bf Claim 2.} We have
$$d(P,Q)\leqslant K.\eqno{(5.5)}$$

\medskip

{\it Proof.}
Denote $A=(p_0)_{-}$, $B=(p_1)_{-}$, $C=(p_{n+1})_{-}$, $D=(p_{n+1})_{+}$, $E=(q_0)_{-}$, see Fig.~5.


\vspace*{-25mm}
\hspace*{-12mm}
\includegraphics[scale=0.7]{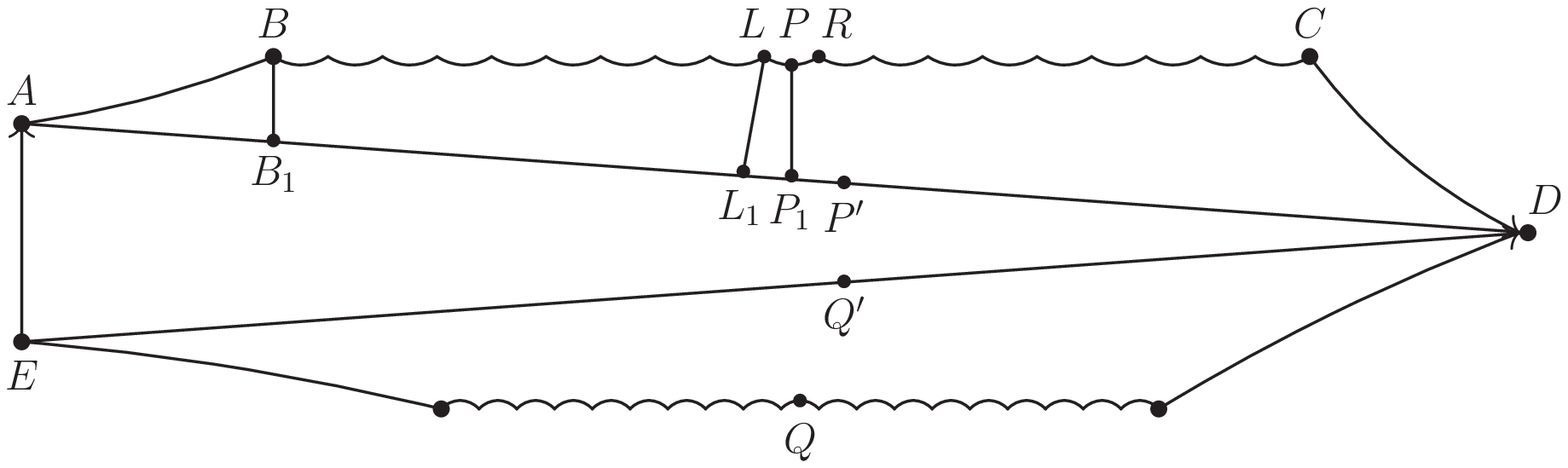}

\vspace*{-13cm}

\begin{center}
Fig. 5.
\end{center}

\medskip

By Lemma~\ref{Hausdorff}, there exists a point $P_1\in [A,D]$ such that
$$d(P,P_1)\leqslant \mu.\eqno{(5.6)}$$
The point $P$ divides the path $p$ into two halves. In the following we define two points $L,R\in p$.
If $n$ is even, we set $L=P=R$. If $n$ is odd, we set $L=(p_{\lceil\frac{n}{2}\rceil})_{-}$
and $R=(p_{\lceil\frac{n}{2}\rceil})_{+}$.
By Lemma~\ref{close_qg} applied to the first half of the quasi-geodesic $p$ and the geodesic $[A,P_1]$,
there exists a point $L_1\in [A,P_1]$ such that
$$d(L,L_1)\leqslant d(P,P_1)+ 2(\mu+\delta)\leqslant 3\mu+2\delta.\eqno{(5.7)}$$
Applying Lemma~\ref{close_qg} once more, we obtain that there exists a point $B_1\in [A,L_1]$ such that
$$d(B,B_1)\leqslant d(L,L_1)+2(\mu+\delta)\leqslant 5\mu+4\delta.\eqno{(5.8)}$$
Using triangle inequality several times, we deduce
$$
\begin{array}{ll}
d(A,P_1) & =  d(A,B_1)+d(B_1,L_1)+d(L_1,P_1)\vspace*{2mm}\\
& \approx_T d(A,B)+d(B,L)+d(L,P),
\end{array}
$$
where $T=2d(B,B_1)+2d(L,L_1)+d(P,P_1)$. It follows from (5.6)-(5.8) that $T\leqslant 17\mu+12\delta$. Hence,
$$
d(A,P_1)\approx_{17\mu+12\delta} d(A,B)+d(B,L)+d(L,P).\eqno{(5.9)}
$$
Analogously, we have
$$
d(D,P_1)\approx_{17\mu+12\delta} d(D,C)+d(C,R)+d(R,P).\eqno{(5.10)}
$$
Since the three summands in (5.9) are equal to the three summands in (5.10), we deduce that
$$d(A,P_1)\approx_{34\mu+24\delta} d(D,P_1).$$
Let $P'$ be the middle point of $[A,D]$. Then $d(P_1,P')\leqslant 17\mu+12\delta$.
We have
$$
d(P,P')\leqslant d(P,P_1)+d(P_1,P')\leqslant 18\mu+12\delta.\eqno{(5.11)}
$$
Analogously, if $Q'$ is the middle point of $[E,D]$, then
$$
d(Q,Q')\leqslant 18\mu+12\delta.\eqno{(5.12)}
$$
By Lemma~\ref{middle_points} applied to the geodesic triangle $ADE$, we have
$$
d(P',Q')\leqslant d(A,E)+2\delta=k+2\delta.\eqno{(5.13)}
$$
Now the claim follows from (5.11)-(5.13).\hfill $\Box$

\medskip

We consider the decomposition $p_0p_1p_2\dots p_np_{n+1}=\alpha_0\alpha_1\alpha_2\alpha_3\alpha_4\alpha_5$, where $\alpha_0=p_0$, $\alpha_1=p_1\dots p_{\lfloor\frac{n}{4}\rfloor}$,
$\alpha_4=p_{n-\lfloor\frac{n}{4}\rfloor+1}\dots p_n$, $\alpha_5=p_{n+1}$,
and $\alpha_2$ and $\alpha_3$ are determined by the condition $(\alpha_2)_{+}=P=(\alpha_3)_{-}$.
We also consider the decomposition $q_0q_1\dots q_mq_{m+1}=\gamma_1\gamma_2\gamma_3\gamma_4$,
where $\gamma_1=q_0$, $\gamma_4=q_{m+1}$, and $\gamma_2$ and $\gamma_3$ are determined by the condition $(\gamma_2)_{+}=Q=(\gamma_3)_{-}$, see Fig. 6.

\medskip

\vspace*{-25mm}
\hspace*{-10mm}
\includegraphics[scale=0.7]{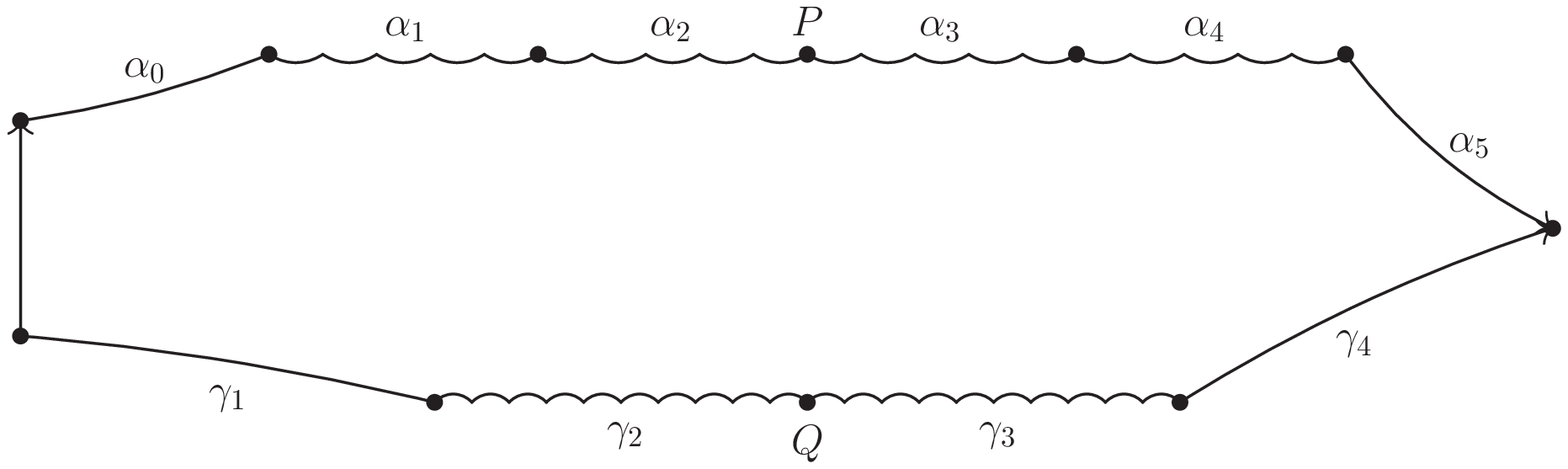}

\vspace*{-13cm}

\begin{center}
Fig. 6.
\end{center}

\medskip

{\bf Claim 3.} There are decompositions $\gamma_1\gamma_2=\beta_0\beta_1\beta_2$ and $\gamma_3\gamma_4=\beta_3\beta_4\beta_5$ such that
$\alpha_i$ and $\beta_i$ are $K_1$-similar for $i=0,\dots ,5$.
In particular, the Hausdorff distance between $\alpha_i$ and $\beta_i$ is at most $K_2$ for $i=0,\dots ,5$.
\medskip

{\it Proof.} Because of symmetry, we show only that the first decomposition exists.
We have
$$
\begin{array}{ll}
d((\alpha_0\alpha_1\alpha_2)_{-},(\gamma_1\gamma_2)_{-}) & = k,\vspace*{2mm}\\  d((\alpha_0\alpha_1\alpha_2)_{+},(\gamma_1\gamma_2)_{+}) & \overset{(5.5)}{\leqslant} K.
$$
\end{array}
$$
By Lemma~\ref{close_qg} applied to the $(\varkappa,\varepsilon_0)$-quasi-geodesics $\alpha_0\alpha_1\alpha_2$ and $\gamma_1\gamma_2$, there exists a point $U\in \gamma_1\gamma_2$ such that
$d((\alpha_1)_{+},U)\leqslant K+2\mu+2\delta$, see Fig.~7. Applying this lemma once more, we obtain
a point $V$ on the segment of $\gamma_1\gamma_2$ from $E$ to $U$ such that
$d((\alpha_1)_{-},V)\leqslant K+4\mu+4\delta=K_1$.
The points $V$ and $U$ divide the path $\gamma_1\gamma_2$ into three consecutive subpaths.
We denote them $\beta_0,\beta_1,\beta_2$.
By construction, the paths $\alpha_i$ and $\beta_i$ are $K_1$-similar.
Again by Lemma~\ref{close_qg}, the Hausdorff distance between $\alpha_i$ and $\beta_i$ is at most
$K_1+2\mu+2\delta\leqslant K_2$.
\hfill $\Box$

\medskip

\noindent
It follows from $\gamma_1\gamma_2=\beta_0\beta_1\beta_2$ that either $\beta_2$ is a subpath of $\gamma_2$, or
$\beta_0\beta_1$ is an initial subpath of $\gamma_1$.
Analogously, it follows from $\gamma_3\gamma_4=\beta_3\beta_4\beta_5$ that either
$\beta_3$ is a subpath of $\gamma_3$, or $\beta_4\beta_5$ is a terminal subpath of $\gamma_4$.

\medskip

\vspace*{-25mm}
\hspace*{-10mm}
\includegraphics[scale=0.7]{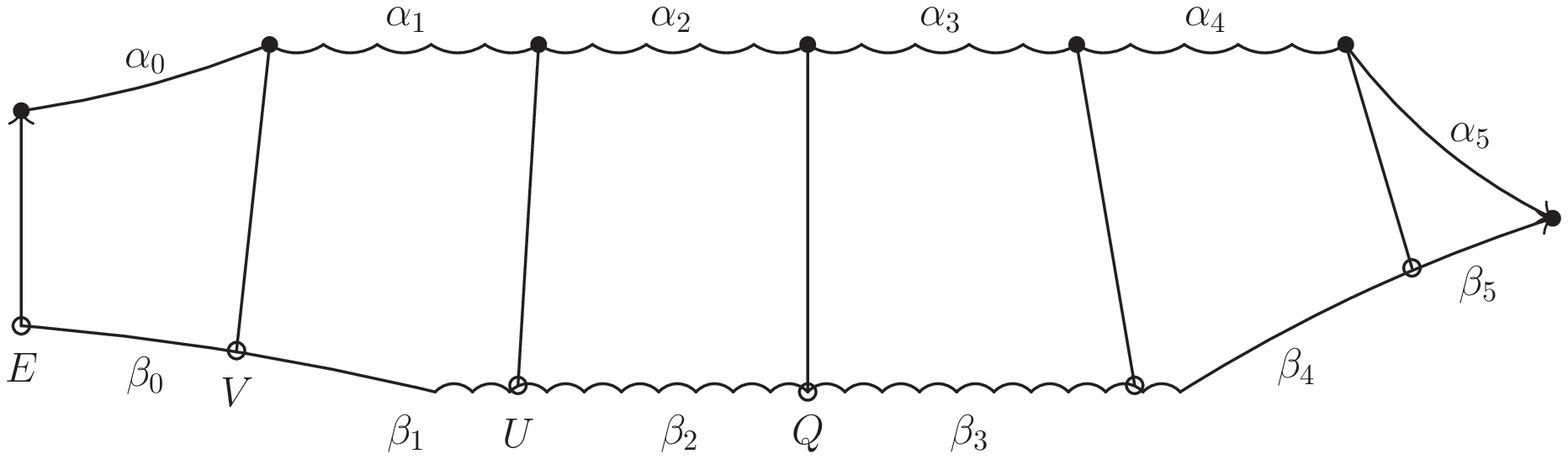}

\vspace*{-13cm}

\begin{center}
Fig. 7.
\end{center}

\medskip

{\it Case 1.} Suppose that $\beta_2$ is a subpath of $\gamma_2$, see Fig. 7.

Observe that $\alpha_2$ and $\beta_2$ satisfy assumptions of Theorem~\ref{acylindric}.
Indeed, in this case $\alpha_2$ and $\beta_2$ are subpaths of the quasi-geodesics $L((p_n)_{+},c_1)$ and
$L((q_m)_{+},d_1^{-1})$, respectively, $|c_1|_X\geqslant |d_1|_X$ by (5.3),
the Hausdorff distance between $\alpha_2$ and $\beta_2$ is at most $K_2$ by Claim~3, and
$\alpha_2$ contains at least $f(K_2)$ $c_1$-periods by Claim 1.

Hence, by~Theorem~\ref{acylindric} and Remark~\ref{appropriate_conjugator},
there exist $s,t\neq 0$ such that
$$
(uv^{-1})d_1^s(vu^{-1})=c_1^t.\eqno{(5.14)}
$$
Indeed, $vu^{-1}$ is the element which corresponds to the label of the path $q_{m+1}\overline{p_{n+1}}$ from the phase vertex $(q_m)_{+}$ to the phase vertex $(p_n)_{+}$.
From (5.14) and from $c=u^{-1}c_1u$, $d=v^{-1}d_1v$, we deduce $c^t=d^s$. Therefore $E_G(c)=E_G(d)$.
A contradiction to assumptions of the lemma.

\medskip

{\it Case 2.} Suppose that $\beta_3$ is a subpath of $\gamma_3$.

This case can be considered analogously.
Thus, it remains to consider the following case.

\medskip

{\it Case 3.} Suppose that
$\beta_0\beta_1$ is an initial subpath of $\gamma_1$ and $\beta_4\beta_5$ is a terminal subpath of $\gamma_4$,
see Fig. 8.

\vspace*{-25mm}
\hspace*{-10mm}
\includegraphics[scale=0.7]{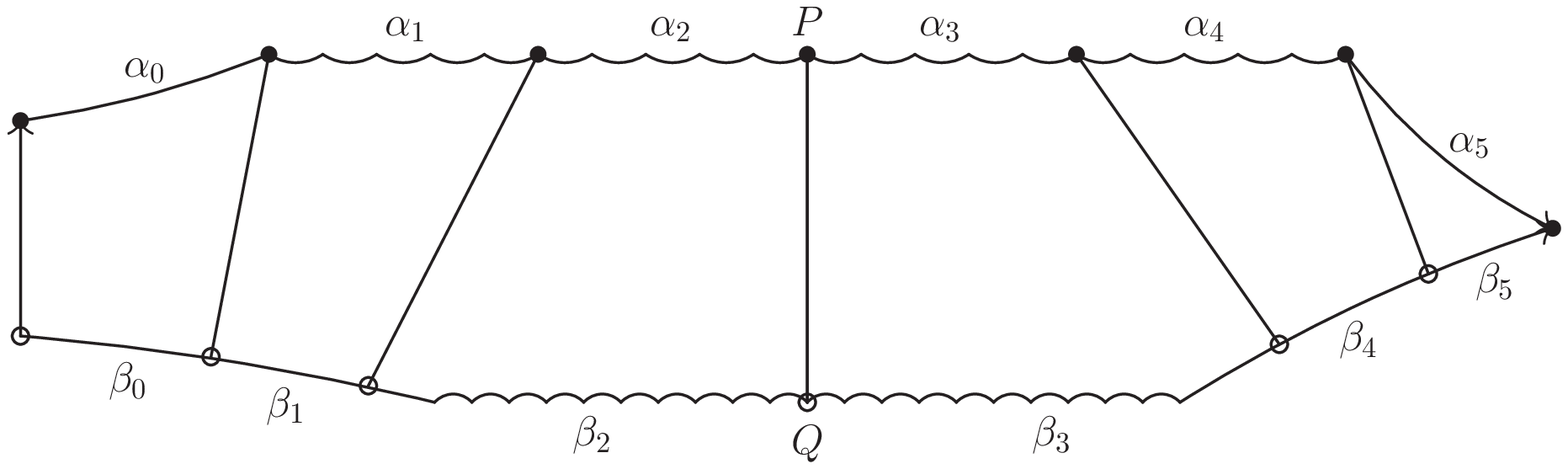}

\vspace*{-13cm}

\begin{center}
Fig. 8.
\end{center}

\medskip

Recall that the geodesics $\gamma_1$ and $\gamma_4$ have mutually inverse labels.
Therefore, there exist $g\in G$ such that $\gamma_1=g(\overline{\gamma_4})$. We denote
$$
\alpha_0'=g(\overline{\alpha_5}),\hspace*{3mm} \alpha_1'=g(\overline{\alpha_4}),\hspace*{3mm} \beta_0'=g(\overline{\beta_5}),\hspace*{3mm} \beta_1'=g(\overline{\beta_4}),
$$
see Fig. 9. Then $\beta_0\beta_1$ and $\beta_0'\beta_1'$ are initial segments of $\gamma_1$.

\bigskip

\vspace*{-25mm}
\hspace*{-10mm}
\includegraphics[scale=0.7]{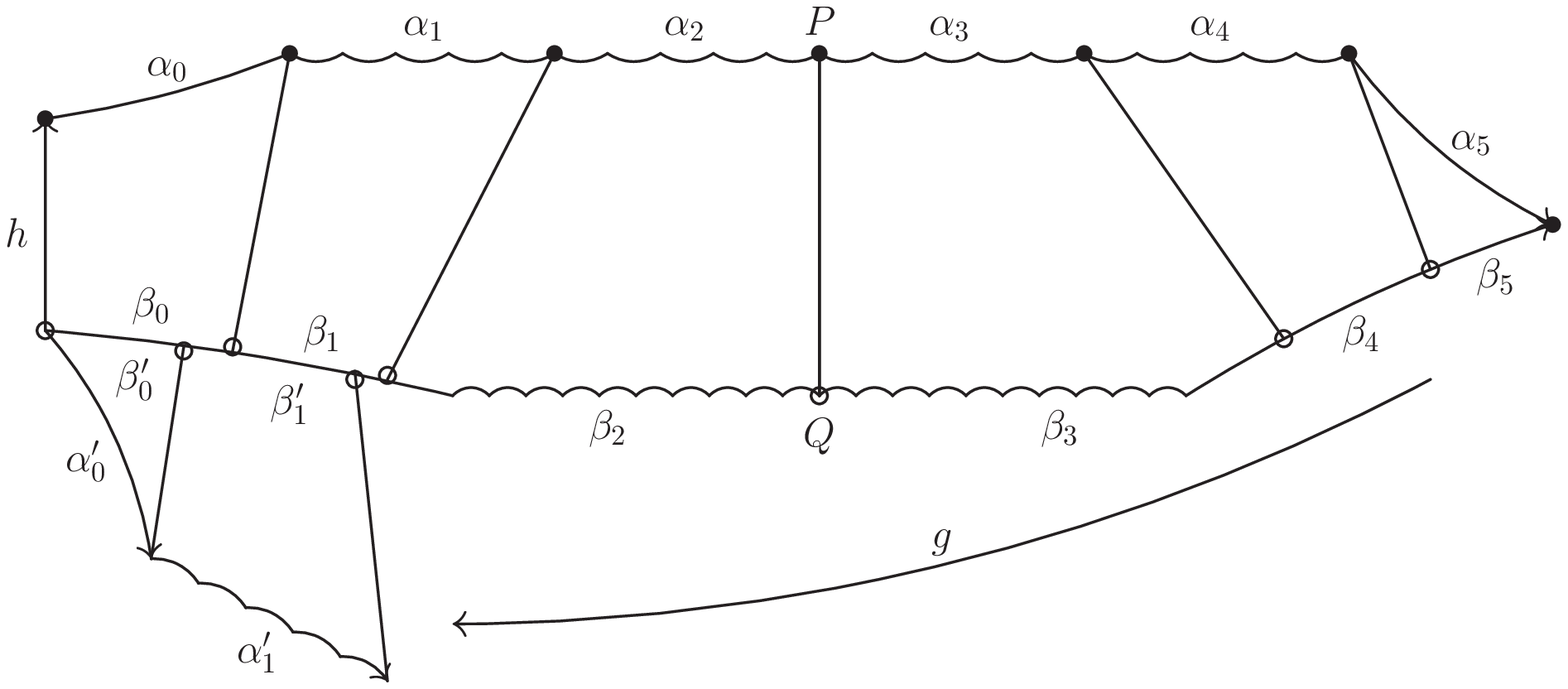}

\vspace*{-11cm}

\begin{center}
Fig. 9.
\end{center}

\medskip

{\bf Claim 4.} The Hausdorff distance between $\alpha_1$ and $\alpha_1'$ is at most $K_2$.

\medskip

{\it Proof.}
Since $\beta_i$ and $\alpha_i$ are $K_1$-similar for all $i$, we have
$$
\begin{array}{l}
d\bigl((\beta_j)_{-},(\beta_j)_{+}\bigr) \approx_{2K_1} d\bigl((\alpha_j)_{-},(\alpha_j)_{+}\bigr),\vspace*{2mm}\\
d\bigl((\beta_j')_{-},(\beta_j')_{+}\bigr) \approx_{2K_1} d\bigl((\alpha_j')_{-},(\alpha_j')_{+}\bigr)
\end{array}\eqno{(5.15)}
$$
for $j=0,1$. Observe that $\alpha_0$ and $\alpha_0'$ have the same labels
and $\alpha_1$ and $\alpha_1'$ have mutually inverse labels.
Therefore the numbers on the right sides of equations (5.15) are equal. This implies
$$
d\bigl((\beta_j)_{-},(\beta_j)_{+}\bigr)\approx_{4K_1} d\bigl((\beta_j')_{-},(\beta_j')_{+}\bigr).
$$
Since $\beta_0\beta_1$ and $\beta_0'\beta_1'$ are initial segments of the geodesic $\gamma_1$,
we deduce that

$$d((\beta_0)_{+},(\beta_0')_{+})\leqslant 4K_1,$$
$$d((\beta_1)_{+},(\beta_1')_{+})\leqslant 8K_1.$$
From here, from the $K_1$-similarity of $\alpha_i$ and $\beta_i$, and from the $K_1$-similarity of $\alpha_i'$ and $\beta_i'$ we deduce
$$
\begin{array}{ll}
d((\alpha_1)_{-},(\alpha_1')_{-}) & =
d((\alpha_0)_{+},(\alpha_0')_{+})\vspace*{2mm}\\
& \leqslant d((\alpha_0)_{+},(\beta_0)_{+})+d((\beta_0)_{+},(\beta_0')_{+})
+d((\beta_0')_{+},(\alpha_0')_{+})\vspace*{2mm}\\
& \leqslant K_1+4K_1+K_1=6K_1
\end{array}
$$
and
$$
\begin{array}{ll}
d((\alpha_1)_{+},(\alpha_1')_{+})
& \leqslant d((\alpha_1)_{+},(\beta_1)_{+})+d((\beta_1)_{+},(\beta_1')_{+})
+d((\beta_1')_{+},(\alpha_1')_{+})\vspace*{2mm}\\
& \leqslant K_1+8K_1+K_1=10K_1.
\end{array}
$$

By Lemma~\ref{close_qg}, the Hausdorff distance between $\alpha_1$ and $\alpha_1'$ is at most $K_2$.\hfill $\Box$

\medskip

Observe that $\alpha_1$ and $\alpha_1'$ satisfy assumptions of Theorem~\ref{acylindric}.
Indeed, $\alpha_1$ and $\alpha_1'$ are subpaths of the quasi-geodesics $L((\alpha_1)_{-},c_1)$
and $L((\alpha_1')_{-},c_1^{-1})$, respectively, each of them contains at least $f(K_2)$ periods by Claim 1,
and the Hausdorff distance between $\alpha_1$ and $\alpha_1'$ is at most $K_2$ by Claim 4.
By this theorem, there exist integers $s,t\neq 0$ such that
$z^{-1}c_1^{s}z=c_1^t$, where $z\in G$ is the element representing the label of any path from
$(\alpha_1')_{-}$ to $(\alpha_1)_{-}$. As such path we take $\overline{\alpha_0'}h \alpha_0$. Then
$z=uwu^{-1}$, and we have $w^{-1}c^sw=c^t$. Hence
$$w\in E_G(c).$$
Since $w=d^{-m}c^{-n}$, we have $d^m=c^{-n}w^{-1}\in E_G(c)$ and hence $E_G(c)=E_G(d)$. A contradiction to assumptions of lemma.

Thus, the inequality (5.2) is impossible, and we are done. \hfill $\Box$



\subsection{Loxodromic-elliptic case}

\begin{rmk}\label{four_delta_plus_1}
{\rm Suppose that $G$ is a group and $X\subseteq G$ is a (possibly infinite) generating set of $G$ such that
$\Gamma(G,X)$ is $\delta$-hyperbolic for some $\delta\geqslant 0$.
A subgroup $H$ of $G$ is called {\it elliptic} (with respect to $X$) if all elements of $H$ are elliptic.

It is well known that any elliptic subgroup $H$ of $G$ can be conjugated into the ball of radius $4\delta+1$ and center 1
in $\Gamma(G,X)$.
The proof of this fact is given in~\cite{BG} for the case where $G$ is a hyperbolic group. It also works under the above assumptions.
An alternative proof is given  in~\cite[Corollary 6.7]{Osin_1}.
}
\end{rmk}

\begin{lem}\label{lem 3.6}
Let $G$ be a group and let $X$ be a generating set of $G$.
Suppose that the Cayley graph $\Gamma(G,X)$ is hyperbolic and acylindrical.
Then there exists a constant $N_1>0$ such that
for any loxodromic (with respect to $X$) element $c\in G$, for any elliptic element $e\in G\setminus E_G(c)$, and for any $n\in \mathbb{N}$ we have that
$$|c^ne|_X>\frac{n}{N_1}.$$
\end{lem}

\medskip

{\it Proof.} The proof is very similar to the proof of Lemma~\ref{lem 3.4}.
However some pieces are new and some are easier.
Therefore we decided to present a complete proof for clearness.
Let $\delta\geqslant 0$ be a constant such that $\Gamma(G,X)$
is $\delta$-hyperbolic.

We use the following constants:

{\small $\bullet$} $\varkappa\geqslant 1$, $\varepsilon_0\geqslant 0$ and $n_0\in \mathbb{N}$
are the constants from Lemma~\ref{blue};

{\small $\bullet$} $\varepsilon_1=\varepsilon_1(\delta)\geqslant 0$ is the constant from Lemma~\ref{small_conjug};

{\small $\bullet$} $\varepsilon_2:=\max\{\varepsilon_0,\varepsilon_1\}$;

{\small $\bullet$} $\mu=\mu(\delta,\varkappa,\varepsilon_2)$, see Lemma~\ref{Hausdorff};

{\small $\bullet$} $\mathcal{C}$ is the constant from Theorem~\ref{acylindric}.

\medskip

We show that the lemma is valid for


$$
N_1=\max\bigl\{n_0,\,\, \frac{400(\mu +\delta+1)}{{\bf inj}(G,X)}+2\mathcal{C}+1,\,\, \varkappa(\varepsilon_0+8\delta+4)\bigr\}.\eqno{(5.16)}
$$


Suppose to the contrary that there exist a loxodromic (with respect to $X$) element $c\in G$, an elliptic element $e\in G\setminus E_G(c)$, and
a number $n\in \mathbb{N}$ such that
$$
n\geqslant N_1k, ,\hspace*{2mm} {\rm where}\hspace*{2mm} k=|c^ne|_X.\eqno{(5.17)}
$$

\medskip

Clearly, $k\neq 0$.

\medskip

{\bf Claim 1.} There exist $x,y\in G$ such that $e=x^{-1}yx$ and the following holds:
\begin{enumerate}
\item[1)] $|y|_X\leqslant 8\delta+3$,

\vspace*{2mm}

\item[2)] any path $q_0q_1q_2$ in $\Gamma(G,X)$, where $q_0$, $q_1$, $q_2$ are geodesics with labels
representing $x^{-1}$, $y$,~$x$, is a $(1,\varepsilon_1)$-quasi-geodesic path.
\end{enumerate}


{\it Proof.} It follows from (5.16) and (5.17) that $n\geqslant N_1k\geqslant N_1\geqslant n_0$. Then, by Lemma~\ref{blue}, we have
$$
|c^n|_X\geqslant \frac{1}{\varkappa}n-\varepsilon_0.\eqno{(5.18)}
$$
From this we deduce

$$
\begin{array}{ll}
k=|c^ne|_X & \geqslant |c^n|_X-|e|_X\, \,\overset{(5.18)}{\geqslant}\,\, \frac{1}{\varkappa}n-\varepsilon_0-|e|_X
\,\, \geqslant\frac{1}{\varkappa}N_1k-\varepsilon_0-|e|_X \vspace*{3mm}\\
& \overset{(5.16)}{\geqslant}
\,\, k(\varepsilon_0+8\delta+4)-\varepsilon_0-|e|_X\,\, \geqslant
\,\, k+8\delta+3-|e|_X.
\end{array}
$$
Hence $|e|_X\geqslant 8\delta+3$. Since $e$ is elliptic, $e$ is conjugate to an element of $G$
of length at most $4\delta+1$ (see Remark~\ref{four_delta_plus_1}).
Thus, the assumption and hence the conclusion of Lemma~\ref{small_conjug} are satisfied.
\hfill $\Box$

\medskip

For any $g\in G$, let $S(g)$ be the set of shortest elements in the conjugacy class
of $g$.
Let $u\in G$ be a shortest element for which there exists $c_1\in S(c)$ with the property $c=u^{-1}c_1u$.




We denote $w=(c^ne)^{-1}$. Then $|w|_X=k$ and we have the equation
$$
u^{-1}\underbrace{c_1c_1\dots c_1}_nux^{-1}yxw=1.\eqno{(5.19)}
$$

\medskip

\noindent
Consider a geodesic $(n+6)$-gon $\mathcal{P}=p_0(p_1p_2\dots p_n)p_{n+1}\bar{q}_{2}\bar{q}_1\bar{q}_0h$
in the Cayley graph $\Gamma(G,X)$ such that the labels of its sides correspond
to the elements in the left side of equation (5.19), see Fig.~10.






\medskip

\vspace*{-35mm}
\hspace*{-9.5mm}
\includegraphics[scale=0.8]{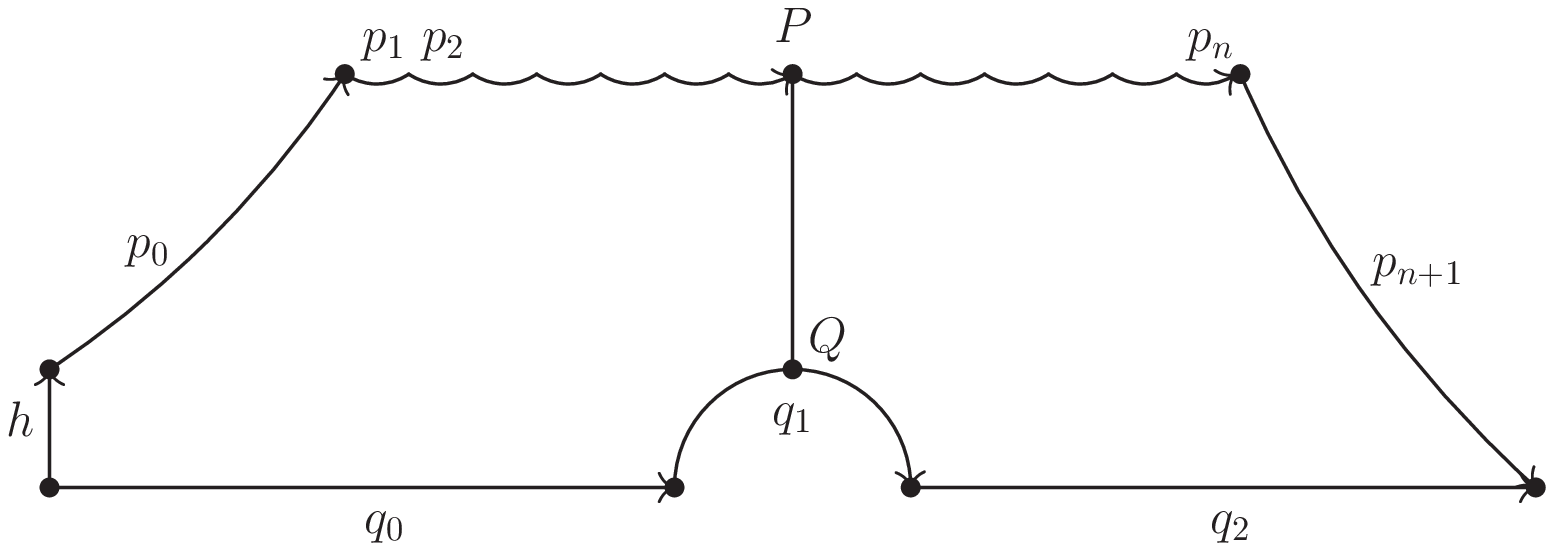}

\vspace*{-15.5cm}

\begin{center}
Fig. 10.
\end{center}

\bigskip

In particular, the labels of the paths $q_0,q_1, q_2$ correspond to the elements
$x^{-1},y^{-1},x$. By Claim 1, the path $q:=q_0q_1q_2$ is a $(1,\varepsilon_1)$-quasi-geodesic.
Since $n\geqslant N_1k\geqslant N_1\geqslant n_0$, we have by  Lemma~\ref{blue}
that the path $p:=p_0p_1p_2\dots p_np_{n+1}$ is a $(\varkappa,\varepsilon_0)$-quasi-geodesic.
Then $p$ and $q$ are $(\varkappa,\varepsilon_2)$-quasi-geodesics.
In the following claims we use the following constants.
$$
K=k+36\mu+26\delta,\hspace*{5mm}
K_1=K+2\mu+6\delta+2,\hspace*{5mm}
K_2=K_1+2\mu+2\delta.
$$



\medskip

{\bf Claim 2.} The quasi-geodesic $p_1p_2\dots p_n$ contains $n\geqslant 2f(2K_2)+1$\,\,\break $c_1$-periods.

\medskip

{\it Proof.} The claim follows straightforward from the definition of function
$f(r)$ in Theorem~\ref{acylindric} and from (5.16) and (5.17).\hfill $\Box$

\medskip

Let $P$ be the middle point of the quasi-geodesic $p_1p_2\dots p_n$ and let
$Q$ be the middle point of the geodesic $q_1$.
As in Claim~2 of the proof of Lemma~\ref{lem 3.4}, we have
$$
d(P,Q)\leqslant K.
$$


\medskip

\vspace*{-35mm}
\hspace*{-9.5mm}
\includegraphics[scale=0.8]{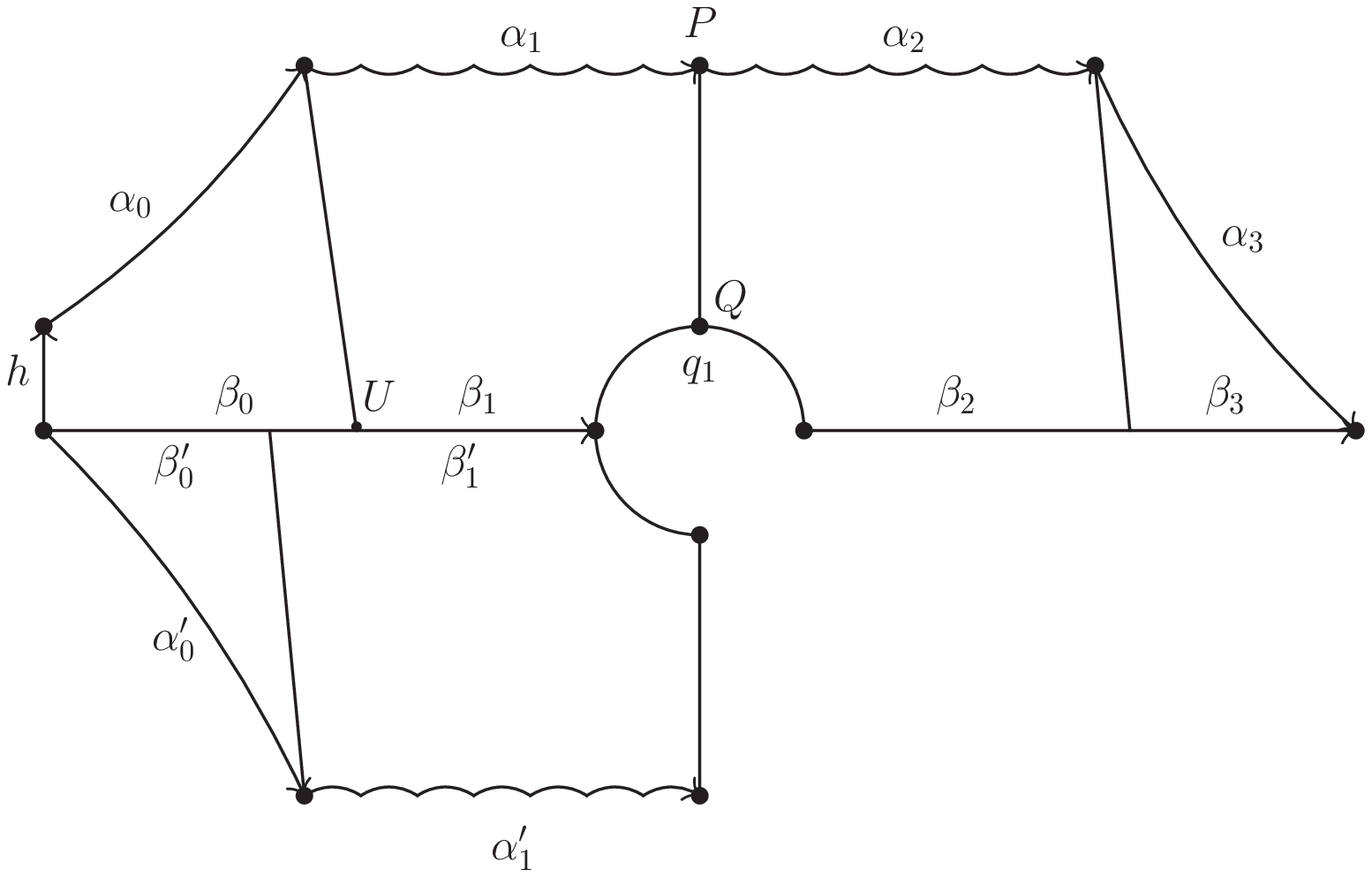}

\vspace*{-11.5cm}

\begin{center}
Fig. 11.
\end{center}

\bigskip

\noindent
We consider the decomposition $p_0p_1p_2\dots p_np_{n+1}=\alpha_0\alpha_1\alpha_2\alpha_3$, where
$\alpha_0=p_0$, $\alpha_3=p_{n+1}$, and $\alpha_1$ and $\alpha_2$ are determined by the condition
$(\alpha_1)_{+}=P=(\alpha_2)_{-}$, see Fig. 11.

\medskip

{\bf Claim 3.} There are decompositions $q_0=\beta_0\beta_1$ and $q_2=\beta_2\beta_3$ such that $\alpha_i$ and
$\beta_i$ are $K_1$-similar for $i=0,\dots ,3$. In particular, the Hausdorff distance between $\alpha_i$ and $\beta_i$ is at most $K_2$.

\medskip

{\it Proof.} Because of symmetry, we show only that the first decomposition exists.
By Claim 1, we have $d((q_1)_{-},(q_1)_{+})=|y|_X\leqslant 8\delta+3$. Hence
$$
\begin{array}{ll}
d((\alpha_0\alpha_1)_{-},(q_0)_{-}) & = k,\vspace*{2mm}\\
d((\alpha_0\alpha_1)_{+},(q_0)_{+}) & \leqslant d(P,Q)+\frac{1}{2}d((q_1)_{-},(q_1)_{+})< K+4\delta+2.
\end{array}
$$
By Lemma~\ref{close_qg} applied to the $(\varkappa,\varepsilon_2)$-quasi-geodesics $\alpha_0\alpha_1$ and $q_0$,
there exists a point $U\in q_0$ such that $d((\alpha_1)_{+},U)\leqslant K+4\delta+2+(2\mu+2\delta)=K_1$.

The point $U$ divides the path $q$ into two subpaths. We denote them $\beta_0,\beta_1$. By construction,
the paths $\alpha_i$ and $\beta_i$ are $K_1$-similar. Again by Lemma~\ref{close_qg}, the Hausdorff distance between $\alpha_i$ and $\beta_i$ is at most $K_1+2\mu+2\delta= K_2$.
\hfill $\Box$

\medskip

Recall that the geodesics $q_0$ and $q_2$ have mutually inverse labels.
Therefore, there exists $g\in G$ such that $q_0=g(\overline{q_2})$. We denote
$$
\alpha_0'=g(\overline{\alpha_3}),\hspace*{3mm} \alpha_1'=g(\overline{\alpha_2}),\hspace*{3mm} \beta_0'=g(\overline{\beta_3}),\hspace*{3mm} \beta_1'=g(\overline{\beta_2}).
$$

Then $\beta_0\beta_1=\beta_0'\beta_1'=q_0$.

\medskip

{\bf Claim 4.} One of $\alpha_1$, $\alpha_1'$ lies in the $2K_2$-neighborhood of the other.

\medskip

{\it Proof.} By Claim 3, the Hausdorff distance between $\alpha_1$ and $\beta_1$
is at most~$K_2$. Also the Hausdorff-distance between $\alpha_1'$ and $\beta_1'$ is at most $K_2$.
The claim follows from the fact that one of $\beta_1$, $\beta_1'$ is a subsegment of the other.
\hfill $\Box$

\medskip

Observe that $\alpha_1$ and $\alpha_1'$ satisfy assumptions of Theorem~\ref{acylindric}.
Indeed, $\alpha_1$ and $\alpha_1'$ are subpaths of the quasi-geodesics $L((\alpha_1)_{-},c_1)$
and $L((\alpha_1')_{-},c_1^{-1})$, respectively, each of them contains at least $f(2K_2)$ periods by Claim 2, and one of them lies in the $2K_2$-neighborhood of the other by Claim 3.
By this theorem, there exist integers $s,t\neq 0$ such that
$z^{-1}c_1^{s}z=c_1^t$, where $z\in G$ is the element representing the label of any path from
$(\alpha_1')_{-}$ to $(\alpha_1)_{-}$. As such path we take $\overline{\alpha_0'}h\alpha_0$.
Then $z=uwu^{-1}$, and we have $w^{-1}c^sw=c^t$.
Hence
$$w\in E_G(c).$$
Since $w=(c^ne)^{-1}$, we have $e\in E_G(c)$. A contradiction to assumptions of the lemma.
Thus, the inequality (5.17) is impossible, and we are done. \hfill $\Box$

\section{Extension of quasi-morphisms from hyperbolically embedded subgroups to the whole group}

\subsection{A sufficient condition for preserving the elipticity under decreasing of a
generating set}

\bigskip


\begin{lem}\label{elliptic-lox}
Let $G$ be a group and $X$ be a generating set of $G$.
Suppose that $\Gamma(G,X)$ is hyperbolic and acylindrical.
Let $a_1,\dots ,a_k\in G$ be finitely many loxodromic elements with respect to $X$.
Let $g\in G$ be an element non-commensurable with elements
of $A=\{a_1,\dots,a_k\}$. If $g$ is elliptic with respect to
$X'=X\cup E_G(a_1)\cup \dots \cup E_G(a_k)$, then $g$ is elliptic with respect to
$X$.

\end{lem}

\medskip

{\it Proof.} 
Let $B_i$ be a finite set of representatives of left cosets of $\langle a_i\rangle$
in $E_G(a_i)$, $i=1,\dots,k$. Enlarging $X$ by a finite set does not change the property of $\Gamma(G,X)$ to be hyperbolic and acylindrical (see~\cite[Lemma~5.1]{Osin_1}) and the property of elements of $G$ to be elliptic or loxodromic. Therefore,
we may assume that $X$ contains $\overset{k}{\underset{i=1}{\cup}} B_i$. Let $\delta\geqslant 0$ be a constant such that $\Gamma(G,X)$
is $\delta$-hyperbolic.



To the contrary, suppose that $g$ is elliptic with respect to
$X'$ and not elliptic with respect to $X$.

Then the following holds.

1) There exists $R>0$ such that $|g^n|_{X'}\leqslant R$ for all $n\in \mathbb{Z}$.

2) $g$ is loxodromic with respect to $X$.

\medskip

We fix an arbitrary $n\in \mathbb{N}$ and write $g^n=x_1x_2\dots x_t$, where $x_i\in X'$ and $t$ is minimal.
In particular, $t\leqslant R$.
Each element of $E_G(a_i)$ can be written in the form $ba_i^k$ for some $b\in B_i\subset X$. Therefore
$g^n$ can be written as $$g^n=u_0a_{i_1}^{k_1}u_1a_{i_2}^{k_2}\dots u_{s-1}a_{i_s}^{k_s}u_s,$$  where $u_0,\dots ,u_{s}$ are words in $X$,  $a_{i_1},a_{i_2}\dots , a_{i_s}\in A$,  $0\leqslant s\leqslant t\leqslant R$, and
$$\overset{s}{\underset{j=0}{\sum}}|u_j|_X\leqslant t\leqslant R.\eqno{(6.1)}$$

Let $f_0,h\in G$ be elements such that $g=f_0^{-1}hf_0$ and $h$ is a shortest element (with respect to $X$) in the conjugacy class of $g$.
For each $a_i\in A$, let $f_i,b_i\in G$ be elements such that $a_i=f_i^{-1}b_if_i$ and $b_i$ is a shortest element (with respect to $X$) in the conjugacy class of $a_i$. Let $F=\max \{|f_i|_X\,|\, i=0,\dots,k\}$.
Then $$h^n=v_0b_{i_1}^{k_1}v_1b_{i_2}^{k_2}\dots v_{s-1}b_{i_s}^{k_s}v_s,$$ where $v_0=f_0u_0f_{i_1}^{-1}$,
$v_j=f_{i_j}u_jf_{i_{j+1}}^{-1}$, $j=1,\dots ,s-1$, and $v_s=f_{i_s}u_sf_0^{-1}$. Using (6.1),
we have
$$\overset{s}{\underset{j=0}{\sum}}|v_j|_X\leqslant
\overset{s}{\underset{j=0}{\sum}}|u_j|_X+2|f_0|_X+2\overset{s}{\underset{j=1}{\sum}}|f_{i_j}|_X
\leqslant R+2(R+1)F.\eqno{(6.2)}$$

Let $\mathcal{P}$ be a geodesic $2(s+1)$-gon in $\Gamma(G,X)$ with sides $p_0,q_0,p_1,\dots,p_s,q_s$
such that $p_0,p_1,\dots,p_s$ are quasi-geodesics representing $h^{-n},b_{i_1}^{k_1},\dots ,b_{i_s}^{k_s}$
and $q_0,\dots,q_s$ are geodesics representing $v_0,\dots ,v_s$.

\vspace*{-30mm}
\hspace*{5mm}
\includegraphics[scale=0.8]{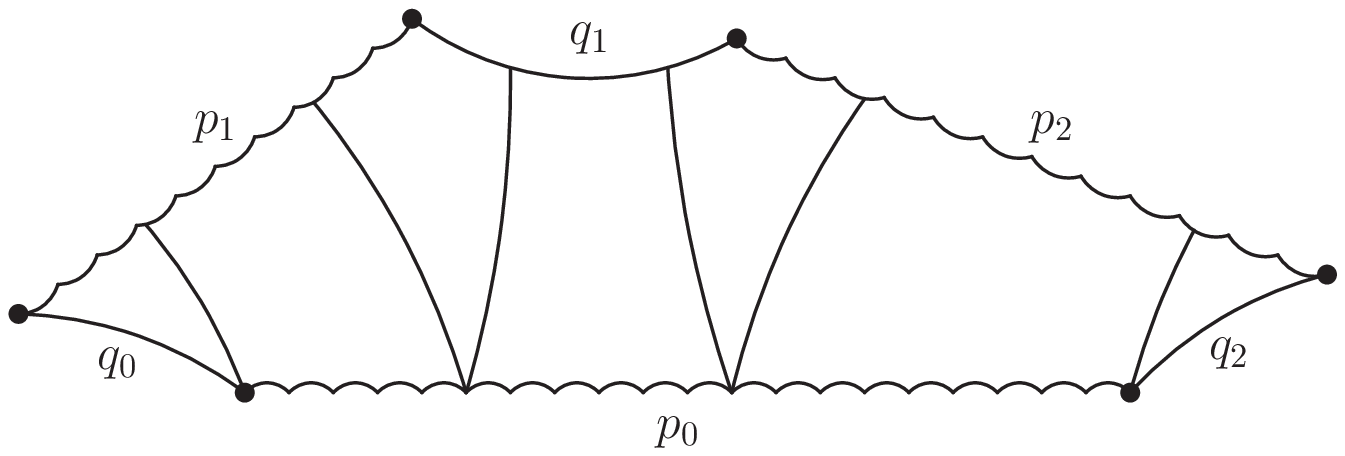}

\vspace*{-15.5cm}

\begin{center}
Fig. 12. Case $s=2$.
\end{center}

\bigskip

Since the elements $h,b_{i_1},\dots,b_{i_s}$ are loxodromic with respect to $X$ and have minimal length in their conjugacy classes, the paths $p_0,p_1,\dots ,p_s$ are $(\varkappa,\varepsilon)$-quasi-geodesics, where
$\varkappa$ and $\varepsilon$ are universal constants from Corollary~\ref{qg}.
We set
$$\alpha=2R\delta+2\mu(\delta,\varkappa,\varepsilon),$$
where $\mu(\delta,\varkappa,\varepsilon)$ is the constant from Lemma~\ref{Hausdorff}.

{\bf Claim 1.} The side $p_0$ lies in the $\alpha$-neighborhood of the union of other sides of $\mathcal{P}$.

\medskip

{\it Proof.}
For each $i=0,\dots,s$, we chose a geodesic segment $\widetilde{p_i}$ such that $(\widetilde{p_i})_{-}=(p_i)_{-}$ and $(\widetilde{p_i})_{+}=(p_i)_{+}$. Consider the geodesic $2(s+1)$-gon  $\widetilde{\mathcal{P}}$ with the sides $\widetilde{p_0},q_0,\widetilde{p_1},\dots,\widetilde{p_s},q_s$.
The side $\widetilde{p_0}$ lies in the $2s\delta$-neighborhood of the other sides of $\widetilde{\mathcal{P}}$.
By Lemma~\ref{Hausdorff}, the Hausdorff distance between $\widetilde{p_i}$ and $p_i$ is at most
$\mu(\delta,\varkappa,\varepsilon)$ for every $i$. This completes the proof. \hfill $\Box$

\medskip

For $i=0,\dots,s$, we set
$$D_i=\varkappa(2\alpha+|v_i|_X+\varepsilon)+1.$$

\medskip

{\bf Claim 2.} For any $i\in \{0,\dots ,s\}$, the $\alpha$-neighborhood of $q_i$ contains at most $D_i$ points of $p_0$.


\medskip

{\it Proof.} Let $Q_i$ be the set of points on $p_0$ which lie in the $\alpha$-neighborhood of $q_i$.
Suppose that $Q_i\neq \emptyset$ and let $z_1$ and $z_2$ be the leftmost and the rightmost points of $Q_i$
on $p_0$.
Then $d_X(z_1,z_2)\leqslant 2\alpha + \ell(q_i)=2\alpha+|v_i|_X$. Therefore the length of the subpath of $p_0$
connecting $z_1$ and $z_2$ is at most $\varkappa(2\alpha+|v_i|_X+\varepsilon)$, and  the claim follows.
\hfill $\Box$

\medskip

We set
$$
\beta=\eta(\delta,\varkappa,\varepsilon,\alpha),
$$
where $\eta$ is the function from Lemma~\ref{close_qg}.
Also, for $j=1,\dots,s$, we set
$$C_j:=\max\{\varkappa\bigl((f(\beta)+2)|b_{i_j}|_X+2\alpha+\varepsilon\bigr), (f(\beta)+2)|h|_X\},\eqno{(6.3)}$$
where $f$ is the function from Theorem~\ref{acylindric}.

\medskip

{\bf Claim 3.}
For any $j\in \{1,\dots ,s\}$, the $\alpha$-neighborhood of $p_j$ contains at most $C_j$ points of $p_0$.

\medskip

{\it Proof.}
Let $c$ be the maximal subpath of $p_0$ for which there exists a subpath $c'$ of  $p_j$ or
$\overline{p_j}$ with the property
$$d_X(c_{-},c_{-}')\leqslant \alpha \hspace*{2mm} {\text{\rm and}}\hspace*{2mm} d_X(c_{+},c_{+}')\leqslant \alpha.\eqno{(6.4)}$$
Suppose that the $\alpha$-neighborhood of $p_j$ contains more than $C_j$ points of $p_0$. Then $$\ell(c)\geqslant C_j.$$
We check the following statements:
\begin{enumerate}
\item[(1)] The Hausdorff distance between $c$ is $c'$ is at most $\beta$.

\item[(2)] The path $c$ contains at least $f(\beta)$ $h$-periods.

\item[(3)] The path $c'$ contains at least $f(\beta)$ $b_{i_j}^{\pm}$-periods.
\end{enumerate}

\medskip

\noindent
Statement (1) follows from (6.4) and Lemma~\ref{close_qg}.
Statement (2) follows from $\ell(c)\geqslant C_j\overset{(6.3)}{\geqslant} (f(\beta)+2)|h|_X.$
Finally, statement (3) follows from
$$
\begin{array}{ll}
\ell(c')\geqslant d_X(c_{-}',c_{+}') & \geqslant d_X(c_{-},c_{+})
-d_X(c_{-},c_{-}')-d_X(c_{+},c_{+}') \geqslant d_X(c_{-},c_{+})-2\alpha\vspace*{2mm}\\
& \geqslant \frac{1}{\varkappa}\ell(c)-\varepsilon-2\alpha \geqslant \frac{1}{\varkappa}C_j-\varepsilon-2\alpha \overset{(6.3)}{\geqslant} (f(\beta)+2)|b_{i_j}|_X.
\end{array}
$$
It follows from statements (1)-(3) and Theorem~\ref{acylindric} that $h$ and $b_{i_j}$
are commensurable. Then $g$ and $a_{i_j}$ are commensurable. A contradiction. \hfill $\Box$

\medskip

It follows from Claims 1-3 that $p_0$ contains at most
$\overset{s}{\underset{i=0}{\sum}} D_i+ \overset{s}{\underset{j=1}{\sum}} C_j$ points.
Using (6.2), we have
$$
\overset{s}{\underset{i=0}{\sum}} D_i\leqslant
(R+1)(\varkappa(2\alpha+\varepsilon)+1)+\varkappa (R+2(R+1)F).
$$
We also have
$$
\overset{s}{\underset{j=1}{\sum}} C_j\leqslant R\max C_j\leqslant
R\max\{\varkappa\bigl((f(\beta)+2)(\underset{i=1,\dots,k}{\max}|a_i|_X)+2\alpha+\varepsilon\bigr), (f(\beta)+2)|h|_X\}.
$$

Therefore $\ell_X(p_0)$ is bounded from above by a constant which does not depend on~$n$.
This contradicts the fact that $\ell_X(p_0)=n|h|_X\geqslant n$. \hfill $\Box$


\medskip

\subsection{Improving relative generating sets for hyperbolically embedded subgroups}

Recall that in the situation
$\{E_1,\dots ,E_k\}\hookrightarrow_h G$, we use notation $\mathcal{E}=E_1\sqcup \dots \sqcup E_k$.
The following lemma follows directly from Definition~\ref{def_hyperb_embedd}.

\begin{lem}\label{change_rel_gener} Suppose that $\{E_1,\dots ,E_k\}\hookrightarrow_h (G,X)$ and let
$Y$ be a subset of $G$ such that $X\subseteq Y\subseteq \langle X\rangle$ and $\underset{y\in Y}{\sup}|y|_X<\infty$.
Then $\{E_1,\dots ,E_k\}\hookrightarrow_h (G,Y)$.
\end{lem}

\medskip

{\it Proof.} The hyperbolicity of $\Gamma(G,Y\sqcup \mathcal{E})$ follows from the
hyperbolicity of $\Gamma(G,X\sqcup \mathcal{E})$ by Lemma~\ref{sup}.
The local finiteness of $(E_i,\widehat{d_i}^Y)$, where the relative metric $\widehat{d_i}^Y$
on $E_i$ is defined using $Y$, follows from the local finiteness of $(E_i,\widehat{d_i}^X)$,
where the relative metric $\widehat{d_i}^X$ on $E_i$ is defined using $X$.
\hfill $\Box$

\medskip

\begin{lem}\label{best_Y}
Suppose that $\{E_1,\dots ,E_k\}\hookrightarrow_h (G,X)$, where $E_1,\dots,E_k$ are infinite virtually cyclic
subgroups of $G$.
Then there exists a subset $Y\subseteq G$ containing $X$ such that the following properties are satisfied:

\begin{enumerate}

\item[(1)] $\langle Y\rangle=G$.

\item[(2)] $\{E_1,\dots ,E_k\}\hookrightarrow_h (G,Y)$.

\item[(3)] $\Gamma(G,Y)$ and $\Gamma(G,Y\sqcup \mathcal{E})$ are hyperbolic and acylindrical.

\item[(4)] If $G$ is not virtually cyclic, then $G$ is acylindrically hyperbolic with respect to $Y$.

\item[(5)] If $g\in G$ is an elliptic element with respect to $Y\cup \overset{k}{\underset{i=1}{\cup}} E_i$ such that $g$
 is non-commen\-sur\-able with elements of $\overset{k}{\underset{i=1}{\cup}} E_i$ of infinite order, then $g$ is elliptic with respect~to~$Y$.
    Moreover, there exists $u\in G$ such that $u^{-1}g^mu\in Y$ for all $m\in \mathbb{Z}$.
\end{enumerate}
\end{lem}

\medskip

{\it Proof.} By Theorem~\ref{enlarging},
there exists a subset $X_1\subseteq G$ such that $X\subseteq X_1$
and the following conditions hold.

\begin{enumerate}
\item[(i)] $\{E_1,\dots ,E_k\}\hookrightarrow_h (G,X_1)$.
\item[(ii)] $\Gamma(G,X_1\sqcup \mathcal{E})$ is hyperbolic and acylindrical.
\end{enumerate}

\medskip

We show that $X_1$ can be chosen so that, additionally, $X_1$ generates $G$.
Let $A_i$ be a finite generating set of $E_i$. We set $X_2:=X_1\cup A_1\cup \dots \cup A_k$.
Since $G=\langle X_1\cup \overset{k}{\underset{i=1}{\cup}} E_i\rangle$, we have $G=\langle X_2\rangle$.

\medskip

{\bf Claim 1.}

\begin{enumerate}
\item[(a)] $\{E_1,\dots ,E_k\}\hookrightarrow_h (G,X_2)$.
\item[(b)] $\Gamma(G,X_2\sqcup \mathcal{E})$ is hyperbolic and acylindrical.

\item[(c)] $\Gamma(G,X_2)$ is hyperbolic and acylindrical.
\end{enumerate}

\medskip

{\it Proof.} Since $|X_1\Delta X_2|<\infty$, we have  ${\rm (i)}\Leftrightarrow {\rm (a)}$ by Lemma~\ref{hyp_embed_1}
and ${\rm (ii)}\Leftrightarrow {\rm (b)}$ by Lemma~\ref{sup}.
To prove (c), we first observe that $\{1\}\hookrightarrow_h (E_i,A_i)$ for $i=1,\dots,k$, and recall that $\{E_1,\dots ,E_k\}\hookrightarrow_h (G,X_1)$.
By~\cite[Proposition 4.35]{DOG}, this implies that $\{1\}\hookrightarrow_h (G,X_2)$.
In particular, by definition this means that $\Gamma(G,X_2)$ is hyperbolic.
The acylindricity of $\Gamma(G,X_2)$ can be proved as in the part of the proof of~\cite[Theorem~1.4]{Osin_1},
starting from the words ``Let us show that $\Gamma(G,X)$ is acylindrical''.\hfill $\Box$

\medskip

Let $\delta>0$ be a number such that  $\Gamma(G,X_2)$ is $\delta$-hyperbolic.
We set
$$Y=\{g\in G\,|\, |g|_{X_2}\leqslant 4\delta+1\}.$$ Since $G=\langle X_2\rangle$, we have $G=\langle Y\rangle$, i.e. (1).
Statement (2) follows from Lemma~\ref{change_rel_gener} and Claim 1 (a).






\medskip

Now we prove (3). The hyperbolicity and acylindricity of $\Gamma(G,Y\sqcup \mathcal{E})$ follows from
the hyperbolicity and acylindricity of $\Gamma(G,X_2\sqcup \mathcal{E})$ by~Lemma~\ref{sup}.
Analogously, the hyperbolicity and acylindricity of $\Gamma(G,Y)$ follows from the hyperbolicity and acylindricity of $\Gamma(G,X_2)$.

For (4) and (5), we first prove the following claim.

\medskip

{\bf Claim 2.} Let $a_i\in E_i$ be an arbitrary element of infinite order.
Then $a_i$ is loxodromic with respect to $Y$ and $E_i=E_G(a_i)$.

\medskip

{\it Proof.} By statement (3), $\Gamma(G,Y)$ is hyperbolic and acylindrical. Therefore any element of $G$ is either elliptic or loxodromic with respect to $Y$. By statement (2), the space $(E_i,\widehat{d_i}^{Y})$ is locally finite, hence $a_i$ cannot be elliptic with respect to $Y$.

Since $E_G(a_i)$ is the maximal virtually cyclic subgroup containing $a_i$, we have $E_i\subseteq E_G(a_i)$.
The inverse inclusion follows from Lemma~\ref{periph_subgr} and the algebraic characterisation of $E_G(a_i)$ in (3.1).
\hfill $\Box$

\medskip

We prove (4). Suppose that $G$ is not virtually cyclic. Since $\Gamma(G,Y)$ is hyperbolic and acylindrical, it suffices to show that the action of $G$ on $\Gamma(G,Y)$ is non-elementary. By Claim 2, this action has unbounded orbits. Thus, cases (a) and (b) of Theorem~\ref{check_elementary} are excluded. The remaining case (c) of this theorem says that $G$ is acylindrically hyperbolic with respect to $Y$.

\medskip

Finally we prove (5).
Suppose that $g\in G$ is an elliptic element with respect to $Y\cup \overset{k}{\underset{i=1}{\cup}} E_i$ and that $g$
is non-commensurable with elements of $\overset{k}{\underset{i=1}{\cup}} E_i$ of infinite order.
Assumptions of Lemma~\ref{elliptic-lox} are valid for $\Gamma(G,Y)$, the elements $a_1,\dots, a_k$ from Claim~2, and the element $g$.
By this lemma, $g$ is elliptic with respect to $Y$.

Since $\underset{y\in Y}{\sup}|y|_{X_2}\leqslant 4\delta+1$, we conclude that $g$ is elliptic with respect to $X_2$.
By Remark~\ref{four_delta_plus_1}, $\langle g\rangle$ is conjugated into $Y$.\hfill $\Box$


\medskip

\subsection{Quasi-morphisms}


Let $G$ be a group.
Recall that a map $q:G\rightarrow \mathbb{R}$ is called a {\it quasi-morphism} if there exists a constant
$\varepsilon>0$ such that for every $f,g\in G$ we have
$$
|q(fg)-q(f)-q(g)|<\varepsilon.
$$

For a quasi-morphism $q:G\rightarrow \mathbb{R}$, we define its {\it defect} $D(q)$ by
$$
D(q)=\underset{f,g\in G}{\sup}|q(fg)-q(f)-q(g)|.
$$

The quasi-morphism $q$ is called {\it homogeneous} if $q(g^m)=m\cdot q(g)$ for all $g\in G$ and $m\in \mathbb{Z}$.

\begin{rmk}
{\rm For every quasi-morphism $q:G\rightarrow \mathbb{R}$, there exists  a unique {\it homogenous} quasi-morphism  $\widetilde{q}:G\rightarrow \mathbb{R}$ which lies at a bounded distance from $q$. This means that the following conditions are satisfies.

\begin{enumerate}
\item[(1)] $\widetilde{q}(g^m)=m\cdot \widetilde{q}(g)$ for all $g\in G$ and $m\in \mathbb{Z}$.

\item[(2)] There exists $C\geqslant 0$ such that $|q(g)-\widetilde{q}(g)|\leqslant C$ for all $g\in G$.
\end{enumerate}
\medskip

The quasi-morphism $\widetilde{q}$ is defined by the formula $$\widetilde{q}(g)=\underset{m\rightarrow \infty}{\lim}\frac{q(g^m)}{m}.\eqno{(6.5)}$$
}
\end{rmk}
\begin{lem}\label{lem 4.5}
Suppose that $q:G\rightarrow \mathbb{R}$ is a homogeneous quasi-morphism. Then $q$ is constant on each conjugacy class of elements of $G$. In particular, if $a,b\in G$ are two commensurable elements and $q(a)=0$, then $q(b)=0$.
\end{lem}

\medskip

{\it Proof.}
Let $u,g,h\in G$ be elements such that $u^{-1}gu=h$. Then
$$
|q(g^m)-q(h^m)|\leqslant 2(q(u)+D(q))
$$
for all $m\in \mathbb{N}$. Hence
$$q(g)=\underset{m\rightarrow \infty}{\lim}\frac{q(g^m)}{m}=
\underset{m\rightarrow \infty}{\lim}\frac{q(h^m)}{m}=q(h).$$
\hfill $\Box$

Statements (a) and (b) of the following corollary follow from~\cite[Theorem~4.2]{HO}).
We show that statement (c) can be deduced from the proof of this theorem combined with Lemma~\ref{best_Y}.

\begin{cor}\label{quasimorphisms}
Suppose that $G$ is a group and $E_1,\dots,E_k$ are infinite cyclic subgroups of $G$ generated by elements $a_1,\dots,a_k$,
respectively. Suppose that
$\{E_1,\dots ,E_k\}\hookrightarrow_h (G,X)$ and denote $E=\overset{k}{\underset{i=1}{\cup}} E_i$.
Then for all $I\subseteq \{1,\dots,k\}$, there exists a homogenous quasi-morphism $\widetilde{q}:G\rightarrow \mathbb{R}$
such that the following hold.

\begin{enumerate}
\item[{\rm (a)}] $\widetilde{q}(a_i)=1$ for all $i\in I$.
\vspace*{2mm}
\item[{\rm (b)}] $\widetilde{q}(a_i)=0$ for all $i\notin I$.
\vspace*{2mm}
\item[{\rm (c)}] $\widetilde{q}(g)=0$ for all $g\in {\text{\rm Ell}}(G, X\cup E)$ that are  non-commensurable
with elements of $E$.
\end{enumerate}
\end{cor}

\medskip

{\it Proof.} By Lemma~\ref{best_Y}, there exists a subset $Y\subseteq G$ such that
$X\subseteq Y$ and the statements (1)-(5) of this lemma are satisfied.
In particular, we have $\{E_1,\dots ,E_k\}\hookrightarrow_h (G,Y)$.
By~\cite[Theorem~4.2]{HO} applied to this hyperbolic embedding,
there exists a quasi-morphism $q:G\rightarrow \mathbb{R}$ (possibly inhomogenous) such that
the following hold.

(a$'$) $q(a_i^n)=n$ for all $i\in I$ and all $n\in \mathbb{Z}$.

(b$'$) $q(a_i^n)=0$ for all $i\notin I$ and all $n\in \mathbb{Z}$.

\medskip

Moreover, by construction of $q$ in the proof of this theorem, we obtain

(c$'$) $q(y)=0$ for all $y\in Y\setminus E$.

\medskip

Let $\widetilde{q}$ be the  homogenous quasi-morphism obtained from $q$ by formula (6.5). Then
$\widetilde{q}$ satisfies conditions (a), (b). We prove that $\widetilde{q}$ satisfies condition (c).

Condition (c) is obviously valid for elements $g\in G$ of finite order. Suppose that $g\in G$ has infinite order,
is elliptic with respect to $X\cup E$ and non-commensurable with elements of $E$. Then $g$ is elliptic with respect to $Y\cup E$.
By statement (5) of Lemma~\ref{best_Y}, there exists $u\in G$ such that $u^{-1}g^mu\in Y$ for all $m\in \mathbb{N}$. Since $g$ is non-commensurable with elements of~$E$, we have $u^{-1}g^mu\in Y\setminus E$ for all $m\in \mathbb{N}$. By (c$'$), we obtain
$
q(u^{-1}g^mu)=0
$
for all $m\in \mathbb{N}$.
Then $\widetilde{q}(u^{-1}gu)=0$.
By Lemma~\ref{lem 4.5}, we have $\widetilde{q}(g)=0$, i.e. (c).
\hfill $\Box$

\section{Equation $x^ny^m=a^nb^m$ in acylindrically hyperbolic groups}

\begin{prop}\label{lem 0.1}
Let $G$ be an acylindrically hyperbolic group with respect to a generating set $Z$.
Suppose that $a$ and $b$ are two non-commensurable special (with respect to $Z$)
elements of $G$.
Then there exists a generating set~$Y$
containing $\mathcal{E}=\langle a\rangle\cup \langle b\rangle$ and there exists a number $N\in \mathbb{N}$ such that for all $n,m>N$ the following holds:

If $(c,d)$ is a solution of the equation $x^ny^m=a^nb^m$, then one of the following holds:

\begin{enumerate}
   \item[1)] $c$ and $d$ are loxodromic with respect to $Y$, and
$E_G(d)=E_G(c)$;

\item[2)] $c$ is loxodromic with respect to $Y$ and $d$ is elliptic, and $d^m\in E_G(c)$;

\item[3)] $d$ is loxodromic with respect to $Y$ and $c$ is elliptic, and
$c^n\in E_G(d)$;

\item[4)] $c$ and $d$ are elliptic with respect to $Y$, and one of the following holds:

\begin{enumerate}
\item[(a)]
$c$ is conjugate to $a$ and $d$ is conjugate to $b$;

\item[(b)]
$c$ is conjugate to $b$ and $d$ is conjugate to $a$, and $|n-m|\leqslant N$.
\end{enumerate}
\end{enumerate}
\end{prop}

\medskip

{\it Proof.} By~\cite[Theorem 6.8]{DOG}, there exists a subset $X_1\subseteq G$ such that
$\{\langle a\rangle,\langle b\rangle\}\hookrightarrow_h (G,X_1)$.
Then, by~\cite[Theorem 5.4]{Osin_1}, there exists a subset $X\subseteq G$ such that $X_1\subseteq X$
and the following conditions hold.
\begin{enumerate}
\item[(1)] $\{\langle a\rangle,\langle b\rangle\}\hookrightarrow_h (G,X)$.
\item[(2)] $\Gamma(G,X\sqcup \mathcal{E})$ is hyperbolic and acylindrical, where $\mathcal{E}=\langle a\rangle\cup \langle b\rangle$.
\end{enumerate}

We set $Y=X\cup \mathcal{E}$. Let us analyze the equation $c^nd^m=a^nb^m$.

\medskip

{\bf Case 1.} Suppose that $c,d$ are loxodromic with respect to $Y$.\\
Then, by Lemma~\ref{lem 3.4}, if $n,m\geqslant 2N_0$ and $E_G(c)\neq E_G(d)$, then $|c^nd^m|_Y>2$.
On the other hand, $|c^nd^m|_Y=|a^nb^m|_Y\leqslant 2$. Therefore we have $E_G(c)=E_G(d)$ if $n,m\geqslant 2N_0$.


\medskip

{\bf Case 2.} Suppose that $c$ is loxodromic and $d$ is elliptic with respect to $Y$.\\
Then, by Lemma~\ref{lem 3.6}, the following holds: if $n\geqslant 2N_1$ and $d^m\notin E_G(c)$, then $|c^nd^m|_Y>2$.
On the other hand, $|c^nd^m|_Y=|a^nb^m|_Y\leqslant 2$.
Therefore we have $d^m\in E_G(c)$ if $n\geqslant 2N_1$.

\medskip

{\bf Case 3.} Suppose that $d$ is loxodromic and $c$ is elliptic with respect to $Y$.
Then, analogously to Case 2, we obtain $c^n\in E_G(d)$ if $m\geqslant 2N_1$.

\medskip

{\bf Case 4.} Suppose that $c,d$ are elliptic with respect to $Y$.

In this case we want to apply Corollary~\ref{quasimorphisms}. Let $q_a:G\rightarrow \mathbb{R}$ be a homogenous quasi-morphism, such that $q_a(a)=1$, $q_a(b)=0$, and $q_a(g)=0$ for all $g\in {\text{\rm Ell}}(G, Y)$, which are non-commensurable with $a$ and $b$.
By Lemma~\ref{lem 4.5}, if $g$ is commensurable with $b$, then $q_a(g)=0$.
Thus, $q_a(g)=0$ for all $g\in {\text{\rm Ell}}(G,Y)$, which are non-commensurable with $a$.

Analogously, let $q_b:G\rightarrow \mathbb{R}$ be a homogenous quasi-morphism, such that $q_b(a)=0$, $q_b(b)=1$, and $q_b(g)=0$ for all $g\in {\text{\rm Ell}}(G, Y)$, which are non-com\-mensurable with $b$.
We set
$$
N_2=2\max\{D(q_a),D(q_b)\}.
$$
and suppose that $n,m> N_2$.
From the definition of a quasi-morphism, we have $$|q_a (a^nb^m)- q_a (a^n)-q_a (b^m)|\leqslant D(q_a).$$
Then
$$|q_a (a^nb^m)- n|\leqslant D(q_a).\eqno{(7.1)}$$
Analogously,
$$|q_b (a^nb^m)- m|\leqslant D(q_b).\eqno{(7.2)}$$

We also have
$$|q_a (c^nd^m)- q_a (c^n)-q_a (d^m)|\leqslant D(q_a).\eqno{(7.3)}$$

and
$$|q_b (c^nd^m)- q_b (c^n)-q_b  (d^m)|\leqslant D(q_b).$$

\medskip

{\bf Subcase 4.1.} Suppose that there exists $x\in \{a,b\}$ such that $c$ and $d$ are non-commensurable
with $x$. Without loss of generality, we assume that $x=a$.
Then (7.3) implies
$$|q_a (c^nd^m)|\leqslant D(q_a).\eqno{(7.4)}$$
It follows from (7.1) and (7.4) that $n\leqslant 2D(q_a)\leqslant N_2$. A contradiction.

\medskip

{\bf Subcase 4.2.}
Suppose that $c$ is commensurable with $a$ and $d$ is commensurable with $b$.
The former means that there exist $u\in G$ and $s_1,s_2\in \mathbb{Z}\setminus \{0\}$ such that
$c^{s_1}=u^{-1}a^{s_2}u$.
Then $c\in E_G(u^{-1}au)$. Since $a$ is special, we have $E_G(u^{-1}au)=\langle u^{-1}au\rangle$. Hence $c=u^{-1}a^su$
for some $s\in \mathbb{Z}\setminus \{0\}$.
Then (7.3) implies
$$|q_a(c^nd^m)-sn|\leqslant D(q_a).$$
Using (7.1), we deduce
$$|n-sn|\leqslant 2D(q_a)\leqslant N_2,$$
Since $n> N_2$, we have $s=1$, i.e. $c=u^{-1}au$. Analogously, $d$ and $b$ are conjugate.

\medskip

{\bf Subcase 4.3.}
Suppose that $c$ is commensurable with $b$ and $d$ is commensurable with $a$.
The latter means that there exist $v\in G$ and $t_1,t_2\in \mathbb{Z}\setminus \{0\}$ such that
$d^{t_1}=v^{-1}a^{t_2}v$.
Then $d\in E_G(v^{-1}av)=\langle v^{-1}av\rangle$. Hence $d=v^{-1}a^{t}v$
for some $t\in \mathbb{Z}\setminus \{0\}$.
Then (7.3) implies
$$|q_a(c^nd^m)-tm|\leqslant D(q_a).$$
Using (7.1), we deduce
$$|n-tm|\leqslant 2D(q_a)\leqslant N_2.\eqno{(7.5)}$$
Analogously, using the commensurability of $c$ and $b$, we deduce
$$
|m-sn|\leqslant N_2\eqno{(7.6)}
$$
for some $s\in \mathbb{Z}\setminus \{0\}$.
It follows from $n,m>N_2$ and (7.5), (7.6) that
$s,t\geqslant 1$.
Without loss of generality, we assume that $m\geqslant n$.
Then
$$N_2\geqslant |n-tm|=tm-n\geqslant (t-1)m\geqslant (t-1)N_2.$$
If $N_2\neq 0$, then we have $t=1$ and $|n-m|\leqslant N_2$.
If $N_2=0$, we deduce from (7.5), (7.6) directly that $n=tm$ and $m=sn$, hence $n=m$ and $|n-m|=N_2$.
\noindent
Taking into account all considered cases, we can set $N=\max\{2N_0,2N_1,\lceil N_2\rceil\}.$
\hfill $\Box$

\medskip


The following lemma says that if $n,m$ in Proposition~\ref{lem 0.1} have a certain common divisor,
then only the subcase (a) in the conclusion of this proposition is possible, i.e. $c$ is conjugate to $a$ and $d$ is conjugate to $b$.
A description of these conjugates will be given in Corollary~\ref{prop}.

\medskip

\begin{lem}\label{lem 0.1.1}
Let $G$ be an acylindrically hyperbolic group with respect to a generating set $Z$.
Suppose that $a,b\in G$ are two non-commensurable special elements (with respect to $Z$).
Then there exists $\ell\in \mathbb{N}$ such that for all $n,m\in \ell \mathbb{N}$, $n\neq m$, the following holds:

If $a^nb^m=c^nd^m$, where $c,d\in G$, then
$c$ is conjugate to $a$ and $d$ is conjugate to $b$.
\end{lem}

{\it Proof.} We use the generating set $Y$ as in the proof of Proposition~\ref{lem 0.1}.
In particular, $\langle a\rangle\cup \langle b\rangle\subseteq Y$. We also use the homogenous quasi-morphisms
$q_a:G\rightarrow \mathbb{R}$ and $q_b:G\rightarrow \mathbb{R}$ defined there. In particular, $q_a(a)=1$, $q_a(b)=0$,
$q_b(a)=0$, and $q_b(b)=1$.



\medskip

$\bullet$ By Lemma~\ref{elem_index},
there exists a number $L\in \mathbb{N}$ such that for every loxodromic element $g\in G$
(with respect to $Y$), the elementary subgroup $E_G(g)$ contains a normal cyclic subgroup of index $L$.

$\bullet$ By Lemma~\ref{blue}, there exist $\varkappa\geqslant 1$, $\varepsilon_0\geqslant 0$ and $n_0\in \mathbb{N}$ such that for any loxodromic (with  respect to $Y$) element $g\in G$ and any $n\geqslant n_0$
we have $$|g^n|_Y\geqslant \frac{1}{\varkappa}n-\varepsilon_0.$$

$\bullet$ We set $$\ell=LMNn_0,$$ where $M=\lceil(2+\varepsilon_0)\varkappa\rceil+1$ and $N\in \mathbb{N}$ is the number from Proposition~\ref{lem 0.1}.
Recall that in the proof of this proposition, we defined $N$ so that we have
$$N\geqslant 2\max\{D(q_a),D(q_b)\}.$$

\medskip

Suppose that $a^nb^n=c^nd^m$, where $n,m\in \ell\mathbb{N}$, $n\neq m$.
We analyze cases in the conclusion of Proposition~\ref{lem 0.1}.
\begin{enumerate}
\item[1)]
$c$ and $d$ are loxodromic with respect to $Y$ and $E_G(c)=E_G(d)$.\\
By Lemma~\ref{elem_index}, there exists $z\in E_G(c)$ such that $c^L$ and $d^L$ are powers of $z$,
say $c^L=z^s$ and $d^L=z^t$.
Then $c^nd^m=z^{sn+mt/L}$. Observe that $sn+mt\neq 0$, otherwise $a^nb^m=c^nd^m=1$ that is impossible by
non-commensurability of $a$ and $b$. Since $\ell=MLNn_0$ is a divisor of $n$ and $m$, we have
(using Lemma~\ref{blue}) that
$$\hspace*{10mm}2\geqslant |a^nb^m|_Y=|z^{sn+mt/L}|_Y\geqslant \frac{1}{\varkappa}\frac{|sn+mt|}{L}-\varepsilon_0\geqslant
\frac{1}{\varkappa}\frac{\ell}{L} -\varepsilon_0\geqslant
\frac{1}{\varkappa}M -\varepsilon_0>2.$$
A contradiction.

\medskip

\item[2)] $c$ is loxodromic with respect to $Y$ and $d$ is elliptic, and $d^m\in E_G(c)$.\\
By definition of $L$, the group $E_G(c)$ contains a normal infinite cyclic subgroup $C$ of index $L$. Hence
$c^L\in C$ and $d^{-m}c^Ld^m=c^{\pm L}$. Since $L$ is a divisor of $n$, this implies $$d^{-m}c^nd^m=c^{\pm n}.$$
The group $C$ is generated by a loxodromic element, because it contains the loxodromic element $c^L$.
Since $d^{mL}\in C$ and $d$ is elliptic, we have $$d^{mL}=1.$$

{\bf Subcase 1.} Suppose that $d^{-m}c^nd^m=c^n$. Then $$(a^nb^m)^L=(c^nd^m)^L=c^{nL}d^{mL}=c^{nL}.$$
It follows
$$\hspace*{10mm}2L\geqslant |(a^nb^m)^L|_Y=|c^{nL}|_Y\geqslant \frac{1}{\varkappa}nL-\varepsilon_0\geqslant
\frac{1}{\varkappa}ML -\varepsilon_0L>2L.$$
A contradiction.

{\bf Subcase 2.} Suppose that $d^{-m}c^nd^m=c^{-n}$. Then
$$(a^nb^m)^{2}=(c^nd^m)^{2}=d^{2m}.$$
Recalling that $d^{mL}=1$, we deduce $(a^nb^m)^{2L}=1$.
Since homogeneous quasi-morphisms vanish on periodic elements, we have $q_a(a^nb^m)=0$.
In view of $|q_a(a^nb^m)-q_a(a^n)-q_a(b^m)|\leqslant~D(q_a)$, this implies that $n\leqslant D(q_a)\leqslant N< \ell$. A contradiction.

\medskip

\item[3)] $d$ is loxodromic with respect to $Y$ and $c$ is elliptic, and $c^n\in E_G(d)$.\\
This case is impossible by the same reason as the previous one.

\medskip

\item[4)] $c$ and $d$ are elliptic with respect to $Y$ and one of the following holds.
\begin{enumerate}
\item[(a)]
$c$ is conjugate to $a$ and $d$ is conjugate to $b$.

\item[(b)]
$c$ is conjugate to $b$ and $d$ is conjugate to $a$, and $|n-m|\leqslant N$.
\end{enumerate}

The subcase (b) is impossible since $n,m\in \ell\mathbb{N}$, $n\neq m$, and $N$ is a proper divisor of $\ell$.
Thus, only the case (a) is possible.
\hfill $\Box$
\end{enumerate}



\begin{rmk}
{\rm
The condition on $\gcd(n,m)$ in Lemma~\ref{lem 0.1.1} cannot be replaced
by the condition that $n,m$ are sufficiently large.
Indeed, if $\gcd(n,m)=1$, then the equation $x^ny^m=a^nb^m$ in the free group $F(a,b)$ of rank 2
has infinitely many solutions $(x,y)=\bigl((a^nb^m)^s,(a^nb^m)^t\bigr)$, where $s,t$ are integers satisfying $ns+mt=1$.
Non of the components of these solutions is conjugate to $a$ or~$b$.
}
\end{rmk}

\section{Isolated components in geodesic polygons}

In the following proof we use Proposition 4.14 from~\cite{DOG}.
Since this proposition and accompanied definitions are crucial in the following proof, we recall them here.

Let $G$ be a group, $\{H_{\lambda}\}_{\lambda\in \Lambda}$ a collection of subgroups of $G$,
$X$ a symmetrised subset of $G$. We assume that $X$ together with $\{H_{\lambda}\}_{\lambda\in \Lambda}$ generates $G$. Let $\mathcal{H}=\bigsqcup_{\lambda\in\Lambda}H_{\lambda}$.

\begin{defn}
{\rm (see~\cite[Definition 4.5]{DOG})
Let $q$ be a path in the Cayley graph $\Gamma(G,X\sqcup \mathcal{H})$. A (non-trivial) subpath $p$ of $q$
is called an {\it $H_{\lambda}$-subpath}, if the label of $p$ is a word in the alphabet $H_{\lambda}$.
An $H_{\lambda}$-subpath $p$ of $q$ is an {\it $H_{\lambda}$-component} if $p$ is not contained in a longer subpath of $q$ with this property. Two $H_{\lambda}$-components $p_1,p_2$ of a path $q$ in $\Gamma(G,X\sqcup \mathcal{H})$ are called {\it connected} if there exists a path $\gamma$ in $\Gamma(G,X\sqcup \mathcal{H})$ that connects some vertex of $p_1$ to some vertex of $p_2$, and ${\text{\bf Lab}}(\gamma)$ is a word consisting only of letters of
$H_{\lambda}$.

Note that we can always assume that $\gamma$ has length at most 1 as every element of $H_{\lambda}$
is included in the set of generators. An $H_{\lambda}$-component $p$ of a path $q$ in $\Gamma(G,X\sqcup \mathcal{H})$ is {\it isolated} if it is not connected to any other component of $q$.

}
\end{defn}

Recall that definitions of a weakly hyperbolic group and of a relative metric $\widehat{d}_{\lambda}$ on $H_{\lambda}$
were given in Section 2.



\begin{defn}\label{n-gon} {\rm (see~\cite[Definition 4.13]{Osin_1})
Let $\varkappa \geqslant 1$, $\varepsilon\geqslant 0$, and $n\geqslant 2$. Let $\mathcal{P}=p_1\dots p_n$ be an $n$-gon in $\Gamma(G, X\sqcup \mathcal{H})$ and let $I$ be a subset of the set of its sides $\{p_1,\dots,p_n\}$
such that:

1) Each side $p_i\in I$ is an isolated $H_{\lambda_i}$-component of $\mathcal{P}$ for some $\lambda_i\in \Lambda$.

2) Each side $p_i\notin I$ is a $(\varkappa,\varepsilon)$-quasi-geodesic.

\medskip
\noindent
We denote $s(\mathcal{P},I)=\underset{p_i\in I}{\sum} \widehat{d}_{\lambda_i}(1,{\text{\rm \bf Lab}}(p_i))$.
}
\end{defn}

\begin{prop}\label{Proposition_Osin}
{\rm (see~\cite[Proposition 4.14]{DOG})}
Suppose that $G$ is weakly hyperbolic relative to $X$ and $\{H_{\lambda}\}_{\lambda\in \Lambda}$.
Then for any $\varkappa\geqslant 1$, $\varepsilon\geqslant 0$, there exists a constant $D(\varkappa,\varepsilon)>0$ such that
for any $n$-gon $\mathcal{P}$ in $\Gamma(G,X\sqcup \mathcal{H})$ and any subset $I$ of the set of its sides satisfying conditions
of Definition~\ref{n-gon}, we have $s(\mathcal{P},I)\leqslant D(\varkappa,\varepsilon)n$.
\end{prop}


\begin{cor}\label{Cor_1.2} Suppose that $G$ is weakly hyperbolic relative to $X$ and $\{H_{\lambda}\}_{\lambda\in \Lambda}$.
Let $\mathcal{P}=p_1p_2p_3$ be a geodesic triangle in $\Gamma(G,X\sqcup \mathcal{H})$,
where $p_3$ is an isolated component of $\mathcal{P}$ or a degenerate path.
Suppose that $q$ is an $H_{\lambda}$-component of $\mathcal{P}$ of the form $q=q_1q_2$, where $q_1$ is a terminal subpath of $p_1$ and $q_2$ is an initial subpath of $p_2$. Then
$\widehat{d}_{\lambda}(1, {\rm \bf Lab}(q))\leqslant 4D(1,0)$.
\end{cor}

\medskip

{\it Proof.}  Let $q_1'$ and $q_2'$ be paths such that $p_1=q_1'q_1$ and $p_2=q_2q_2'$.
Consider the 4-gon $\mathcal{P}'=q_1'qq_2'p_3$. The $H_{\lambda}$-component $q$ of $\mathcal{P}'$
cannot be connected to $p_3$ by assumption and it cannot be connected to an $H_{\lambda}$-component of $q_1'$ or $q_2'$,
since $p_1$ and $p_2$ are geodesics. Therefore $q$ is an isolated component in $\mathcal{P}'$,
and we are done by Proposition~\ref{Proposition_Osin}.\hfill $\Box$

\section{Perfect equations of kind $x^ny^m=a^nb^m$ in acylindrically hyperbolic groups}

The first proposition in this section describes conjugators in the conclusion of Lemma~\ref{lem 0.1.1}.
From these lemma and proposition, we deduce Corollary~\ref{prop}, which gives a clear description of solutions of the equation $x^ny^m=a^nb^m$
in acylindrically hyperbolic groups in the case where  $a,b$ are non-commensurable and special, and $n,m$ have a certain common divisor.

\medskip

\begin{prop}\label{lem 1.1}
Let $G$ be an acylinrically hyperbolic group with respect to a generating set $Z$.
Suppose that $a,b\in G$ are two non-commensurable special elements (with respect to $Z$).
Then there exists $N\in \mathbb{N}$ such that for any $n,m>N$ and any $u,v\in G$ satisfying
$$(u^{-1} a^{n} u) (v^{-1} b^m v) = a^nb^m,$$
there exists $r\in \mathbb{Z}$ such that $u\in \langle a\rangle(a^nb^m)^r$
and $v\in \langle b\rangle (a^nb^m)^r$.
\end{prop}

\medskip

We give a proof of this proposition after introducing of some auxiliary definitions and lemmas.
These lemmas will be proved at the end of this section.

We set $H_a=\langle a\rangle$, $H_b=\langle b\rangle$, and $\mathcal{H}=H_a\sqcup H_b$.
By~\cite[Theorem 6.8]{DOG}, there exists a symmetrized subset $X$ of $G$ such that $\{H_a,H_b\}\hookrightarrow_h(G,X)$.
In particular, $G$ is weakly hyperbolic relative to $X$ and $\{H_a,H_b\}$.
The associated relative metrics on $H_a$ and on  $H_b$ are denoted by $\widehat{d}_a$
and  $\widehat{d}_b$, respectively.

For an element $g=a^i\in H_a$, the number $|i|$ is called the {\it $a$-length} of $g$.
Analogously we define the {\it $b$-length} of an element $g\in H_b$.

A word $w$ in the alphabet $X\sqcup \mathcal{H}$ is called {\it geodesic} if
it has minimal length among all words, representing the same element in $G$ as $w$.
In particular, $w$ does not contain two consecutive letters which both
lie in $H_a$ or both lie in~$H_b$.

\medskip

{\bf Definition of the complexity of a word.}
Given a geodesic word $w$ in the alphabet $X\sqcup \mathcal{H}$, we define its {\it complexity} ${\it Compl}(w)$ as the pair
$(|w|_{X\sqcup {\mathcal{H}}}, s)$, where $s$ is the sum of $a$-lengths of its $H_a$-components
plus the sum of $b$-lengths of its $H_b$-components.
We order the pairs lexicographically: $(t_1',t_2') \prec(t_1,t_2)$ if $t_1'<t_1$ or $t_1=t_1'$ and $t_2'<t_2$.

\medskip

For an element $g\in G$, we define its complexity ${\it Compl}(g)$ as the
minimum of complexities of geodesic words
$w$ in the alphabet $X\sqcup \mathcal{H}$ representing $g$. Note that there is only finitely many elements in each descending chain of complexities.

For any pair $(u,v)$ of elements of $G$, we define its complexity as follows:
$$Compl(u,v)=(Compl(u),Compl(v)).$$
We write $$Compl(u_1,v_1)< Compl(u,v)$$ if $Compl(u_1)\prec Compl(u)$ and $Compl(v_1)\prec Compl(v)$.

\medskip

{\bf Definition of the number $N$.}

We have observed that $G$ is weakly hyperbolic relative to $X$ and $\{H_a,H_b\}$.
Let $D=D(1,0)$ be the constant from Proposition \ref{Proposition_Osin} for parameters $(\varkappa,\varepsilon)=(1,0)$.
We set
$$
\begin{array}{l}
N_a=\max\{i\,|\, \widehat{d}_a(1,a^i)\leqslant 9D\}, \hspace*{3mm} N_b=\max\{i\,|\, \widehat{d}_b(1,b^i)\leqslant 9D\},\\
\end{array}
$$
Since the spaces $(H_a, \widehat{d}_a)$ and $(H_b, \widehat{d}_b)$ are locally finite, the numbers $N_a$ and $N_b$ are finite. We set
$$N=4\cdot \max\{N_a,N_b\}.$$

\medskip

We will prove that this $N$ satisfies Proposition~\ref{lem 1.1}.
For the rest of the proof we assume that $G$, $a$ and $b$ satisfy assumptions of this proposition and that $n,m> N$.

\medskip

Consider the following equations in variables $x,y$:
$$(x^{-1} a^n x) (y^{-1} b^m y) = a^nb^m,\eqno{(9.1)}$$

$$(x^{-1} a^n x) (y^{-1} b^m y) = b^ma^n.\eqno{(9.2)}$$

\begin{lem}\label{auxil_1}
Suppose that $(u,v)$ is a solution of (9.1) such that $u\notin \langle a\rangle$, $v\notin \langle b\rangle$.
We set $(u_1,v_1):=(ua^n,va^n)$ and $(u_2,v_2):=(ub^{-m},vb^{-m})$.
Then $(u_1,v_1)$ and $(u_2,v_2)$ are solutions of (9.2), and we have
$$
Compl(u_1,v_1)<Compl(u,v)\hspace*{3mm}{\text{\rm or}}\hspace*{3mm}Compl(u_2,v_2)<Compl(u,v).
$$
\end{lem}

\medskip

The following lemma is dual to Lemma~\ref{auxil_1}.

\begin{lem}\label{auxil_2}
Suppose that $(p,q)$ is a solution of (9.2) such that $p\notin \langle a\rangle$, $q\notin \langle b\rangle$.
We set $(p_1,q_1):=(pa^{-n},qa^{-n})$ and $(p_2,q_2):=(pb^m,qb^m)$.
Then $(p_1,q_1)$ and $(p_2,q_2)$ are solutions of (9.1), and we have
$$
Compl(p_1,q_1)<Compl(p,q)\hspace*{3mm}{\text{\rm or}}\hspace*{3mm}Compl(p_2,q_2)<Compl(p,q).
$$
\end{lem}

\medskip

Proofs of these lemmas will be given later.

\medskip

{\it Proof of Proposition~\ref{lem 1.1}}.
Suppose that $(u,v)$ is a solution of (9.1).
If $u\in \langle a\rangle$, then $v^{-1}b^mv=b^m$, hence $v\in E_G(b)=\langle b\rangle$,
and we are done. Analogously, if $v\in \langle b\rangle$, then $u\in \langle a\rangle$, and we are done.
Thus, we may assume that $u\notin \langle a\rangle$ and $v\notin \langle b\rangle$.

By Lemma~\ref{auxil_1}, $(ua^n,va^n)$ and $(ub^{-m},vb^{-m})$ are solutions of (9.2) and we have
$
Compl(ua^n,va^n)<Compl(u,v)\hspace*{3mm}{\text{\rm or}}\hspace*{3mm}Compl(ub^{-m},vb^{-m})<Compl(u,v).
$
We consider only the first case $$Compl(ua^n,va^n)<Compl(u,v),\eqno{(9.3)}$$
since the second case can be considered analogously.

We may assume that $(ua^n,va^n)$ satisfies assumption of Lemma~\ref{auxil_2}.
Indeed, the first assumption $ua^n\notin \langle a\rangle$ is satisfied, since $u\notin \langle a\rangle$. Suppose that the second assumption is not satisfied, i.e. $va^n\in \langle b\rangle$. Then $v\in \langle b\rangle \cdot (a^nb^m)^{-1}$, and we deduce from (9.1) that $u^{-1}a^nu \cdot a^nb^ma^{-n}=a^nb^m$. It follows that $ua^nb^m\in E_G(a)=\langle a\rangle$. Hence, $u\in \langle a\rangle (a^nb^m)^{-1}$, and we are done.

Thus, we assume that $(ua^n,va^n)$ satisfies assumption of Lemma~\ref{auxil_2}. By (9.3),
the first case in the conclusion of this lemma cannot happen. Therefore we have the second case, i.e.
$(ua^nb^m, va^nb^m)$ satisfies (9.1) and
$$Compl(ua^nb^m,va^nb^m)<Compl(ua^n,va^n).$$
This formula and (9.3) imply that $Compl(ua^nb^m,va^nb^m)<Compl(u,v)$,
and the statement of Proposition~\ref{lem 1.1} follows by induction.\hfill $\Box$

\bigskip

{\it Proof of Lemma~\ref{auxil_1}}.
By an abuse of notation, for any element $w\in G$, we also denote by $w$ a geodesic word of minimal complexity among all geodesic words in $X\sqcup \mathcal{H}$ representing the element $w$.

Let $(u,v)$ be a solution of (9.1) satisfying assumption of Lemma~\ref{auxil_1}.
Obviously, $(u_1,v_1)$ and $(u_2,v_2)$ are solutions of (9.2). Thus, it suffices to prove that
one of the following holds:

(a) Both $u$ and $v$ end with a power of $a$ which is  smaller than $-n/2$.

(b) Both $u$ and $v$ end with a power of $b$ which is  larger than $m/2$.

\medskip

Since $n,m>N$, we have
$$
n>4N_a=4\max\{i\,|\, \widehat{d}_a(1,a^i)\leqslant 9D\}, \eqno{(9.4)}
$$
$$
m>4N_b=4\max\{i\,|\, \widehat{d}_b(1,b^i)\leqslant 9D\}.\eqno{(9.5)}
$$
Let $\mathcal{P}=p_1p_2\dots p_8$ be a geodesic 8-gon in $\Gamma(G,X\sqcup \mathcal{H})$ with sides $p_i$
labelled by consecutive syllables of the word
$$u^{-1}a^{n}u\, v^{-1}b^{m}vb^{-m}a^{-n}.$$
Observe that the sides $p_2,p_5,p_7,p_8$ of $\mathcal{P}$ are edges labelled by powers of $a$ and~$b$,
see Fig. 13.

\hspace*{2mm}
\includegraphics[scale=0.7]{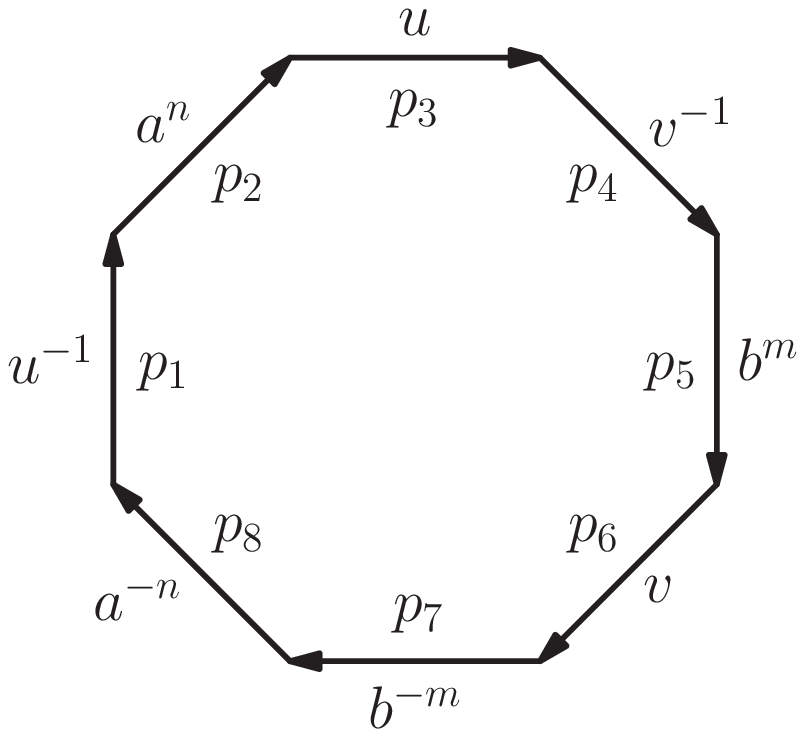}

\vspace*{-120mm}
\begin{center}
Fig. 13.
\end{center}

\medskip

By assumption of lemma, the words $u$ and $v$ are nonempty.
Write $u=a^iu'$ and $v=b^jv'$, where $i,j\in \mathbb{Z}$, the first letter of $u'$ does not lie in $H_a$,
and the first letter of $v'$ does not lie in $H_b$. As $(u,v)$, the pair $(u',v')$ is also a solution of equation (9.1) satisfying assumption of lemma. Obviously, if we prove it for $(u',v')$, then we prove it for $(u,v)$ too.

Thus, we may assume that the first letter of $u$ is not a nontrivial power of $a$, and
the first letter of $v$ is not a nontrivial power of $b$.


Then $p_2$ is an $H_a$-component of $\mathcal{P}$ and $p_5$
is an $H_b$-component of $\mathcal{P}$. By (9.4) and (9.5),
$p_2$ and $p_5$ are not isolated in $\mathcal{P}$ (see Proposition~\ref{Proposition_Osin}).
Then $p_2$ is connected to a component of $p_4$ or $p_6$ and $p_5$ is connected to a component of $p_1$
or~$p_3$.

\medskip

{\it Case 1.} Suppose that $p_2$ is connected to an $H_a$-component of $p_6$.

Then there exists a geodesic rectangle $\mathcal{P}_1=p_2r_1o_2r_2$, where $o_2$ is an $H_a$-com\-ponent of $p_6$ (see Fig.~14a). Let $o_1,o_3$ be subpaths of $p_6$ such that $p_6=o_1o_2o_3$. We consider two complementary geodesic 5-gons  $\mathcal{P}_2=p_3p_4p_5o_1\overline{r_1}$ and $\mathcal{P}_3=p_1\overline{r_2}o_3p_7p_8$.

The path $p_5$ is not an isolated $H_b$-component of $\mathcal{P}_2$ (by Proposition~\ref{Proposition_Osin} and formula (9.5)).
Therefore $p_5$ is connected to some $H_b$-component $\beta$ of $p_3$ (see Fig.~14b).
Then $p_3=\alpha\beta\gamma$ for some subpaths $\alpha,\gamma$ of $p_3$.
Let $\delta$ be a geodesic $b$-path from $(p_5)_{-}$ to $\beta_{+}$.
We consider the triangle $\Delta_1=\gamma p_4\delta$.
The path $\delta$ is an isolated $H_b$-component of $\Delta_1$ or a degenerate path.

\vspace*{-20mm}
\hspace*{-35mm}
\includegraphics[scale=0.7]{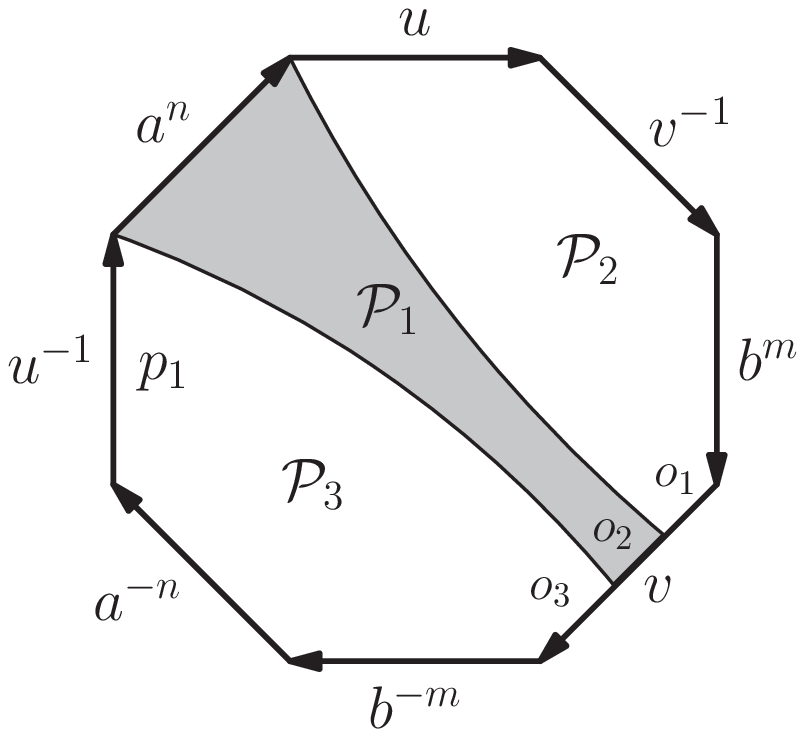}

\vspace*{-208mm}
\hspace*{40mm}
\includegraphics[scale=0.7]{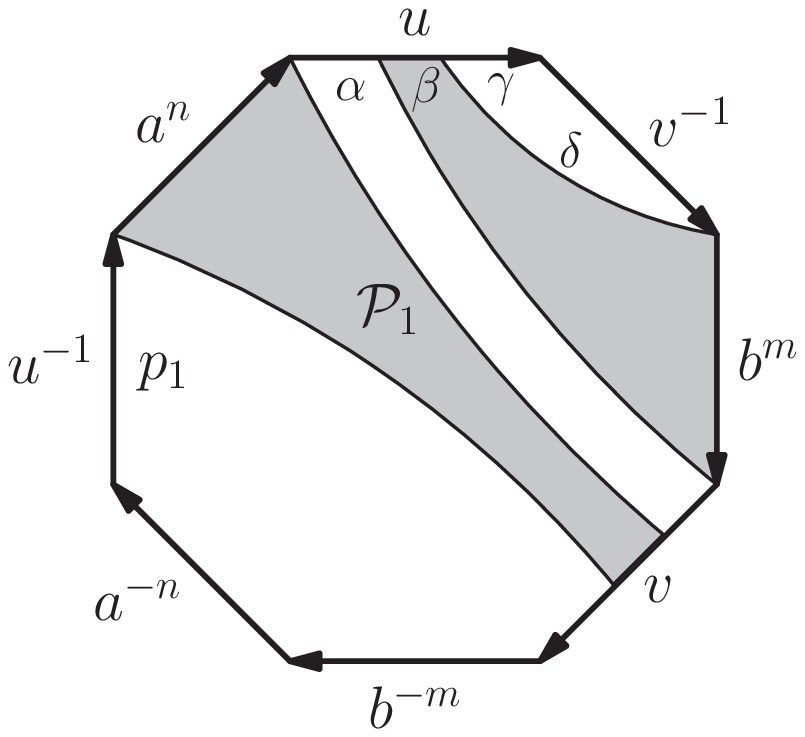}

\vspace*{-120mm}
\hspace*{40mm}
\begin{center}
Fig. 14a and 14b.
\end{center}

\medskip




In the rest of the proof we use the following notation.
For any nontrivial path $p$ in the Cayley graph $\Gamma(G,X\sqcup \mathcal{H})$, let $p^{\circ}$ and $p^{\bullet}$ denote the first and the last edges  of $p$, respectively.


\medskip

Let $t$ be the $H_b$-component of $\mathcal{P}$ containing the edge $p_7$. Then $t$ is contained in $o_3p_7$ (see Fig. 15).

\medskip

{\it Case 1.1.} Suppose that $t$ is not connected to a component of $p_1$.

Then $t$ is isolated in the the 5-gon $\mathcal{P}_3$.

\medskip

\vspace*{-25mm}
\hspace*{2mm}
\includegraphics[scale=0.7]{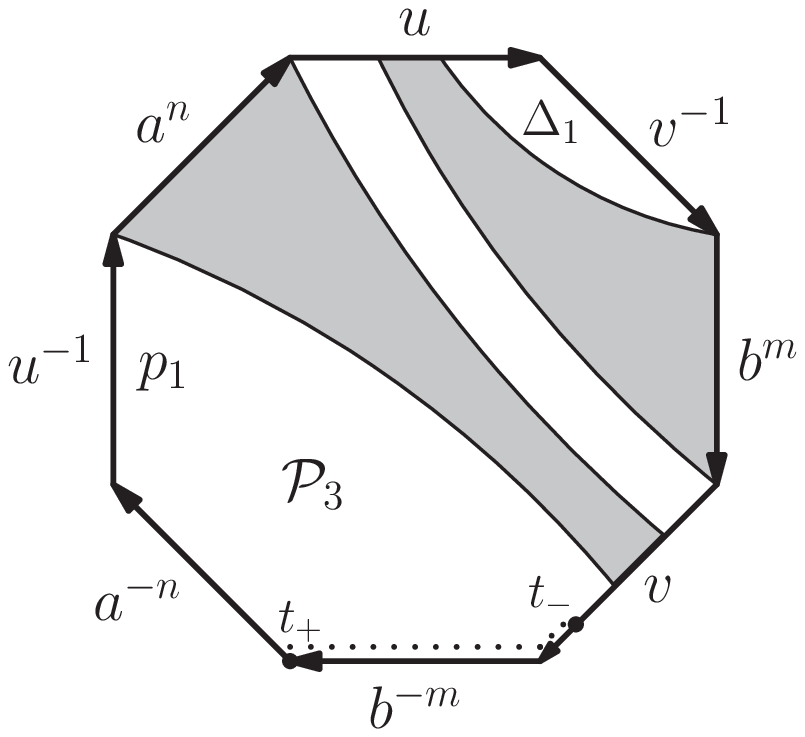}

\vspace*{-123mm}
\begin{center}
Fig. 15.
\end{center}

\medskip

By Proposition~\ref{Proposition_Osin} applied to $\mathcal{P}_3$, we obtain
$$\widehat{d}_b(1, {\rm \bf Lab}(t))\leqslant 5D.\eqno{(9.6)}$$
By (9.5), we have
$$\widehat{d}_b(1, {\rm \bf Lab}(p_7))=\widehat{d}_b(1, b^{-m})> 9D.$$
Therefore $p_7$ is a proper subpath of $t$, and hence $t=p_6^{\bullet}p_7$
with
$$
\widehat{d}_b(1, {\rm \bf Lab}(p_6^{\bullet}))>4D.
$$
Since $\mathbf{Lab}(p_6)=\mathbf{Lab}(\overline{p_4})$,
we have $\mathbf{Lab}(p_6^{\bullet})=(\mathbf{Lab}(p_4^{\circ}))^{-1}$. Hence
$$
\widehat{d}_b(1, {\rm \bf Lab}(p_4^{\circ}))>4D.\eqno{(9.7)}
$$

We claim that the $b$-path $p_4^{\circ}$ cannot be a component
of $\Delta_1$. Indeed, if it were, we could apply Corollary~\ref{Cor_1.2} to the triangle $\Delta_1$, its side $\delta $ (which is an isolated $H_b$-component of $\Delta_1$ or a degenerate path), and to the component $p_4^{\circ}$,
and get a contradiction to (9.7).


Hence, $p_3^{\bullet}p_4^{\circ}$ is a component of $\Delta_1$ and, by Corollary~\ref{Cor_1.2}, we have
$$\widehat{d}_b(1,{\rm\bf Lab}(p_3^{\bullet}p_4^{\circ}))\leqslant 4D.\eqno{(9.8)}$$

Now we estimate $y,z\in \mathbb{Z}$ such that ${\rm \bf Lab}(p_3^{\bullet})=b^y$ and ${\rm \bf Lab}(p_6^{\bullet})=b^z$. Since ${\rm \bf Lab}(t)=b^{z-m}$ and ${\rm \bf Lab}(p_3^{\bullet}p_4^{\circ})=b^{y-z}$, we deduce from (9.6) and (9.8) that

$$
|z-m|\leqslant N_b \hspace*{2mm}{\rm and} \hspace*{2mm} |y-z|\leqslant N_b.
$$

Since $m>4N_b$, we deduce that $z\geqslant m-N_b>\frac{3}{4}m$ and $y\geqslant z-N_b>\frac{1}{2}m$.

\medskip

Thus, both $u$ and $v$ end with a power of $b$ which is  larger than $m/2$, and we are done.

\medskip

{\it Case 1.2.} Suppose that $t$ is connected to some component $\beta'$ of $p_1$ (see Fig.~16).
Then $p_1=\alpha'\beta'\gamma'$ for some subpaths $\alpha',\gamma'$ of $p_1$.
Let $\delta'$ be a geodesic $b$-path from $\beta'_{-}$ to $(p_7)_{+}$.
We consider the triangle $\Delta_2=p_8\alpha'\delta'$. The path $\delta'$ is either an isolated $H_b$-component of $\Delta_2$ or a degenerate path.
Let $q$ be an $H_a$-component of $\Delta_2$ containing $p_8$. Then $q$ is also isolated in $\Delta_2$.


\vspace*{-20mm}
\hspace*{2mm}
\includegraphics[scale=0.7]{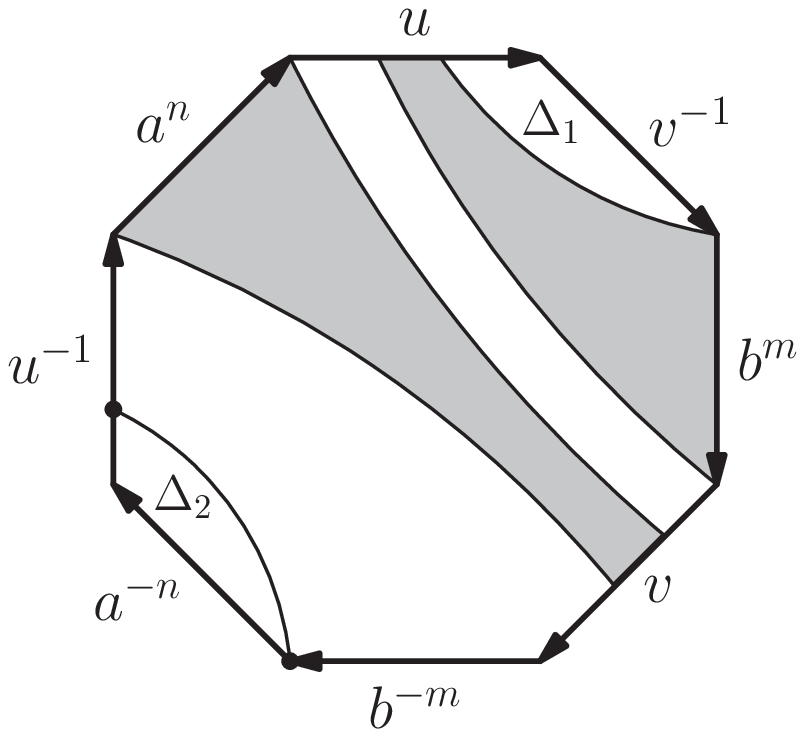}

\vspace*{-120mm}
\hspace*{40mm}
\begin{center}
Fig. 16.
\end{center}

\medskip

By Proposition~\ref{Proposition_Osin}, we have $$\widehat{d}_a(1,{\rm \bf Lab}(q))\leqslant 3D.\eqno{(9.9)}$$

By (9.4) we have $$\widehat{d}_a(1,{\rm\bf Lab}(p_8))=\widehat{d}_a(1,a^{-n})>9D.$$
Therefore $p_8$ is a proper subpath of $q$, and hence $q=p_8p_1^{\circ}$ with
$$\widehat{d}_a(1,\mathbf{Lab}(p_1^{\circ}))>6D.$$
Since $\mathbf{Lab}(p_1)=\mathbf{Lab}(\overline{p_3})$, we have $\mathbf{Lab}(p_1^{\circ})=(\mathbf{Lab}(p_3^{\bullet}))^{-1}$. Hence
$$\widehat{d}_a(1,\mathbf{Lab}(p_3^{\bullet}))>6D.$$

By Corollary~\ref{Cor_1.2} applied to the triangle $\Delta_1$, the $a$-path $p_3^{\bullet}$ cannot be a component
of $\Delta_1$. Hence, $p_3^{\bullet}p_4^{\circ}$ is one and, by this corollary, we have
$$\widehat{d}_a(1,{\rm\bf Lab}(p_3^{\bullet}p_4^{\circ}))\leqslant 4D.\eqno{(9.10)}$$

Now we estimate $z,y\in \mathbb{Z}$ such that ${\rm \bf Lab}(p_1^{\circ})=a^z$ and ${\rm \bf Lab}(p_4^{\circ})=a^y$. Since ${\rm \bf Lab}(q)=a^{-n+z}$ and ${\rm \bf Lab}(p_3^{\bullet}p_4^{\circ})=a^{-z+y}$, we deduce from (9.9) and (9.10) that

$$
|-n+z|\leqslant N_a \hspace*{2mm}{\rm and} \hspace*{2mm} |-z+y|\leqslant N_a.
$$

Since $n>4N_a$, we deduce that $z\geqslant n-N_a>\frac{3}{4}n$ and $y\geqslant z-N_a>\frac{1}{2}n$.
Thus, both $u$ and $v$ end with a power of $a$ which is  smaller than $-n/2$, and we are done.




\medskip

{\it Case 2.} Suppose that $p_2$ is connected to an $H_a$-component of $p_4$.

Arguing as in Case 1, we can prove that $p_5$ is connected to a component of~$p_1$.
After that, renaming $a,b,u,v,n,m$ by $b^{-1},a^{-1},v,u,m,n$, respectively, we reduce to Case 1.
\hfill $\Box$

\medskip

To simplify formulations, we introduce the following definition.

\begin{defn}
{\rm
Let $g\in G$ and $n,m\in \mathbb{Z}$. The equation $x^ny^m=g$ in variables $x,y$ is called {\it perfect} if it has a solution $(x_0,y_0)$ in $G$
and any solution of this equation has the form $(x_0^{g^{\alpha}},y_0^{g^{\alpha}})$ for some $\alpha\in \mathbb{Z}$.
}
\end{defn}

The following corollary directly follows from Lemma~\ref{lem 0.1.1} and Proposition~\ref{lem 1.1}.

\begin{cor}\label{prop}
Let $G$ be an acylinrically hyperbolic group with respect to a generating set $S$.
Suppose that $a,b\in G$ are two non-commensurable special elements (with respect to $S$).
Then there exists a number $\ell=\ell(a,b)\in \mathbb{N}$ such that for all $n,m\in \mathcal{\ell}\mathbb{N}$, $n\neq m$, the equation $x^ny^m=a^nb^m$ is perfect.
\end{cor}







\section{Special generating sets for finitely generated acylindrically hyperbolic groups}

The main purpose of this section is Proposition~\ref{special_generators}.
Other statements of this section will be also used in Sections~11 and~12.
The following lemma proven in~\cite{DOG} is crucial in many proofs.
Therefore we reproduce it here.

\begin{lem} {\rm (see~\cite[Lemma 4.21]{DOG})}\label{second_important}
Let $G$ be a group weakly hyperbolic relative to $X$ and $\{H_{\lambda}\}_{\lambda\in \Lambda}$
and let $\mathcal{W}$ be the set consisting of all words $U$ in $X\sqcup \mathcal{H}$ such that:

\begin{enumerate}
\item[$(W_1)$] $U$ contains no subwords of type $xy$, where $x,y\in X$.

\item[$(W_2)$] If $U$ contains a letter $h\in H_{\lambda}$ for some $\lambda\in \Lambda$, then $\widehat{d}_{\lambda}(1,h)>50 D$, where $D=D(1,0)$ is given by Proposition~\ref{Proposition_Osin}.

\item[$(W_3)$] If $h_1xh_2$ (respectively $h_1h_2$) is a subword of $U$, where $x\in X$, $h_1\in H_{\lambda}$, $h_2\in H_{\mu}$, then either $\lambda\neq \mu$ or the element represented by $x$ in $G$
    does not belong to $H_{\lambda}$ (respectively $\lambda\neq \mu$).
\end{enumerate}

\medskip

Then the following hold.

\begin{enumerate}
\item[(a)] Every path in $\Gamma(G,X\sqcup \mathcal{H})$ labelled by a word from $\mathcal{W}$ is $(4,1)$-quasi-geodesic.


\item[(b)] For every $\varepsilon>0$ and every integer $K>0$, there exist $R=R(\varepsilon,K)>0$ satisfying
the following condition. Let $p,q$ be two paths in $\Gamma(G,X\sqcup \mathcal{H})$ such that $\ell(p)\geqslant R$, ${\rm \bf Lab}(p), {\rm \bf Lab}(q)\in \mathcal{W}$, and $p,q$ are oriented $\varepsilon$-close, i.e.,
$$
\max\{d(p_{-},q_{-}),d(p_{+},q_{+})\}\leqslant \varepsilon.
$$
Then there exist $K$ consecutive components of $p$ which are connected to $K$ consecutive components of $q$.
\end{enumerate}
\end{lem}

\medskip

\begin{cor}\label{free_group}
Suppose that $G$ is a group, $Y\subset G$ a subset, and $E_1,\dots,E_k$ are subgroups of $G$ such that
$$\{E_1,\dots ,E_k\}\hookrightarrow_h (G,Y).$$
Suppose that $a_1,\dots,a_k$ are elements of infinite order from $E_1,\dots,E_k$, respec\-tively.
Then there exists $N\in \mathbb{N}$ such that if $n_1,\dots,n_k\geqslant N$, then
every cyclically reduced word $W$ of syllable length at least~2 in the alphabet $\{a_1^{n_1},\dots ,a_k^{n_k}\}^{\pm}$ represents a loxodromic element of $G$ with respect to $Y\sqcup \mathcal{E}$, where $\mathcal{E}=E_1\sqcup \dots \sqcup E_k$.
In particular, $\langle a_1^{n_1},\dots ,a_k^{n_k}\rangle$ is a free group of rank~$k$.

Moreover, if $\langle Y\rangle=G$, then each word $W$ as above represents a loxodromic element of $G$
with respect to $Y$.

\end{cor}

\medskip

{\it Proof.} For  $i=1,\dots ,k$,
let $\widehat{d}_i$ be the metric on $E_i$ associated with the embedding
$\{E_1,\dots ,E_k\}\hookrightarrow_h (G,Y)$.
Let
$$
N_i:=\max\{n\in \mathbb{N}\,|\, \widehat{d}_i(1,a_i^n)\leqslant 50D\},
$$
where $D=D(1,0)$ is given by Proposition~\ref{Proposition_Osin}.
We claim that
$$
N=\max \{N_i\,|\, i=1,\dots ,k\}+1,
$$
satisfies the corollary. Let $n_1,\dots,n_k\geqslant N$ and
let $W$ be a cyclically reduced word in the alphabet $\{a_1^{n_1},\dots ,a_k^{n_k}\}^{\pm}$
such that the syllable length of $W$ is at least 2.
Using conjugation, we may assume that the first and the last letters of $W$ are not coincide and
not inverse to each other.

Let $U$ be the word in the alphabet
$\langle a_1^{n_1}\rangle\sqcup \dots \sqcup \langle a_k^{n_k}\rangle$
obtained from $W$ by replacing each syllable of $W$ of kind $\underbrace{a_i^{\pm n_i}\dots a_i^{\pm n_i}}_s$ by the unique letter $a_i^{\pm sn_i}$.
Then $U^m$ satisfies conditions $(W_1)-(W_3)$ of Lemma~\ref{second_important} for any $m\in \mathbb{N}$.
Let $p_m$ be the path in $\Gamma(G,Y\sqcup \mathcal{E})$ labelled by $U^m$, such that $(p_m)_{-}=1$. Since, by this lemma, the path $p_m$ is $(4,1)$-quasi-geodesic, we have $d(1,(p_m)_{+})\geqslant \ell(p_m)/4-1\geqslant m/2-1$. Then $U$, and hence $W$, represent a loxodromic element of $G$ with respect to $Y\sqcup \mathcal{E}$.
In particular, $W\neq 1$ in $G$.
If $W$ is of syllable length 1, i.e. $W=a_i^{n_im}$ for some
$i\in \{1,\dots,k\}$ and $m\neq 0$, then, obviously, $W\neq 1$. Therefore
$\langle a_1^{n_1},\dots ,a_k^{n_k}\rangle$ is a free group of rank $k$.

The last statement of corollary obviously follows from the first one.
\hfill $\Box$

\medskip

The following lemma is closely related to~\cite[Corollary 6.12]{DOG}.

\begin{lem}\label{many_non-commensurable}
Let $G$ be a group, $X\subseteq G$, $H\hookrightarrow_h(G,X)$ a finitely generated infinite subgroup.
Then for any finite collection of elements $a_1,\dots,a_s\in G\setminus H$ and any infinite
subset $\widetilde{H}\subseteq H$, there exist elements $h_1,\dots ,h_s\in \widetilde{H}$
such that $a_1h_1,\dots, a_sh_s$ are paarwise non-commensurable
loxodromic elements with respect to the action of $G$ on $\Gamma(G,X\sqcup H)$.
\end{lem}

{\it Proof.}
First we show that conditions $(a')$ and $(b)$ of~\cite[Theorem 6.11]{DOG} are satisfied for some extended relative generating set $X_1$.
Though conclusion of this theorem is not sufficient for our aims, the proof is sufficient and can be easily adopted to obtain the desired statement.

{\it Definition of $X_1$.} Let $B$ be a finite generating set of $H$. We set $$X_1=X\cup B\cup\{a_1,\dots ,a_s\}^{\pm}.$$
Since $|X\bigtriangleup X_1|<\infty$, we have $H\hookrightarrow_h(G,X_1)$. Let $\widehat{d}$ be the relative metric on $H$ associated with this embedding.

\medskip

\noindent
{\it Verification of condition $(a')$.} We shall show that

$\bullet$ $\widehat{d}(1,h)<\infty$ for any $h\in H$ and

$\bullet$ $H$ is unbounded with respect to~$\widehat{d}$.

The former follows from the fact that the relative generating set $X_1$ contains a finite generating set of $H$, namely $B$.
The latter follows from  the fact that the metric space $(H,\widehat{d})$ is locally finite (since $H\hookrightarrow_h(G,X_1)$) and the assumption that $H$ is infinite.

\medskip

\noindent
{\it Verification of condition $(b)$ for the elements $a_1,\dots,a_s$.} We shall check that

$\bullet$ these elements lie in $X_1$ and

$\bullet$ $|H^{a_i}\cap H|<\infty$ for $i=1,\dots,s$.

The former is valid by definition of $X_1$, the latter follows from the assumption that
$a_i\in G\setminus H$ (see Lemma~\ref
{periph_subgr}).

\medskip

Thus conditions $(a')$ and $(b)$ of~\cite[Theorem 6.11]{DOG} are satisfied for the extended relative generating set $X_1$.
Now we look in the proof of this theorem.
Condition $(a')$ and the local finiteness of $(H,\widehat{d})$ enable to choose $h_1,\dots,h_s$ in~$\widetilde{H}$
such that
$$
\begin{array}{ll}
\widehat{d}(1,h_1) & >50D,\vspace*{2mm}\\
\widehat{d}(1,h_{i+1}) & >\widehat{d}(1,h_i)+8D,\hspace*{2mm} i=1,\dots,s-1,
\end{array}
$$
where $D=D(1,0)$ is provided by Proposition~\ref{Proposition_Osin}. We set $f_i=a_ih_i$. Then the proof that $f_1,\dots,f_s$ are non-commensurable and loxodromic
with respect to the action of $G$ on $\Gamma(G,X_1\sqcup H)$ is the same as in~\cite[Theorem 6.11]{DOG}.
Since $X\subseteq \,X_1$,
these elements remain loxodromic with respect to the action of $G$ on $\Gamma(G,X\sqcup H)$.
\hfill $\Box$

\begin{lem}\label{many_special}
Let $G$ be an acylindrically hyperbolic group with respect to a generating set $Z$ and let $a,b\in G$ be
two non-commensurable loxodromic with respect to $Z$ elements, where, additionally, $a$ is special.
Then there exists a positive integer $n_0$ such that for any $n,m\geqslant n_0$ the element $g=a^nb^m$
is special with respect to some generating set, in particular $E_G(g)=\langle g\rangle$.
\end{lem}


{\it Proof.} Since $G$ acts acylindrically on the hyperbolic space $\Gamma(G,Z)$
and $a,b$ are non-commensurable and loxodromic with respect to $Z$, we can apply~\cite[Theorem 6.8]{DOG} which says that in this situation there exists a subset $X\subset G$ such that
$\{E_G(a),E_G(b)\}\hookrightarrow_h(G,X)$. Denote $\mathcal{E}=E_G(a)\sqcup E_G(b)$. By Theorem~\ref{enlarging}, there exists $Y\subseteq G$ such that
$X\subseteq Y$ and the following conditions hold.
\begin{enumerate}
\item[(a)] $\{E_G(a),E_G(b)\}\hookrightarrow_h (G,Y)$. In particular, the Cayley graph
$\Gamma(G, Y\sqcup \mathcal{E})$ is hyperbolic.

\item[(b)] The action of $G$ on $\Gamma(G,Y\sqcup \mathcal{E})$ is acylindrical.
\end{enumerate}

Let $\widehat{d}_1$ and $\widehat{d}_2$ be the relative metrics on $E_G(a)$ and on $E_G(b)$,
respectively, associated with the embedding $\{E_G(a),E_G(b)\}\hookrightarrow_h(G,Y)$.
Then there exists $n_0$ such that for any $l\geqslant n_0$ we have $\widehat{d}_1(1,a^l)>50D$
and $\widehat{d}_2(1,b^l)>50D$, where
$D=D(1,0)$ is the constant from Proposition~\ref{Proposition_Osin}.  Let $n,m\geqslant n_0$ and $g=a^nb^m$.
Observe that conditions $(W_1)-(W_3)$ of Lemma~\ref{second_important} are satisfied for $g$ considered as a word of length $2$ in the alphabet $\mathcal{E}$.
By part (a) of Lemma~\ref{second_important}, the element $g$ is loxodromic with respect to $Y\sqcup \mathcal{E}$.
Since $G$ acts acylindrically on the hyperbolic space $\Gamma(G,Y\sqcup \mathcal{E})$ and contains a loxodromic element with respect to this action and $G$ is not virtually cyclic,
we conclude that this action is non-elementary (see Theorem~\ref{check_elementary}).
Therefore $G$ is acylindrically hyperbolic with respect to $Y\sqcup \mathcal{E}$ and the subgroup $E_G(g)$ is well defined.


The rest of the proof is very similar to the second part of the proof of \cite[Lemma 6.18]{DOG}.
Let $t\in E_G(g)$. Then $tg^{\sigma k}=g^kt$ for some $k\in \mathbb{N}$ and $\sigma\in \{-1,1\}$.
Consider the paths $p$ and $q$ in $\Gamma(G,Y\sqcup \mathcal{E})$ labelled by $(a^nb^m)^k$ and
$(a^nb^m)^{\sigma k}$, respectively, such that $p_{-}=1$, $q_{-}=t$. We have
$d_{Y\sqcup \mathcal{E}}(p_{-},q_{-})=d_{Y\sqcup \mathcal{E}}(p_{+},q_{+})=\varepsilon$, where
$\varepsilon=|t|_{Y\sqcup \mathcal{E}}$.

Let $R=R(\varepsilon,3)$ be as in the part (b) of Lemma~\ref{second_important}. Passing to a multiple of $k$ if necessary, we may assume that $\ell(p)\geqslant R$.
Then by statement (b) of Lemma~\ref{second_important}, there exist~3 consecutive components $p_1,p_2,p_3$ of $p$ that are connected to 3 consecutive components $q_1,q_2,q_3$ of $q$.

Without loss of generality we may assume that $p_1,p_3,q_1,q_3$ are $E_G(a)$-compon\-ents while $p_2,q_2$ are $E_G(b)$-components. Let $e_j$ be a path connecting $(p_j)_{+}$ to $(q_j)_{+}$ in $\Gamma(G,Y\sqcup\mathcal{E})$
and let $z_j$ be the element of $G$ represented by ${\rm \bf Lab}(e_j)$, $j=1,2$.
Then $z_j\in E_G(a)\cap E_G(b)$.
But $E_G(a)\cap E_G(b)=1$ since $a,b$ are non-commensurable and $E_G(a)=\langle a\rangle$ by assumption. Thus, $z_j=1$. Reading the label of the closed path $e_1q_2\overline{e_2}\,\overline{p_2}$ we obtain $b^{\sigma m}b^{-m}=1$.
Therefore $\sigma=1$ and the label of $q$ is $(a^nb^m)^k$.
Reading the labels of the segment of $p$ from 1 to $(p_1)_{+}$, $e_1$, and the segment of $\overline{q}$ from $(q_1)_{+}$ to $t$, we obtain $t=g^l$ for some $l\in \mathbb{Z}$. Therefore $E_G(g)\leqslant \langle g\rangle$. Hence $E_G(g)=\langle g\rangle$ and $g$ is special with respect to $Y\sqcup \mathcal{E}$.
\hfill $\Box$

\begin{lem}\label{one_special}
Suppose that $G$ is an acylindrically hyperbolic group
without nontrivial finite normal subgroups. Then $G$ contains at least one special element.

Moreover, there exist an element $g\in G$ and a generating set $Y$ of $G$ such that
$g$ is special with respect to $Y$ and $\langle g\rangle\hookrightarrow_h(G,Y)$.

\end{lem}


{\it Proof.}
By~\cite[Lemma 6.18]{DOG}, there exist a subset $X\subseteq G$, a subgroup $E\hookrightarrow_h (G,X)$, and an element $g\in G$ such that $E=\langle g\rangle\times K(G)$, where $K(G)$ is the maximal finite normal subgroup of $G$.
By assumption, $K(G)=1$. Then $\langle g\rangle \hookrightarrow_h (G,X)$.

It was shown in the proof of this lemma that $g$ can be chosen to be loxodromic with respect
to some generating set. In particular, we may assume that $g$ has infinite order.
Note that $G$ is not virtually cyclic, since $G$ is acylindrically hyperbolic.

By Lemma~\ref{best_Y}, there exists a generating set $Y$ of $G$ such that
$G$ is acylindrically hyperbolic with respect to $Y$ and $\langle g\rangle\hookrightarrow_h(G,Y)$.
Let $\widehat{d}$ be the relative metric on $\langle g\rangle$ associated with this embedding.
Since the metric space $(\langle g\rangle,\widehat{d}\,)$ is locally finite, $g$ cannot be elliptic with respect to $Y$.
Then $g$ is loxodromic with respect to $Y$.
Since $\langle g\rangle\hookrightarrow_h(G,Y)$, we deduce from Lemma~\ref{periph_subgr} and (3.1) that $E_G(g)=\langle g\rangle$. Thus $g$ is special with respect to $Y$.
\hfill $\Box$

\medskip

\begin{prop}\label{special_generators}
Suppose that $G$ is a finitely generated acylindrically hyperbolic group
without nontrivial finite normal subgroups. Then there exists
a generating set $Y$ of $G$ such that $G$ is acylindrically hyperbolic
with respect to $Y$ and the following holds:

For any $n\in \mathbb{N}$, $G$ can be generated by a finite set $A$ such that $|A|\geqslant n$
and the elements of $A$ are pairwise non-commensurable and special with respect to~$Y$.
\end{prop}

\medskip

{\it Proof.}
By Lemma~\ref{one_special},
there exist an element $g\in G$ and a generating set~$Y$ of $G$ such that
$g$ is special with respect to $Y$
and $E\hookrightarrow_h(G,Y)$, where $E=\langle g\rangle$.
It follows from Definition~\ref{exact_defn_special} that $G$ is acylindrically hyperbolic with respect to $Y$.


Suppose that $G=\langle a_1,\dots,a_l\rangle$.
Removing those $a_i$, which are powers of $g$, we may assume that $G=\langle g,a_1,\dots,a_k\rangle$
for some $1\leqslant k\leqslant l$, where $a_i\notin E$ for $i=1,\dots ,k$.

\medskip



{\bf Step 1.}
We show how to find a finite generating set $B$ of $G$ such that $g\in B$ and the elements of $B$ are pairwise non-commensurable and loxodromic with respect to $Y$.

We set $G_0=\langle g\rangle$ and $G_i=\langle g,a_1,\dots,a_i\rangle$ for $i=1,\dots ,k$. Note that $G=G_k$.
Arguing inductively, we fix $i\in \{0,\dots ,k-1\}$ and suppose that we have found a finite generating set $B_i$ of $G_i$ such that $g\in B_i$ and the elements of $B_i$ are pairwise non-commensurable and loxodromic with respect to $Y$.

We set $s=|B_i|+1$.
By~\cite[Corollary 6.12]{DOG} (or by Lemma~\ref{many_non-commensurable}), there exist positive integers $n_1<n_2<\dots < n_s$ such that the elements
of the set $\{a_{i+1}g^{n_1}, a_{i+1}g^{n_2}, \dots ,a_{i+1}g^{n_s}\}$ are  pairwise non-commensurable and loxodromic with respect to $Y\sqcup E$. It follows that they are loxodromic with respect to $Y$.
Since the number of these elements is $|B_i|+1$, there exists $j\in \{1,\dots,s\}$ such that the elements of $B_{i+1}:=B_i\cup \{a_{i+1}g^{n_j}\}$ are pairwise non-commensurable. Since $g\in B_i$, we deduce that
$G_{i+1}=\langle B_{i+1}\rangle$.
Finally, we set $B=B_k$.

\medskip

Thus, we may assume from the beginning that $G=\langle g,a_1,\dots,a_k\rangle$,
where $g,a_1,\dots,a_k$ are pairwise non-commensurable and loxodromic with respect to~$Y$.
Recall that $g$ is special with respect to $Y$.

\medskip

{\bf Step 2.} We show how to find a finite generating set of $G$ consisting of pairwise non-commensurable and special elements with respect to $Y$.

We set $A_0=\{ g\}$.
Arguing inductively, we fix $i\in \{0,\dots ,k-1\}$ and suppose that we have found a finite subset $A_i\subset G$
such that $g,a_1,\dots,a_i\in \langle A_i\rangle$ and the elements of $A_i$ are pairwise non-commensurable and special with respect to~$Y$. We set $s=2|A_i|+2$ and construct a finite set $A_{i+1}\subset G$ with analogous
properties.

By Lemma~\ref{many_special}, there exists  $n_0\in \mathbb{N}$ such that for any $n,m\geqslant n_0$ the element $a_{i+1}^ng^m$ is special with respect to some generating set of $G$. In particular, we have
$$
E_G(a_{i+1}^ng^m)=\langle a_{i+1}^ng^m\rangle\eqno{(10.1)}
$$
for any $n,m\geqslant n_0$. By Lemma~\ref{many_non-commensurable}, there exist $s$ integers $m_1,m_2,\dots,m_s\geqslant n_0$ such that the elements
$a_{i+1}^{n_0+1}g^{m_1}, a_{i+1}^{n_0+2}g^{m_2},\dots ,a_{i+1}^{n_0+s}g^{m_s}$ are pairwise non-commensurable and loxodromic with respect to $Y\sqcup E$ (and hence with respect to $Y$).
By (10.1) they are special with respect to $Y$.
Since the number of these elements is $2|A_i|+2$, there exists an odd $j\in \{1,\dots ,s-1\}$ such that
the elements of $A_{i+1}:=A_i\cup \{a_{i+1}^{n_0+j}g^{m_j},a_{i+1}^{n_0+j+1}g^{m_{j+1}}
\}$ are pairwise non-commensurable. By construction, all elements of $A_{i+1}$ are special with respect to $Y$.
Since $g\in A_i$, we have $a_{i+1}\in \langle A_{i+1}\rangle$.

Observe that $\langle A_k\rangle=G$ and $|A_k|=2k+1$, where $k\geqslant 1$ is fixed before step~2. Repeating the construction of step 2 several times, we can obtain a finite generating
set $A$ of $G$ with desired properties and of arbitrary large finite cardinality.
\hfill $\Box$

\medskip

The following proposition is a generalization of Proposition~\ref{special_generators}.
It will be used only in Section~14.



\begin{prop}\label{subgroup_gener}
Let $G$ be an acylindrically hyperbolic group
and $K$ a noncyclic finitely generated subgroup of $G$.
If $K$ contains at least one special element of $G$, then
for any $n\in \mathbb{N}$ the subgroup $K$ can be generated
by a finite set $A$ such that $|A|\geqslant n$ and the elements of $A$ are pairwise
non-commensurable and jointly special in $G$.

\end{prop}

\medskip

{\it Proof.}
The proof can be obtained from the proof of Proposition~\ref{special_generators} if we substitute there
$K$ instead of $G$ in appropriate places.
\hfill $\Box$

\medskip

\section{A special case of the first main theorem}

In this section we illustrate the main idea of the proof
of Theorem~\ref{prop 4.1} in a special case, see Proposition~\ref{prop 2.1}.
In the proof of this proposition we use Lemma~\ref{lem 2.2}.


\begin{prop}\label{prop 2.1}
Let $G$ be a group having a finite presentation with one relation,
$G=\langle g_1,\dots, g_n\,|\, R\rangle$, and let $H$ be a subgroup of $G$ satisfying the following properties.

\begin{enumerate}
\item[1)]
$H$ is acylindrically hyperbolic with respect to a generating set $Z$.

\item[2)] $H$ is generated by three elements $a_1,a_2,a_3$ which are pairwise non-commen\-surable and special with respect to $Z$.

\item[3)] $H$ does not have a nontrivial normal finite subgroup.
\end{enumerate}

\medskip

Then $H$ is verbally closed in $G$ if and only if $H$ is a retract of~$G$.
\end{prop}

\medskip

{\it Proof.}
We suppose that $H$ is verbally closed in $G$ and prove that $H$ is a retract of $G$. The converse
statement is obvious.
Let $F=F(x_1,\dots,x_n)$ be a free group of rank $n$.
Let $v_1,v_2,v_3,w$ be words in $F$ such that  $v_i(g_1,\dots,g_n)=~\!a_i$, $i=1,2,3,$ and $w(g_1,\dots,g_n)=R$.
Consider the following equation in variables $x_1,\dots,x_n$, where the 10 exponents are chosen as in
Lemma~\ref{lem 2.2}.
$$
\begin{array}{ll}
  &
\Bigl(\bigl(a_1^{k_1}a_3^{l_1}\bigr)^{m_1}\bigl(a_2^{k_2}a_3^{l_2}\bigr)^{m_2}\Bigr)^{s}
\bigl(a_2^{p}a_3^{q}\bigr)^{t}\vspace*{2mm}\\
= &
\Bigl(\bigl(v_1^{k_1}v_3^{l_1}\bigr)^{m_1}\bigl(v_2^{k_2}v_3^{l_2}\bigr)^{m_2}\Bigr)^{s}
\bigl(v_2^{p}(v_3w)^{q}\bigr)^{t}.
\end{array}
$$

This equation has the solution $(x_1,\dots,x_n)=(g_1,\dots ,g_n)$ in $G$, hence it has a solution
$(h_1,\dots, h_n)$ in $H$.
We set $V_i=v_i(h_1,\dots,h_n)$, $W=w(h_1,\dots,h_n)$. Then we have

$$
\begin{array}{ll}
  &
\Bigl(\bigl(a_1^{k_1}a_3^{l_1}\bigr)^{m_1}\bigl(a_2^{k_2}a_3^{l_2}\bigr)^{m_2}\Bigr)^{s}
\bigl(a_2^{p}a_3^{q}\bigr)^{t}\vspace*{2mm}\\
= &
\Bigl(\bigl(V_1^{k_1}V_3^{l_1}\bigr)^{m_1}\bigl(V_2^{k_2}V_3^{l_2}\bigr)^{m_2}\Bigr)^{s}
\bigl(V_2^{p}(V_3W)^{q}\bigr)^{t}.
\end{array}
$$

Let $U$ be the left side of this equation. Observe that $U\in H$. By Lemma~\ref{lem 2.2}, there exists $\alpha\in \mathbb{Z}$ such that

$$V_1=a_1^{U^{\alpha}},\hspace*{2mm} V_2=a_2^{U^{\alpha}},\hspace*{2mm} V_3=a_3^{U^{\alpha}}, \hspace*{2mm} W=1.\eqno{(11.1)}$$

Since $W=1$, there is a homomorphism $\varphi: G\rightarrow H$, sending $g_j$ to $h_j$, $j=1,\dots,n$. Moreover, the homomorphism $\widehat{(U^{\alpha})^{-1}}\circ \varphi:G\rightarrow H$ is a retraction,
since for $i=1,2,3$, we have
$$
\begin{array}{ll}
\widehat{(U^{\alpha})^{-1}}\circ \varphi(a_i)& =\widehat{(U^{\alpha})^{-1}}\circ \varphi(v_i(g_1,\dots,g_n))\vspace*{2mm}\\
& = \widehat{(U^{\alpha})^{-1}}(v_i(h_1,\dots,h_n))=\widehat{(U^{\alpha})^{-1}}(V_i)\overset{(11.1)}{=}a_i.
\end{array}
$$
\hfill $\Box$

\medskip

The last statement of the following lemma (starting from ``Moreover'') is not needed for Proposition~\ref{prop 2.1}; it will be used later in Section~12.


\begin{lem}\label{lem 2.2} Let $H$ be an acylindrically hyperbolic group without nontrivial normal finite subgroup
and let $a_1,\dots, a_k\in H$  ($k\geqslant 3$) be jointly special and pairwise non-commensurable elements.

Then there exist 10 positive integers
$k_1$, $l_1$, $m_1$, $k_2$, $l_2$, $m_2$, $s$, $p$, $q$, $t$ such that the following holds.

Let $U$ be the left side of the equation
$$
\begin{array}{ll}
  &
\Bigl(\bigl(a_1^{k_1}a_3^{l_1}\bigr)^{m_1}\bigl(a_2^{k_2}a_3^{l_2}\bigr)^{m_2}\Bigr)^{s}
\bigl(a_2^{p}a_3^{q}\bigr)^{t}\vspace*{2mm}\\
= &
\Bigl(\bigl(x_1^{k_1}x_3^{l_1}\bigr)^{m_1}\bigl(x_2^{k_2}x_3^{l_2}\bigr)^{m_2}\Bigr)^{s}
\bigl(x_2^{p}(x_3y_3)^{q}\bigr)^{t}.
\end{array}
$$

Then, for any solution $(x_1,x_2,x_3,y_3)=(b_1,b_2,b_3,c_3)$ of this equation in $H$, there exists an integer number $\alpha$ such that
$$b_1=a_1^{U^{\alpha}},\hspace*{2mm} b_2=a_2^{U^{\alpha}},\hspace*{2mm} b_3=a_3^{U^{\alpha}}, \hspace*{2mm} c_3=1.$$

Moreover, the 10 exponents can be chosen so that, additionally
to the above statement, the elements $U, a_1,\dots,a_k,$ became jointly special and pairwise non-commensurable.
\end{lem}

\medskip

{\it Proof.} First we find an appropriate generating set of $H$.



{\bf Claim 1.} There exists a generating set $Y$ of $H$ such that the following properties are satisfied.

(i) The group $H$ is acylindrically hyperbolic with respect to $Y$.

(ii) The elements $a_1,\dots, a_k$ are special with respect to $Y$.

(iii) $\{\langle a_1\rangle,\dots , \langle a_k\rangle\}\hookrightarrow_h(H,Y)$.

\medskip

{\it Proof.}
Conditions of lemma imply that $E_H(a_j)=\langle a_j\rangle$, $j=1,\dots,k$,
and
$$\{\langle a_1\rangle,\dots ,\langle a_k\rangle\}\hookrightarrow_h H,\eqno{(11.2)}$$
(see~\cite[Theorem 6.8]{DOG}).
Applying Lemma~\ref{best_Y} to this hyperbolic embedding, we obtain a generating set $Y$ of $H$ such that
$H$ is acylindrically hyperbolic with respect to $Y$ and
$\{\langle a_1\rangle,\dots , \langle a_k\rangle\}\hookrightarrow_h(H,Y)$.
Thus, the properties (i) and (iii) are satisfied.

We prove (ii). Property (i) implies that any element of $H$ is either elliptic or loxodromic with respect to $Y$.
For $j=1,\dots ,k$, let $\widehat{d_j}^{Y}$ be the relative metric on $\langle a_j\rangle$ associated with the hyperbolic embedding (11.2). By definition, the space $(\langle a_j\rangle,\widehat{d_j}^{Y})$ is locally finite. Therefore $a_j$ cannot be elliptic with respect to $Y$. Thus, $a_j$ is loxodromic with respect to $Y$ and satisfies $E_G(a_j)=\langle a_j\rangle$. Hence, $a_j$ is special with respect to $Y$.
\hfill $\Box$

\medskip

We use the following

{\bf Notation.} Given $a,b,c,d\in H$, we say that the pair $(a,b)$ is {\it conjugate} to the pair $(c,d)$
if there exists $g\in H$ such that
$g^{-1}ag=c$ and $g^{-1}bg=d$.
In this case we write $(a,b)\sim (c,d)$.

\medskip

Let $k_1$, $l_1$, $m_1$, $k_2$, $l_2$, $m_2$, $s$, $p$, $q$, $t$ be arbitrary 10 positive integers
(we call them exponents)
and let $(b_1,b_2,b_3,c_3)$ be a solution of equation in Lemma~\ref{lem 2.2}:
$$\Bigl(\bigl(a_1^{k_1}a_3^{l_1}\bigr)^{m_1}\bigl(a_2^{k_2}a_3^{l_2}\bigr)^{m_2}\Bigr)^{s}
\bigl(a_2^{p}a_3^{q}\bigr)^{t}
=
\Bigl(\bigl(b_1^{k_1}b_3^{l_1}\bigr)^{m_1}\bigl(b_2^{k_2}b_3^{l_2}\bigr)^{m_2}\Bigr)^{s}
\bigl(b_2^{p}(b_3c_3)^{q}\bigr)^{t}.\eqno{(11.3)}
$$

The diagrams on Fig. 17 reflect the nested structure of the left and the right side of this
equation.

\vspace*{-20mm}
\hspace*{-5mm}
\includegraphics[scale=0.65]{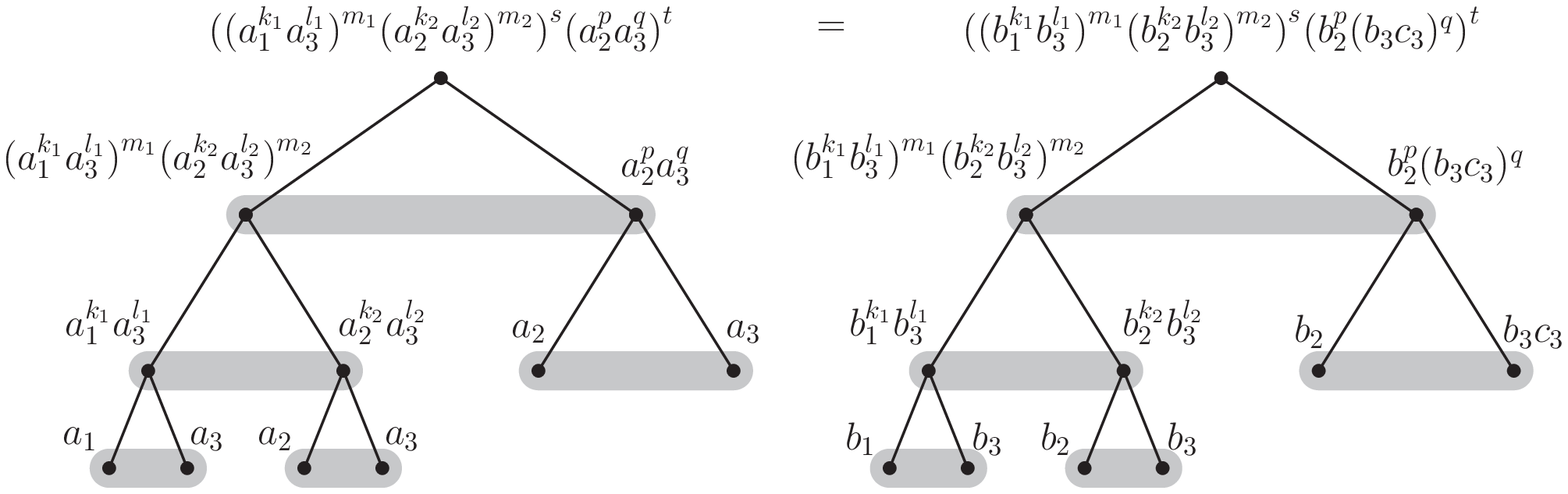}

\vspace*{-120mm}
\hspace*{40mm}
\begin{center}
Fig. 17.
\end{center}

\medskip

\noindent

We explain the structure of forthcoming proof.

$\bullet$ In Step 1 we will choose 10 exponents so that assumptions of Corollary~\ref{prop}
became applicable to 5 pairs of labels of the left diagram (we put them in 5 shadowed regions).

$\bullet$ In Step 2 we will start from the root equation and deduce from Corollary~\ref{prop} consequently
the following formulas:

\begin{enumerate}
\item[(1)] $\Bigl(\bigl(a_1^{k_1}a_3^{l_1}\bigr)^{m_1}\bigl(a_2^{k_2}a_3^{l_2}\bigr)^{m_2},\, \,
a_2^{p}a_3^{q}\Bigr)\sim \Bigl(\bigl(b_1^{k_1}b_3^{l_1}\bigr)^{m_1}\bigl(b_2^{k_2}b_3^{l_2}\bigr)^{m_2},\,\,
b_2^{p}(b_3c_3)^{q}\Bigr)$;
\vspace*{2mm}
\item[(2)] $\bigl(a_1^{k_1}a_3^{l_1},\,\, a_2^{k_2}a_3^{l_2}\bigr)\sim \bigl(b_1^{k_1}b_3^{l_1},\,\, b_2^{k_2}b_3^{l_2}\bigr)$ and $(a_2,a_3)\sim (b_2,b_3c_3)$;
\vspace*{2mm}
\item[(3)] $(a_1,a_3)\sim (b_1,b_3)$ and $(a_2,a_3)\sim (b_2,b_3)$.
\end{enumerate}

$\bullet$ In Step 3 we will analyze these formulas and deduce the statement of lemma.

\medskip

We fix $m\in \mathbb{N}$ such that $\langle a_1^m,a_2^m,a_3^m\rangle$ is a free group of rank 3
(see  Corollary~\ref{free_group}). Let $Y$ be the generating set of $H$ from Claim 1.

\medskip

{\bf Step 1.}
In the following, we will use

-- Corollary~\ref{prop} (to provide perfectness of equations),

-- Corollary~\ref{free_group}
(to construct many loxodromic elements with respect to $Y$),

-- Lemma~\ref{many_non-commensurable} (to construct many non-commensurable elements), and

-- Lemma~\ref{many_special} (to provide $E_G(g)=\langle g\rangle$ for each constructed element $g$).

\medskip

We will also use the principle, that if $u,v\in H$ are non-commensurable, then any element $g\in H$ is non-commensurable with at least one of $u,v$.


\begin{enumerate}
\item[(a)] We choose $k_1,l_1\in m\mathbb{N}$ so that
\begin{enumerate}
\item[(1)] the equation $a_1^{k_1}a_3^{l_1}=x^{k_1}y^{l_1}$ is perfect;
\item[(2)] the element $a_1^{k_1}a_3^{l_1}$ is special with respect to $Y$.
\end{enumerate}

\medskip

In details: By Corollary~\ref{free_group}, the elements $a_1^ia_3^j$ are loxodromic with respect to $Y$
for all sufficiently large $i,j$. By Lemma~\ref{many_special}, $E_G(a_1^ia_3^j)=\langle a_1^ia_3^j\rangle$ for all sufficiently large $i,j$. Thus, $a_1^ia_3^j$ is special with respect to $Y$ for all sufficiently large $i,j$. Then we apply Corollary~\ref{prop} to provide the perfectness.

\medskip

\item[(b)] We choose $k_2,l_2\in m\mathbb{N}$ so that
\begin{enumerate}
\item[(1)] the equation $a_2^{k_2}a_3^{l_2}=x^{k_2}y^{l_2}$ is perfect;
\item[(2)] the element $a_2^{k_2}a_3^{l_2}$ is special with respect to $Y$ and non-commensurable with $a_1^{k_1}a_3^{l_1}$.
\end{enumerate}

\medskip

\item[(c)] We choose $m_1,m_2\in \mathbb{N}$ so that
\begin{enumerate}
\item[(1)] the equation $(a_1^{k_1}a_3^{l_1})^{m_1}(a_2^{k_2}a_3^{l_2})^{m_2}=x^{m_1}y^{m_2}$ is perfect;
\item[(2)] the element $(a_1^{k_1}a_3^{l_1})^{m_1}(a_2^{k_2}a_3^{l_2})^{m_2}$ is special with respect to $Y$.
\end{enumerate}

\medskip

\item[(d)] We choose $p,q\in m\mathbb{N}$ so that
\begin{enumerate}
\item[(1)] the equation $a_2^pa_3^q=x^py^q$ is perfect;
\item[(2)] the element $a_2^pa_3^q$ is special with respect to $Y$ and non-commensurable with $(a_1^{k_1}a_3^{l_1})^{m_1}(a_2^{k_2}a_3^{l_2})^{m_2}$.
\end{enumerate}

\medskip

\item[(e)] We choose $s,t\in \mathbb{N}$ so that
\begin{enumerate}
\item[(1)] the following equation is perfect:
$$
\Bigl(\bigl(a_1^{k_1}a_3^{l_1}\bigr)^{m_1}\bigl(a_2^{k_2}a_3^{l_2}\bigr)^{m_2}\Bigr)^s\Bigl(a_2^pa_3^q\Bigr)^t
=x^sy^t;
$$
\item[(2)] the element on the left side of this equation is special with respect to $Y$
and non-commensurable with elements $a_1,\dots,a_k$.
\end{enumerate}
\end{enumerate}

{\bf Notation.} Let  $A$, $B$, $C$, $D$, $E$ denote the left sides of equations in (a), (b), (c), (d), (e),
respectively.

\medskip

{\bf Step 2.} By (11.3) and (e), there exists $\varepsilon\in \mathbb{Z}$ such that
$$
\bigl(a_1^{k_1}a_3^{l_1}\bigr)^{m_1}\bigl(a_2^{k_2}a_3^{l_2}\bigr)^{m_2}=
\Bigl(\bigl(b_1^{k_1}b_3^{l_1}\bigr)^{m_1}\bigl(b_2^{k_2}b_3^{l_2}\bigr)^{m_2}\Bigr)^{E^{\varepsilon}},\eqno{(11.4)}
$$

$$
a_2^pa_3^q =
\bigl(b_2^p(b_3c_3)^q\bigr)^{E^{\varepsilon}}.\eqno{(11.5)}
$$

\medskip

By (11.4) and (c), there exists $\gamma\in \mathbb{Z}$ such that
$$
a_1^{k_1}a_3^{l_1}=\bigl(b_1^{k_1}b_3^{l_1}\bigr)^{E^{\varepsilon}C^{\gamma}},
\eqno{(11.6)}
$$

$$
a_2^{k_2}a_3^{l_2}=\bigl(b_2^{k_2}b_3^{l_2}\bigr)^{E^{\varepsilon}C^{\gamma}}.
\eqno{(11.7)}
$$

By (11.5) and (d), there exists $\delta\in \mathbb{Z}$ such that
$$
a_2=b_2^{E^{\varepsilon}D^{\delta}},\hspace*{2mm} a_3=(b_3c_3)^{E^{\varepsilon}D^{\delta}},
\eqno{(11.8)}
$$

By (11.6) and (a), there exists $\alpha\in \mathbb{Z}$ such that
$$
a_1=b_1^{E^{\varepsilon}C^{\gamma}A^{\alpha}},\hspace*{2mm} a_3=b_3^{E^{\varepsilon}C^{\gamma}A^{\alpha}},
\eqno{(11.9)}
$$

By (11.7) and (b), there exists $\beta\in \mathbb{Z}$ such that
$$
a_2=b_2^{E^{\varepsilon}C^{\gamma}B^{\beta}},\hspace*{2mm} a_3=b_3^{E^{\varepsilon}C^{\gamma}B^{\beta}},
\eqno{(11.10)}
$$

{\bf Step 3.} From the last equations in (11.9) and (11.10), we deduce that $A^{-\alpha}B^{\beta}$
centralizes $a_3$. We claim that $\alpha=\beta=0$.

Indeed, let $H_1$ be the subgroup of $H$ generated by $a_1^m,a_2^m,a_3^m$. By the choice of $m$,
$H_1$ is free of rank 3.
Since  $A^{-\alpha}B^{\beta}$ lies in $H_1$ and centralises $a_3^{m}$ (which is primitive in $H_1$), the element  $A^{-\alpha}B^{\beta}$ is a power of $a_3^{m}$.
Consider the homomorphism
$$
\varphi: H_1\rightarrow \mathbb{Z}\times \mathbb{Z},\,\,\, a_1^{m}\mapsto (1,0),\,\,\,a_2^{m}\mapsto (0,1),\,\,\,a_3^{m}\mapsto (0,0).
$$
Then $$\varphi(A^{-\alpha}B^{\beta})=
\Bigl(-\frac{k_1}{m}\alpha,\,\,\frac{k_2}{m}\beta\Bigr)=(0,0).$$
Hence $\alpha=\beta=0$.


Using this, we deduce from the first equations in (11.8) and (11.10) that $C^{-\gamma}D^{\delta}$ centralizes $a_2$. We claim that $\gamma=\delta=0$. Indeed, as above we deduce that $C^{-\gamma}D^{\delta}$ is a power of $a_2^m$.
Consider the homomorphism
$$
\psi: H_1\rightarrow \mathbb{Z}\times \mathbb{Z},\,\,\, a_1^{m}\mapsto (1,0),\,\,\,a_2^{m}\mapsto (0,0),\,\,\,a_3^{m}\mapsto (0,1).
$$
Then
$$\psi(C^{-\gamma}D^{\delta})=\Bigl(-\frac{k_1m_1}{m}\gamma,\,\,\,
\frac{q}{m}\delta-\bigl(\frac{l_1m_1}{m}+\frac{l_2m_2}{m}\bigr)\gamma\Bigr)=(0,0).$$
Hence $\gamma=\delta=0$. Then the last equations in (11.8) and (11.9) imply that $c_3=1$.
Moreover, the equations (11.8)-(11.9) imply that $a_1=b_1^{E^{e}}$, $a_2=b_2^{E^{e}}$, $a_3=b_3^{E^{e}}$.

Finally note that the elements $E,a_1,\dots,a_k$ are pairwise non-commensurable and jointly special by (e) and Claim 1~(ii).
\hfill $\Box$

\medskip

\section{Test words in acylindrically hyperbolic groups}

The history of test words in free groups is illuminated in~\cite{Ivanov}. In this paper Ivanov constructed the so-called $C$-test words in free groups. In~\cite{Lee} Lee constructed $C$-test words with some additional property.
In~\cite{MR}, Myasnikov and Roman'kov used Lee's test words to study verbally closed subgroups of free groups.
In this section we construct certain test words in acylindrically hyperbolic groups.

\begin{defn}\label{def_text-word}
{\rm Let $H$ be a group and let $a_1,\dots,a_k$ be some elements of~$H$.
A word $W(x_1,\dots ,x_k)$ is called an
{\it $(a_1,\dots,a_k)$-test word} if for every solution $(b_1,\dots,b_k)$ of the equation
$$
W(a_1,\dots ,a_k)=W(x_1,\dots ,x_k)
$$
in $H$, there exists a number $\alpha\in \mathbb{Z}$ such that
$b_i=a_i^{U^{\alpha}}$ for $i=1,\dots, k$, where $U=W(a_1,\dots ,a_k)$.
}
\end{defn}

\begin{rmk}\label{test_special_words_1}
{\rm Suppose that $H$ is a group and $a_1,a_2,a_3\in H$ are three elements,
satisfying assumptions of Lemma~\ref{lem 2.2} for $k=3$.
This lemma says, in particular, that there exist 10 positive integers $k_1$, $l_1$, $m_1$, $k_2$, $l_2$, $m_2$, $s$, $p$, $q$, $t$ such that
$$
W_3=\Bigl(\bigl(x_1^{k_1}x_3^{l_1}\bigr)^{m_1}\bigl( x_2^{k_2}x_3^{l_2}\bigr)^{m_2}\Bigr)^{s}
(x_2^p(x_3y_3)^q)^t
$$
is an $(a_1,a_2,a_3,1)$-test word.

}
\end{rmk}

{\bf Notation.} We write $\bold{1}^k$ for the tuple $(\underbrace{1,\dots,1}_{k})$.

The aim of this section is to prove the following proposition, which will
be used in Section 13.

\begin{prop}\label{test_special_words}
Let $H$ be an acylindrically hyperbolic group without nontrivial normal finite subgroups
and let $a_1,\dots, a_k\in H$ (where $k\geqslant 3$) be jointly special and pairwise non-commensurable elements.
Then there is an $(a_1,\dots,a_k,\bold{1}^{k-2})$-test word $W_k(x_1,\dots,x_k,y_3,\dots,y_k)$.
\end{prop}

However, to prove this proposition, we need the following its generalisation,
which is simultaneously a generalization of Lemma~\ref{lem 2.2}.

\begin{prop}\label{test_special_words_1}
Let $H$ be an acylindrically hyperbolic group without nontrivial normal finite subgroups
and let $a_1,\dots, a_k\in H$ (where $k\geqslant 3$) be jointly special and pairwise non-commensurable elements.
Then, for any $n=3,\dots,k$, there is an $(a_1,\dots,a_n,\bold{1}^{n-2})$-test word $W_n(x_1,\dots,x_n,y_3,\dots,y_n)$
such that the elements $W_n(a_1,\dots,a_n,\bold{1}^{n-2})$, $a_1,\dots,a_k$ are jointly special and pairwise non-com\-men\-surable.




\end{prop}


\medskip

{\it Proof.}
We fix $k\geqslant 3$ and proceed by induction on $n$. For $n=3$, the statement is valid by Lemma~\ref{lem 2.2}.
Suppose that for some $3\leqslant n<k$, we
have constructed the desired word $W_n=W_n(x_1,\dots,x_n,y_3,\dots,y_n)$.
We show how to construct $W_{n+1}$.

Denote $A=W_n(a_1,\dots,a_n,\bold{1}^{n-2})$.
Since the elements  $A$, $a_1,\dots ,a_k$ are jointly special and pairwise non-commensurable, they satisfy the assumption of Lemma~\ref{lem 2.2}. By this lemma, there
exist positive integers $k_1$, $l_1$, $m_1$, $k_2$, $l_2$, $m_2$, $s$, $p$, $q$, $t$
such that the following holds.

\begin{enumerate}
\item[{\bf (a)}] The word
$$
\frak{M}(X,x_n,x_{n+1},y_{n+1})=\Bigl(\bigl(X^{k_1}x_{n+1}^{l_1}\bigr)^{m_1}\bigl( x_n^{k_2}x_{n+1}^{l_2}\bigr)^{m_2}\Bigr)^{s}
(x_n^p(x_{n+1}y_{n+1})^q)^t
$$
in variables $(X,x_n,x_{n+1},y_{n+1})$ is an $(A,a_n,a_{n+1},1)$-test word and

\medskip

\item[{\bf (b)}] The elements
$\frak{M}(A,a_n,a_{n+1},1)$, $A$, $a_1,\dots,a_k$ are jointly special and pairwise non-commensurable.
\end{enumerate}

We define
$$
W_{n+1}(x_1,\dots,x_{n+1},y_3,\dots,y_{n+1}) = \frak{M}(W_n,x_n,x_{n+1},y_{n+1}).
$$
First we prove that $W_{n+1}$ is an $(a_1,\dots,a_{n+1},\bold{1}^{n-1})$-test word.

Suppose that for some elements $b_1,\dots,b_{n+1}$, $c_3,\dots,c_{n+1}$ in $H$ we have
$$
\begin{array}{ll}
& \Bigl(\bigl(W_n^{k_1}(a_1,\dots,a_n,1,\dots ,1)a_{n+1}^{l_1}\bigr)^{m_1}\bigl( a_n^{k_2}a_{n+1}^{l_2}\bigr)^{m_2}\Bigr)^{s}
(a_n^p(a_{n+1}\cdot 1)^q)^t\vspace*{2mm}\\
= &
\Bigl(\bigl(W_n^{k_1}(b_1,\dots,b_n,c_3,\dots ,c_n)b_{n+1}^{l_1}\bigr)^{m_1}\bigl( b_n^{k_2}b_{n+1}^{l_2}\bigr)^{m_2}\Bigr)^{s}
(b_n^p(b_{n+1}c_{n+1})^q)^t.
\end{array}
$$
Denote $B:=W_n(b_1,\dots,b_n,c_3,\dots ,c_n)$ and write this equation shorter:
$$
\begin{array}{ll}
& \Bigl(\bigl( A^{k_1}a_{n+1}^{l_1}\bigr)^{m_1}\bigl( a_n^{k_2}a_{n+1}^{l_2}\bigr)^{m_2}\Bigr)^{s}
(a_n^p(a_{n+1}\cdot 1)^q)^t\vspace*{2mm}\\
= &
\Bigl(\bigl( B^{k_1}b_{n+1}^{l_1}\bigr)^{m_1}\bigl( b_n^{k_2}b_{n+1}^{l_2}\bigr)^{m_2}\Bigr)^{s}
(b_n^p(b_{n+1}c_{n+1})^q)^t.
\end{array}
$$
Let $U$ be the left side of this equation.

\medskip

Since, by statement {\bf (a)}, $\frak{M}(X,x_n,x_{n+1},y_{n+1})$ is an $(A,a_n,a_{n+1},1)$-test word,
there exists $\alpha\in \mathbb{Z}$ such that
$$
B =A^{U^{\alpha}},
\eqno{(12.1)}
$$

$$
b_n =a_n^{U^{\alpha}},\hspace*{3mm}
b_{n+1}=a_{n+1}^{U^{\alpha}},\hspace*{3mm}\\
c_{n+1}=1;
\eqno{(12.2)}
$$
From (12.1) we deduce
$$
W_n\bigl(b_1^{U^{-\alpha}},\dots,b_n^{U^{-\alpha}},c_3^{U^{-\alpha}},\dots ,c_n^{U^{-\alpha}}\bigr) =W_n(a_1,\dots,a_n,\bold{1}^{n-2})=A.
$$
Since $W_n$ is an $(a_1,\dots,a_n,\bold{1}^{n-2})$-test word, there exists $\beta\in \mathbb{Z}$
such that

$$
\begin{array}{lll}
b_1^{U^{-\alpha}} & \hspace*{-1mm}=\hspace*{-1mm} & a_1^{A^{\beta}},\vspace*{2mm}\\
\vdots & & \vspace*{2mm}\\
b_n^{U^{-\alpha}} & \hspace*{-1mm}=\hspace*{-1mm} & a_n^{A^{\beta}},
\end{array}
\eqno{(12.3)}
$$
and
$$
c_3=\dots =c_n=1. \eqno{(12.4)}
$$
From the first equation in (12.2) and the last equation in (12.3), we deduce that $a_n^{A^{\beta}}=a_n$.
Since $A$ and $a_n$ are jointly special and non-commensurable, we have $\beta=0$.
Then (12.2)-(12.4) imply that
$$
(b_1,\dots,b_{n+1},c_3,\dots,c_{n+1})=(a_1^{U^{\alpha}},\dots, a_{n+1}^{U^{\alpha}},\bold{1}^{n-1}),\\
$$
i.e. the word $W_{n+1}$ is an $(a_1,\dots,a_{n+1},\bold{1}^{n-1})$-test word.

It remains to show that the elements $W_{n+1}(a_1,\dots,a_{n+1},\bold{1}^{n-1})$, $a_1,\dots,a_k$ are jointly special and pairwise non-commensurable. This follows from statement~{\bf (b)} and the fact that
$W_{n+1}(a_1,\dots,a_{n+1},\bold{1}^{n-1})=\frak{M}(A,a_n,a_{n+1},1)$.
\hfill $\Box$

\section{Verbal closedness for finitely generated acylindrically hyperbolic subgroups}

\begin{thm}\label{prop 4.1}
Suppose that $G$ is a finitely presented group and $H$ is a finitely generated acylindrically
hyperbolic subgroup of $G$ such that $H$ does not have nontrivial finite normal subgroups.
Then $H$ is verbally closed in $G$ if and only if $H$ is a retract of $G$.
\end{thm}

\medskip

{\it Proof.} We suppose that $H$ is verbally closed in $G$ and prove that $H$ is a retract of $G$. The converse
statement is obvious. Let $G=\langle g_1,\dots,g_n\,|\, R_1,\dots,R_m\rangle$ and $H=\langle a_1,\dots, a_k\rangle$.
Enlarging $m$ by adding trivial relations and using Proposition~\ref{special_generators}, we may assume that $m\geqslant 1$, $k=m+2$, and that $a_1,\dots, a_k$ are jointly special
and pairwise non-commensurable.
By Proposition~\ref{test_special_words}, there exists an $(a_1,\dots,a_k,\bold{1}^{k-2})$-test word
$W_k(x_1,\dots,x_k, y_3,\dots,y_k)$.

Let $v_i$ and $u_j$ be words in variables $x_1,\dots,x_n$ such that  $a_i=v_i(g_1,\dots,g_n)$ for $i=1,\dots,k$ and $R_j=u_j(g_1,\dots,g_n)$  for $j=1,\dots,m$.
Consider the following equation:
$$
\begin{array}{ll}
& W_k(a_1,\dots,a_k,1,\dots,1)\vspace*{2mm}\\
= & W_k(v_1(x_1,\dots,x_n),\dots,v_k(x_1,\dots,x_n),
u_1(x_1,\dots,x_n),\dots,u_m(x_1,\dots,x_n)).
\end{array}
$$

This equation has the solution $(g_1, \dots ,g_n)$ in $G$, hence it has a solution $(h_1,\dots, h_n)$ in $H$.
We set $V_i=v_i(h_1,\dots,h_n)$, $U_j=u_j(h_1,\dots,h_n)$. Then we have

$$
 W_k(a_1,\dots,a_k,1,\dots,1)= W_k(V_1,\dots,V_k,U_1,\dots,U_m).
$$

Let $U$ be the left side of this equation. Observe that $U\in H$. By Definition~\ref{def_text-word}, there exists $\alpha\in \mathbb{Z}$ such that

$$V_i=a_i^{U^{\alpha}},\hspace*{2mm} i=1,\dots,k,\eqno{(13.1)}$$
$$
U_j=1,\hspace*{2mm} j=1,\dots,m.\eqno{(13.2)}
$$

Because of (13.2), there is a homomorphism $\varphi: G\rightarrow H$, sending $g_i$ to $h_i$, $i=1,\dots,n$. Moreover, the homomorphism $\widehat{(U^{\alpha})^{-1}}\circ \varphi:G\rightarrow H$ is a retraction,
since
$$
\begin{array}{ll}
\widehat{(U^{\alpha})^{-1}}\circ \varphi(a_i)& =\widehat{(U^{\alpha})^{-1}}\circ \varphi(v_i(g_1,\dots,g_n))\vspace*{2mm}\\
& = \widehat{(U^{\alpha})^{-1}}(v_i(h_1,\dots,h_n))=\widehat{(U^{\alpha})^{-1}}(V_i)\overset{(13.1)}{=}a_i.
\end{array}
$$
\hfill $\Box$

\begin{rmk}\label{example}
{\rm In Theorem~\ref{prop 4.1},
the assumption that $H$ does not have nontrivial finite normal subgroups cannot be omitted.
Indeed, consider two copies of the dihedral group $D_4$:
$$
\begin{array}{ll}
A=\langle a,b\,|\, a^4=1, b^2=1, b^{-1}ab=a^{-1} \rangle,\vspace*{2mm}\\
B=\langle c,d\,|\, c^4=1, d^2=1, d^{-1}cd=c^{-1} \rangle.
\end{array}
$$

Let $\varphi:B\rightarrow A$ be the isomorphism sending $c$ to $a$ and $d$ to $b$.
We write $A\underset{a^2=c^2}{\times} B$ for the quotient of the direct product $A\times B$ by the cyclic subgroup
$\langle (a^2,c^2)\rangle$.
We identify $A$ and $B$ with their canonical isomorphic images in this quotient.
Then the elements of $A\underset{a^2=c^2}{\times} B$ can be written as $pq$, where $p\in A$, $q\in B$. If
$p,p_1\in A$ and  $q,q_1\in B$, then $pq=p_1q_1$ in this quotient if and only if $p=p_1$ and $q=q_1$, or
$p_1=pa^2$ and $q_1=qc^2$.
Let $F$ be the free group of rank~2. We set
$$
G=F\times (A\underset{a^2=c^2}{\times} B)=(F\times A)\underset{a^2=c^2}{\times} B
$$
and consider $H=F\times A$ as a subgroup of $G$.
Clearly, $H$ is hyperbolic. Since $H$ is not virtually cyclic, it is acylindrically hyperbolic.

\medskip

{\it Claim.} The following statements hold.
\begin{enumerate}
\item[(a)] $A$ is verbally closed in $A\underset{a^2=c^2}{\times} B$.
\vspace*{1mm}
\item[(b)] $H$ is verbally closed in $G$.
\vspace*{3mm}
\item[(c)] $H$ is not a retract of $G$.
\vspace*{3mm}
\item[(d)] $H$ is not algebraically closed in $G$.
\end{enumerate}
}
\end{rmk}

{\it Proof.} (a) Suppose that an equation $W(x_1,\dots,x_n)=v1$, where $v\in A$,
has a solution $x_1=p_1q_1,\dots ,x_n=p_nq_n$ in $A\underset{a^2=c^2}{\times} B$.
We shall find a solution in~$A$.
Using commutativity, we deduce $W(p_1,\dots,p_n)W(q_1,\dots,q_n)=v1$. Then we have two cases.



{\it Case 1.} $W(p_1,\dots,p_n)=v$ and $W(q_1,\dots,q_n)=1$.

Then $(p_1,\dots,p_n)$ is the desired solution.

\medskip

{\it Case 2.} $W(p_1,\dots,p_n)=va^2$ and $W(q_1,\dots,q_n)=c^2$.

If the exponent sum of $W$ in some letter $x_i$, say in the letter $x_1$, is odd, then $(p_1a^2,\dots ,p_n)$
is the desired solution.
Suppose that the exponent sum of $W$ in any letter $x_i$ is even. Then $W(p_1,\dots,p_n)\in A^2=\{1,a^2\}$, hence $v=1$ or $v=a^2$.
If $v=1$, then $(1,\dots,1)$ is the desired solution, and
if $v=a^2$, then $(\varphi(q_1),\dots,\varphi(q_n))$ is one.

 \medskip

Statement (b) follows from (a) by using the fact that $F$ is a complementary direct summand to $A\underset{a^2=c^2}{\times} B$ in $G$.

\medskip

(c) Suppose that $\psi:G\rightarrow H$ a retraction. Then $[\psi(B),H]=[\psi(B),\psi(H)]=\psi([B,H])=1$. Therefore $\psi(B)\subseteq \langle a^2\rangle$. Since $\langle a^2\rangle=H\cap B$, we have $\psi(B)=\langle a^2\rangle$.
We obtain that $\langle c^2\rangle$ is a retract of $B$. A contradiction.

\medskip

Statement (d) follows from (c) by~\cite[Proposition 1.1]{MR}. This proposition says that if $P$ is a finitely generated subgroup of a finitely presented group $Q$, then $P$ is algebraically closed in $Q$ if and only if $P$ is a retract of $Q$.
\hfill $\Box$


\medskip


\medskip

Our nearest aim is to prove Proposition~\ref{cyclic_and_dyhedral}, which will enable us to assume that $H$ is non-(virtually cyclic) in the forthcoming proofs.
We will prove this proposition with the help of two following lemmas.
The first lemma is simple; we extracted its proof from the proof of~\cite[Lemma 3.1]{MR}.
The second lemma is nontrivial; its proof was given by Klyachko, Mazhuga, and Miroshnichenko in~\cite{KMM}.

\begin{lem}\label{cyclic}
Let $G$ be a group such that its abelianization $G^{ab}$ is finitely generated. Let
$H$ be an infinite cyclic subgroup of $G$. Then $H$ is verbally closed in $G$ if and only it is a retract of $G$.
\end{lem}


{\it Proof.} Suppose that $H=\langle h\rangle$ is an infinite cyclic verbally closed subgroup of $G$.
Let ${\rm Tor}(G^{ab})$ be the subgroup of $G^{ab}$ consisting of all elements of finite order and let $\varphi:G\rightarrow G^{ab}/{\rm Tor}(G^{ab})$ be the canonical homomorphism.

We claim that $H\cap\, [G,G]=1$. Indeed, suppose that $h_1=[g_1,g_2]\dots [g_{2k-1},g_{2k}]$ for some $h_1\in H$ and $g_i\in G$, $i=1,\dots ,2k$. Consider the equation $h_1=[x_1,x_2]\dots [x_{2k-1},x_{2k}]$. Since this equation has
a solution in $G$, it has a solution in $H$. Since $H$ is abelian, we have $h_1=1$.

Thus, $\varphi$ embeds $H$ into $G^{ab}/{\rm Tor}(G^{ab})$. We claim that $\varphi(h)$ is primitive in the free abelian group $G^{ab}/{\rm Tor}(G^{ab})$.
Indeed, otherwise we would have $\varphi(h)=\varphi(g)^t$ for some $g\in G$ and $t\geqslant 2$.
Let $s$ be the order of ${\rm Tor}(G^{ab})$. Then $h^s=g^{st}[g_1,g_2]\dots [g_{2k-1},g_{2k}]$ for some $g_i\in G$, $i=1,\dots ,2k$. Consider the equation $h^s=x^{st}[x_1,x_2]\dots [x_{2k-1},x_{2k}]$.
Since this equation has a solution in $G$, it has a solution in $H$. Then $h^s=h_1^{st}$ for some $h_1\in H$.
Hence $h=h_1^t$ and we have $t=\pm 1$. A contradiction.

Thus, $H$ is embedded into $G^{ab}/{\rm Tor}(G^{ab})$ as a direct summand. Hence $H$ is a retract of $G$.
The converse statement is obvious.
\hfill $\Box$






\begin{lem}\label{dihedergroup} {\rm (see~\cite[Theorem 2]{KMM})}
Suppose that $H$ is an infinite diheder subgroup of a finitely generated group $G$.
Then $H$ is verbally closed in $G$ if and only if $H$ is a retract of $G$.
\end{lem}

\begin{prop}\label{cyclic_and_dyhedral}
Let $H$ be a virtually cyclic subgroup of a finitely generated group $G$.
Suppose that $H$ does not have nontrivial finite normal subgroups.
Then $H$ is verbally closed in $G$ if and only if $H$ is a retract of $G$.
\end{prop}


{\it Proof.} We may assume that $H$ is nontrivial.
It is well-known (see, for example,~\cite[Lemma 2.5]{FJ})
that every virtually cyclic group has a finite-by-cyclic subgroup of index at most 2.
Thus, there exists a subgroup $H_0\leqslant H$ of index at most 2 and a finite normal subgroup
$K\leqslant H_0$ such that $H_0/K$ is cyclic. By assumptions, $H$ cannot be finite.
Therefore $H_0/K\cong \mathbb{Z}$, which implies that $K$ is the largest finite normal subgroup of $H_0$.
Hence $K$ is normal in $H$.
Since $H$ does not have nontrivial finite normal subgroups, we obtain that $K=1$.
Then $H$ is either infinite cyclic or infinite dihedral, and the statement follows from
Lemmas~\ref{cyclic} and~\ref{dihedergroup}.\hfill $\Box$

\medskip

The following lemma says that, in some sense, generic subgroups of relatively hyperbolic groups
are acylindrically hyperbolic.
In this section we use only special cases of this lemma; in Section 14 we use it in full generality.
For terminology concerning relatively hyperbolic groups we refer to~\cite{Osin_0}.


\begin{lem}\label{from rel_hyp_to ax_hyp}
Suppose that $G$ is a group that is relatively hyperbolic with respect to a collection of subgroups $\{P_{\lambda}\}_{\lambda\in \Lambda}$.
Suppose that $H$ is a non-(virtually cyclic) non-parabolic subgroup of $G$.
Then $H$ is acylindrically hyperbolic.

In particular, the following groups are acylindrically hyperbolic:

\begin{enumerate}
\item[(1)] non-(virtually cyclic) groups that are hyperbolic relative to a
collection of proper subgroups (see also~\cite{Osin_3});

\item[(2)] non-(virtually cyclic) subgroups of  hyperbolic groups.
\end{enumerate}
\end{lem}


{\it Proof.}
Let $X$ be a finite relative generating set of $G$ and let $\mathcal{P}=\underset{\lambda\in \Lambda}{\sqcup} P_{\lambda}$.
Then the Cayley graph $\Gamma(G,X\sqcup \mathcal{P})$ is hyperbolic by~\cite[Corollary 2.54]{Osin_0} and the action of $G$ on $\Gamma(G,X\sqcup \mathcal{P})$ is acylindrical by~\cite[Proposition 5.2]{Osin_1}.\break
In particular, $H$ acts acylindrically on $\Gamma(G,X\sqcup \mathcal{P})$.

By~\cite[Lemma 2.9]{Bog_Bux}), a subgroup of a relatively hyperbolic group contains a loxodromic element
if and only if it is infinite and non-parabolic.
Thus, $H$ contains a loxodromic element. Hence $H$ has unbounded orbits acting on $\Gamma(G,X\sqcup \mathcal{P})$.
Then the statement follows from Theorem~\ref{check_elementary}.
\hfill $\Box$

\medskip


\begin{cor}\label{corol_rel_hyp_1}
Let $G$ be a finitely presented group and $H$ be a finitely generated subgroup of~$G$.
Suppose that $H$ is hyperbolic relative to a collection of proper subgroups and
does not have nontrivial finite normal subgroups.
Then $H$ is verbally closed in $G$ if and only if $H$ is a retract of~$G$.
\end{cor}


{\it Proof.} By Proposition~\ref{cyclic_and_dyhedral}, we may assume that $H$ is non-(virtually cyclic).
Then $H$ is acylindrically hyperbolic by Lemma~\ref{from rel_hyp_to ax_hyp}.
Finally the statement follows from Theorem~\ref{prop 4.1}.\hfill $\Box$

\medskip

Corollary~\ref {corollary_hyp_finitely_generated_1} solves Problem 5.2 in~\cite{MR} in the case of finitely generated subgroups. Its proof is the same as the proof of the previous corollary.

\begin{cor}\label{corollary_hyp_finitely_generated_1}
Let $G$ be a hyperbolic group and $H$ be a finitely generated
subgroup of $G$. Suppose that $H$ does not have nontrivial finite normal subgroups.
Then $H$ is verbally closed in $G$ if and only if $H$ is a retract of $G$.
\end{cor}


In the next section, we solve Problem 5.2 in~\cite{MR} in the general case, see Corollary~\ref{solution_Myasnikov_1}.

\newpage

\section{Verbal closedness for equationally noetherian acylindrically hyperbolic subgroups}


\begin{thm}\label{Theorem_2}
Let $G$ be a group and let $H$ be a subgroup~of~$G$ such that $G$ is finitely generated over $H$. Suppose that
$H$ is equationally noetherian, acylindrically hyperbolic,
and does not have nontrivial finite normal subgroups.
Then $H$ is verbally closed in $G$ if and only if $H$ is a retract of~$G$.
\end{thm}


{\it Proof.} We suppose that $H$ is verbally closed in $G$ and prove that $H$ is a retract of $G$. The converse
statement is obvious. Let $\{g_1,\dots,g_n\}$ be a finite set of elements of $G$ which generates $G$ relatively to subgroup $H$, i.e. $G=\langle g_1,\dots,g_n,H\rangle$.
Let $\mathcal{R}$ be the set of all relations of $G$ for this generating set, i.e. $\mathcal{R}$ is the set of all reduced words $r(y_1,\dots,y_n;H)$ over
the alphabet $\{y_1,\dots,y_n\}^{\pm}\sqcup (H\setminus \{1\})$ such that $r(g_1,\dots,g_n;H)=1$ in $G$.
Since $H$ is equationally noetherian, the system of equations
$$
S=\{r(y_1,\dots,y_n;H)=1\,|\, r\in \mathcal{R}\}
$$
over $H$ is equivalent to a finite subsystem
$$
S_0=\{r_i(y_1,\dots,y_n;H)=1\,|\, i=1,\dots,m\}.
$$
Let $h_1,\dots,h_{\ell}\in H$ be all constants of this  subsystem. It is convenient to write $S_0$ in the
following form:
$$
S_0=\{r_i(y_1,\dots,y_n,h_1,\dots,h_{\ell})=1\,|\, i=1,\dots,m\}.
$$


\medskip

Increasing the set $\{ h_1,\dots,h_{\ell}\}$ if necessary, we may assume that
the subgroup $K=\langle h_1,\dots,h_{\ell}\rangle$ of $H$ is noncyclic and contains a special element
(see Lemma~\ref{one_special}). By Proposition~\ref{subgroup_gener},
$K$ can be generated by a finite number (which can be arbitrary large) of pairwise
non-commensurable and jointly special elements, say $K=\langle a_1,\dots,a_k\rangle$,
where $k\geqslant 3$. Extending $S_0$ by repeating some equations (for instance by setting $r_{m+1}=r_m$ and so on) and increasing $k$ if necessary, we may assume that $k=m+2$.

We express $h_1,\dots ,h_{\ell}$ as words in these generators: $h_i=v_i(a_1,\dots,a_k)$, $i=1,\dots,m$.
Then we have
$$
S_0=\{r_i\bigl(y_1,\dots,y_n,v_1(a_1,\dots,a_k),\dots,v_{\ell}(a_1,\dots,a_k)\bigr)=1\,|\, i=1,\dots,m\}.
$$

Let $W_k$ be an $(a_1,\dots,a_k,\bold{1}^{k-2})$-test word for $H$ from Proposition~\ref{test_special_words}.
We consider the equation
$$
\begin{array}{lll}
& W_k(a_1,\dots,a_k, & \bold{1}^{k-2})\vspace*{2mm}\\
= &W_k\bigl(x_1,\dots,x_k,
& r_1\bigl(y_1,\dots,y_n, v_1(x_1,\dots,x_k),\dots ,v_{\ell}(x_1,\dots,x_k)\bigr),\vspace*{2mm}\\
& & \dots,\vspace*{2mm}\\
& & r_m\bigl(y_1,\dots,y_n, v_1(x_1,\dots,x_k),\dots ,v_{\ell}(x_1,\dots,x_k)\bigr)
\bigr).
\end{array}
$$
This equation has a solution in $G$, namely
$$
(x_1,\dots ,x_k,y_1,\dots,y_n)=(a_1,\dots, a_k,g_1,\dots,g_n).
$$
Hence it has a solution in $H$, say
$$
(x_1,\dots ,x_k,y_1,\dots,y_n)=(b_1,\dots, b_k,e_1,\dots,e_n).
$$
Hence we have
$$
\begin{array}{lll}
& W_k(a_1,\dots,a_k, & \bold{1}^{k-2})\vspace*{2mm}\\
= &W_k\bigl(b_1,\dots,b_k,
& r_1\bigl(e_1,\dots,e_n, v_1(b_1,\dots,b_k),\dots ,v_{\ell}(b_1,\dots,b_k)\bigr),\vspace*{2mm}\\
& & \dots ,\vspace*{2mm}\\
& & r_m\bigl(e_1,\dots,e_n, v_1(b_1,\dots,b_k),\dots ,v_{\ell}(b_1,\dots,b_k)\bigr)
\bigr).
\end{array}
$$
Let $U$ be the left side of this equation. Observe that $U\in H$.
Since $W_k$ is an $(a_1,\dots,a_k,\bold{1}^{k-2})$-test word,
there exists $\alpha\in \mathbb{Z}$ such that

$$b_j=a_j^{U^{\alpha}},\hspace*{2mm} j=1,\dots,k,\eqno{(14.1)}$$
$$
r_i\bigl(e_1,\dots,e_n, v_1(b_1,\dots,b_k),\dots ,v_{\ell}(b_1,\dots,b_k)\bigr)=1,\hspace*{2mm} i=1,\dots,m.\eqno{(14.2)}
$$
Applying the conjugation by $U^{-\alpha}$ to (14.2) and using (14.1), we deduce
$$
r_i\bigl(e_1^{U^{-\alpha}},\dots,e_n^{U^{-\alpha}}, v_1(a_1,\dots,a_k),\dots ,v_{\ell}(a_1,\dots,a_k)\bigr)=1,\hspace*{2mm} i=1,\dots,m,
$$
i.e.
$$
r_i\bigl(e_1^{U^{-\alpha}},\dots,e_n^{U^{-\alpha}}, h_1,\dots ,h_{\ell}\bigr)=1,\hspace*{2mm} i=1,\dots,m,
$$
Then
$$(y_1,\dots,y_n)=(e_1^{U^{-\alpha}},\dots ,e_n^{U^{-\alpha}})$$
is a solution in $H$ of the system $S_0$, and hence of $S$. Thus, the map $\varphi:\{g_1,\dots,g_n\}\sqcup \,H\rightarrow H$ sending
$g_i$ to $e_i^{U^{-\alpha}}$, $i=1,\dots,n$, and $h$ to $h$ for any $h\in H$ preserves all relations of $G$.
Therefore it can be extended to a homomorphism $\varphi:G\rightarrow H$.
Since $\varphi|_{H}={\rm id}$, the homomorphism $\varphi:G\rightarrow H$  is a retraction.\hfill $\Box$

\bigskip

We derive two corollaries about relatively hyperbolic groups
with the help of the following remarkable result of Groves and Hull.

\medskip

\begin{thm}\label{Groves_and_Hull} {\rm(see~\cite[Theorem D]{Groves_1})}
Suppose that $G$ is a relatively hyperbolic group with respect to a finite collection of subgroups $\{H_1,\dots ,H_n\}$. Then $G$ is equationally noetherian if and only if each $H_i$ is equationally noetherian.
\end{thm}

\newpage

\begin{cor}\label{noeth_rel_hyp_a}
Let $G$ be a group and let $H$ be a subgroup~of~$G$ such that $G$ is finitely generated over $H$. Suppose that
$H$ is hyperbolic relative to a finite collection of equationally noetherian proper subgroups and does not have nontrivial finite normal subgroups. Then $H$ is verbally closed in $G$ if and only if $H$ is a retract of $G$.
\end{cor}

{\it Proof.}
By Proposition~\ref{cyclic_and_dyhedral}, we may assume that $H$ is non-(virtually cyclic).
Then, by Lemma~\ref{from rel_hyp_to ax_hyp},
$H$ is acylindrically hyperbolic.
Moreover, $H$ is equationally noetherian by the result of Groves and Hull~\cite[Theorem D]{Groves_1}.
Then the statement follows from Theorem~\ref{Theorem_2}.
\hfill $\Box$





\medskip

The proof of the following corollary is similar to that of the previous one;
we give it for completeness.


\begin{cor}\label{general_solution_Myasnikov_1}
Let $G$ be a relatively hyperbolic group with respect to a finite collection of finitely generated
equationally noetherian subgroups. Suppose that $H$ is a non-parabolic subgroup
of $G$ such that $H$ does not have nontrivial finite normal subgroups.
Then $H$ is verbally closed in $G$ if and only if $H$ is a retract of $G$.
\end{cor}

{\it Proof.} It follows from the assumptions that $G$ is finitely generated.
By Proposition~\ref{cyclic_and_dyhedral}, we may assume that $H$ is non-(virtually cyclic).
Then, by Lemma~\ref{from rel_hyp_to ax_hyp}, $H$ is acylindrically hyperbolic.
By the result of Groves and Hull~\cite[Theorem D]{Groves_1}, $G$ is equationally noetherian.
Any subgroup of an equationally noetherian group is equationally noetherian.
Therefore $H$ is equationally noetherian.
Then the statement follows from Theorem~\ref{Theorem_2}.
\hfill $\Box$

\medskip

The following corollary follows directly from the previous one.
Indeed, every hyperbolic group is relatively hyperbolic with respect to the trivial subgroup.

\begin{cor}\label{solution_Myasnikov_1}
{\rm (Solution to Problem 5.2 in~\cite{MR})}\\
Let $G$ be a hyperbolic group and $H$ be a subgroup of $G$.
Suppose that $H$ does not have nontrivial finite normal subgroups.
Then $H$ is verbally closed in $G$ if and only if $H$ is a retract of $G$.
\end{cor}


\begin{rmk}\label{easier_proof}
{\rm This corollary can be proved without reference to Corollary~\ref{general_solution_Myasnikov_1}.
A direct proof can be obtained from the above proof if we
recall
the result of Reinfeldt and Weidmann~\cite[Corollary 6.13]{RW} that all hyperbolic groups
are equationally noetherian.
}
\end{rmk}

\newpage

\section{Appendix}

The aim of this appendix is to recall a general concept of {\it closedness} of a structure in a class of structures
with respect to a class of first-order formulas.
Certainly, this general concept, coming from the model theory, is very helpful to understand relations between different notions and theorems.
Moreover, one can hope that techniques developed for certain classes of groups and
classes of formulas can be adopted for other ones.
Starting from a very early history, we will formulate a few distinguished results in this area,
not pretending to list all of them.


In~1943, B.H. Neumann~\cite{Neumann_0} considered systems of equations over arbitrary groups and, taking in mind field theory, introduced for groups such notions as adjoining of solutions,  algebraic and transcendent extensions. He gave a simple (very general) criterion for a system $S$ of equations with coefficients in a group $G$ to be solvable in some overgroup of $G$.
Using it, Neumann showed that if at least one such overgroup exists, then
all of them (up to isomorphism) are overgroups of certain quotients of a distinguished overgroup $G_S$ in which $S$ is solvable.
Inspired by the paper~\cite{Neumann_0}, W.R.~Scott~\cite{Scott} introduced the following notion of algebraically closed groups:

\begin{defn}\label{Def_Scott}
{\rm A group $G$ is called {\it algebraically closed} (in the class of all groups) if any finite system of equations and inequalities with coefficients in~$G$,
$$
\begin{array}{lll}
u_i(x_1,\dots ,x_k;G) & =1 &  (i=1,\dots ,n)\vspace*{2mm}\\
v_j(x_1,\dots ,x_k;G) &\neq 1 & (j=1,\dots ,m),
\end{array}\eqno{(\star)}
$$
which can be solved in some overgroup of $G$ can also be solved in $G$ itself.
}
\end{defn}


\begin{rmk}\label{remark_Neumann}
{\rm According to the point of view of the modern model theory, such groups may be named
{\it existentially closed}, and the groups which satisfy Definition~\ref{Def_Scott} with only equations permitted in $(\star)$ may be named algebraically closed.
Fortunately, for non-trivial groups, two variants of Definition~\ref{Def_Scott} (with or without inequalities)
are equivalent.
This surprising fact was proved by B.H.~Neumann in~\cite{Neumann_1}.

\medskip

We list some interesting properties of algebraically closed groups. Some other interesting results (including those
on elementary equivalence) are listed
in the book of Hodges~\cite[Appendix A.4]{Hodges}.

}
\end{rmk}
\begin{enumerate}
\item[(1)] Every (countable) group can be embedded in an algebraically closed (countable) group~\cite{Scott}.
There exists $2^{\aleph_0}$ mutually non-isomorphic countable algebraically closed groups~\cite{Neumann_2}.

\medskip

\item[(2)] Any algebraically closed group is simple~\cite{Neumann_1}, it cannot be finitely generated~\cite{Neumann_2}, and it does not have an infinite recursive presentation~\cite{LS}.

\medskip

\item[(3)] A finitely generated group has a solvable word problem if and only if it can be embedded in every
algebraically closed group.
The ``only if part'' was proved by Neumann in~\cite{Neumann_2} and the ``if part'' was proved by
Macintyre in~\cite{Macintyre}.
\end{enumerate}

\medskip

The following very general definition of closedness goes back to B.H.~Neumann~\cite{Neumann_2}.
It works for all algebraic structures, not only for groups.
For simplicity, we give here a weaker version of his definition.

Recall that a structure is called {\it algebraic} if its signature does not contain relation symbols
(see~\cite{Hodges} for other model-theory terminology used below). The class of groups can be considered as
a subclass of structures with the signature consisting of $\cdot$ (multiplication symbol), $^{-1}$ (inversion symbol), and $1$ (the identity symbol).

\medskip

\begin{defn}
{\rm Let $\frak{A}$ be some class of algebraic structures with the same signature $L$.
We associate with $L$ the language $\mathcal{L}$ consisting of the elements of $L$, the logical symbols $=$, $\neg$, $\bigwedge$, $\bigvee$, $\exists$, $\forall$, the variables $x_1,x_2,\dots$, and the punctuation signs.

We use the following notation. Suppose that $\phi$ is a first-order formula in the language~$\mathcal{L}$.
A variable $x$ in $\phi$ is called {\it free} in $\phi$ if neither $\forall x$ nor $\exists x$ occur in~$\phi$.
We denote $\phi$ by $\phi(x_1,\dots,x_n)$ if and only if $x_1,\dots,x_n$ are all free variables of~$\phi$.

Let $A$ be a structure from the class $\frak{A}$ and let $\frak{S}$ be a set of formulas in the language $\mathcal{L}$.
The structure $A$ is called {\it $\frak{S}$-closed in $\frak{A}$} if for any formula $\phi(x_1,\dots,x_n)\in \frak{S}$, any elements $a_1,\dots,a_n\in A$,
and any structure $B\in \frak{A}$ containing $A$
if $\phi(a_1,\dots,a_n)$ holds in $B$, then it holds already in $A$.
We write this as follows:
$$
B\models \phi(a_1,\dots,a_n) \Rightarrow A\models \phi(a_1,\dots,a_n).
$$

In the case, where $H$ is a substructure of a structure $G$ and
$H$ is $\frak{S}$-closed in the class $\{H,G\}$,
we simplify notation by saying that {\it $H$ is $\frak{S}$-closed in $G$}.

}
\end{defn}

\medskip

\noindent
{\bf Notations.}
Let $\mathcal{G}$ be the class of all groups and $\Phi$ the set of all finite formulas
in the language $\mathcal{L}$ associated with the signature of groups. We define the following subsets of $\Phi$:

\begin{enumerate}
\item[$\bullet$] $\Phi_0$ is the subset of $\Phi$ consisting of all formulas without free variables
(such formulas are called {\it sentences}).

\item[$\bullet$] ${\exists\!\!\!\!\exists}$\, is the
subset of $\Phi$ consisting of all {\it existential} formulas,
i.e. of the formulas which have the form $\exists\, x_1\dots \exists\, x_n \,\phi(x_1,\dots,x_n,x_{n+1},\dots,x_k)$,
where $\phi\in \Phi$ is a quantifier free formula.

\item[$\bullet$] ${\exists\!\!\!\!\exists}^{+}$ is
the subset of $\Phi$ consisting of all {\it positive existential} formulas, i.e.
of the existential formulas without the negation symbol.

\item[$\bullet$] $\mathcal{V}$ is the subset of $\Phi$ defined as follows:
$$
\mathcal{V}=\underset{n\in \mathbb{N}}{\cup}\, \{\exists\, x_1\dots \exists\, x_n: w(x_1,\dots,x_n)=x_{n+1}\,|\, w(x_1,\dots,x_n)\in F(x_1,\dots,x_n)\}.
$$
Here $F(x_1,\dots,x_n)$ is the free group with basis $x_1,\dots,x_n$.
\end{enumerate}

\medskip

Recall that the {\it elementary theory} of a group $G$, denoted by~${\rm \bf Th}(G)$,
consists of all sentences, which hold in $G$.
The {\it existential theory} of a group $G$, denoted by~${\rm \bf Th_{\exists}}(G)$,
consists of all existential sentences which hold in $G$.


\begin{rmk}\label{remarks_closed} \hspace*{-5mm}.

{\rm
\begin{enumerate}
\item[(1)] Suppose that $H$ is a subgroup of a group $G$. Then
$H$ is algebraically closed (respectively, verbally closed) in $G$ in the sense of Definition~\ref{alg,verb,retr}
if and only if $H$ is ${\exists\!\!\!\!\exists}^{+}$-closed (respectively, $\mathcal{V}$-closed) in $G$.

\medskip

\item[(2)]
A nontrivial group $G$ is algebraically closed in the sense of Definition~\ref{Def_Scott}
if and only if $G$ is ${\exists\!\!\!\!\exists}^{+}$-closed in the class $\mathcal{G}$
(see Remark~\ref{remark_Neumann}).

\medskip

\item[(3)]
By analogy to algebraic closedness, a group $G$ is called {\it verbally closed} if $G$ is $\mathcal{V}$-closed in
the class $\mathcal{G}$.
Note that for any nontrivial reduced word $w(x_1,\dots,x_n)$ in variables $x_1,\dots ,x_n$ and any element $g\in G$, the equation $w(x_1,\dots,x_n)=g$
is solvable in the amalgamated product
$$
G\underset{g=w(x_1,\dots,x_n)}{\ast} \langle x_1,\dots,x_n\,|\, w^r=1\rangle,
$$
where $r\in \mathbb{N}\cup \{\infty\}$ is the order of $g$ and $w^{\infty}$ denotes $1$.
Therefore if $G$ is verbally closed, then $G$ is {\it verbally complete}, i.e. each equation
$w(x_1,\dots,x_n)=g$, where $w$ and $g$ as above, is solvable in $G$.

Mikhajlov and Ol'shanskii~\cite{MO} proved that there exist finitely generated verbally complete groups.
This contrasts to the fact that any algebraically closed group cannot be finitely generated.
Osin~\cite[Corollary 1.6]{Osin_4} proved that there exists an uncountable set of pairwise non-isomorphic 2-generated verbally complete groups.

\medskip

\item[(4)]






Recall that a subgroup $H$ of a group $G$ is called
{\it elementarily embedded} in $G$ if for any formula $\phi(x_1,\dots ,x_n)\in \Phi$ and any
elements $h_1,\dots,h_n\in H$ we have
$$
H\models \phi(h_1,\dots ,h_n) \Leftrightarrow G\models \phi(h_1,\dots ,h_n)
$$

The following claim follows from the fact that $\neg\Phi=\Phi$ and $\Phi_0\subseteq \Phi$.

\medskip

{\it Claim.} Let $H$ be a subgroup of a group $G$. Then the following holds:

\vspace*{0.5mm}

\begin{enumerate}
\item[(a)] $H$ is elementarily embedded in $G$ if and only if it is $\Phi$-closed in~$G$.

\vspace*{0.5mm}
\item[(b)] If $H$ is elementarily embedded in $G$, then ${\rm \bf Th}(H)={\rm \bf Th}(G)$.
\end{enumerate}
\medskip

In the course of solving of Tarski's problem Sela and also Kharlampovich and Myasnikov
proved the following theorem.

\medskip

\noindent
{\bf Theorem A.} {\rm (see ~\cite[Theorem 4]{Sela_6} and~\cite[a claim in Introduction]{KM_1})}
{\it For $n\geqslant k\geqslant 2$, the subgroup $F(x_1,\dots, x_k)$ is elementarily embedded (equivalently $\Phi$-closed) in $F(x_1,\dots,x_n)$.
}

\medskip

Developing Sela's technique, Perin proved the following theorem
(compare with our Corollary~\ref{solution_Myasnikov}).

\medskip

\noindent
{\bf Theorem B.}\label{Perin} {\rm (see~\cite[Theorem 1.2]{Perin_2} and~\cite{Perin_3})}
{\it Let $H$ be a subgroup of a torsion-free hyperbolic group $G$.
If $H$ is elementarily embedded (equivalently $\Phi$-closed) in $G$,
then $G$ admits the structure of a hyperbolic tower over $H$ (in particular $H$ is a retract of~$G$).
}



\medskip

\item[(5)] A subgroup $H$ of a group $G$ is called {\it existentially closed in $G$} if
$H$ is ${\exists\!\!\!\!\exists}$\,-closed in $G$. Note that if $H$ is existentially closed in $G$,
then ${\rm \bf Th}_{\exists}(H)={\rm \bf Th}_{\exists}(G)$.

Existentially closed subgroups in finitely generated free nilpotent groups
were characterized by Roman'kov and Khisamiev in~\cite{RK2}.
A lot of information on existentially closed groups in specific classes of groups is contained in the
book of Higman and Scott~\cite{Higman_Scott} and in the survey of Leinen~\cite{Leinen}.
These classes include the class of all groups, some classes of nilpotent and solvable groups, and the class of locally finite groups.
On the other hand, not too much is known on existentially
closed subgroups of groups with a kind of negative curvature.


\medskip




Theorem C below describes existentionally closed subgroups of tor\-sion-free hyperbolic groups.
This theorem follows directly from the following known results.

$\bullet$ We use a terminology of~\cite{MR}.  If $H$ is a subgroup of $G$, the expression $H\leqslant G$ is called an {\it extension} of $H$ to $G$. An extension $H\leqslant G$ is called {\it discriminating} if for any finite subset $K\subseteq G$ there exists a retraction $\varphi:G\rightarrow H$ which is injective on $K$.
Similarly, we call an extension $H\leqslant G$ {\it separating} if for any nontrivial element $g\in G$
there exists a retraction $\varphi:G\rightarrow H$ such that $\varphi(g)\neq 1$.
These two notions coincide in the case where the group $G$ is torsion-free hyperbolic and
$H$ is a non-cyclic subgroup of~$G$ (see~\cite[Theorem C1 and Section 1.3]{BMR}).

$\bullet$ In~\cite[Proposition 2.3]{MR} Myasnikov and Roman'kov proved the following general statement:
Let $H$ be a subgroup of a group $G$. Suppose that $G$ is finitely generated relative to $H$ and
$H$ is equationally noetherian. Then $H$ is existentially closed in $G$ if and only if the extension
$H\leqslant G$ is discriminating.


$\bullet$ Hyperbolic groups are equationally noetherian.
Sela proved this in the torsion-free case~\cite[Theorem 1.22]{Sela_2};
the general case was considered by Reinfeldt and Weidmann~\cite[Corolary 6.13]{RW}.

\medskip

\noindent
{\bf Theorem C.} {\rm (follows directly from~\cite{MR,BMR,Sela_2})}
{\it
Let $H$ be a subgroup of a torsion-free hyperbolic group $G$.
The following statements are equivalent:

\begin{enumerate}
\item[(a)] $H$ is existentially closed in $G$.

\item[(b)] For any finite subset $K\subseteq G$, there exists a retraction $\phi: G\rightarrow H$ that is injective on $K$.

\item[(c)] For any nontrivial $g\in G$, there exists a retraction $\phi: G\rightarrow H$ such that $\phi(g)\neq 1$.
\end{enumerate}
}
\end{enumerate}
}
\end{rmk}


\begin{thebibliography}{AA}

\bibitem{A} S.I. Adian, {\it The Burnside problem and identities in groups}, v. {\bf 95} of
``Ergebnisse der Mathematik und ihrer Grenzgebiete'' [Results in Mathematics and Related Areas],
Springer-Verlag, Berlin, 1979. (Translated from Russian: The Burnside problem and identities in groups,
Nauka, Moscow, 1975.)

\bibitem{BMR} G. Baumslag, A. Myasnikov, and V. Remeslennikov, {\it Algebraic geometry over groups~I. Algebraic
sets and ideal theory}, J. Algebra, {\bf 219}, Issue 1 (1999), 16-79.

\bibitem{BMRom} G. Baumslag, A. Myasnikov, and V. Roman'kov,
{\it Two theorems about equationally noetherian groups}, J. Algebra, {\bf 194}, Issie 2 (1997), 654-664.

\bibitem{BL} B. Baumslag, F. Levin, {\it Algebraically closed torsion-free nilpotent groups of class 2},
Communications in Algebra, {\bf 4}, no. 6 (1976), 533-560.

\bibitem{BF} M. Bestvina, K. Fujiwara, {\it Bounded cohomology of subgroups of mapping class groups},
Geom. Topol., {\bf 6} (2002), 69-89.

\bibitem{Bog_1} O. Bogopolski, {\it A periodicity theorem for acylindrically hyperbolic groups},
ArXiv, 2019. Available at
https://arxiv.org/pdf/1805.05941v2.pdf

\bibitem{Bog_Bux} O. Bogopolski, K.-U. Bux, {\it From local to global conjugacy of subgroups
of relatively hyperbolic groups}, Intern. J. of Algebra and Computation, {\bf 27}, no. 3 (2017), 299-314.

\bibitem{BG} O.V. Bogopolskii, V.N. Gerasimov, {\it Finite subgroups of hyperbolic groups}, Algebra and Logic,
{\bf 34} (1995), 343-345.

\bibitem{BH} M. Bridson, A. Haeffliger, {\it Metric spaces of non-positive curvature},
 Springer, 1999.

\bibitem{Bryant} R. Bryant, {\it The verbal topology of a group}, J. Algebra, {\bf 48} (1977), 340-346.

\bibitem{Bowditch} B. Bowditch, {\it Tight geodesics in the curve complex}, Invent. Math.,
{\bf 171}, no. 2 (2008), 281-300.




\bibitem{CDP} M. Coornaert, T. Delzant, A. Papadopoulos, {\it Geometrie et th\'{e}orie des groupes. Les groupes hyperboliques de Gromov.} Lecture Notes in Mathematics, {\bf 1441}. Springer-Verlag, Berlin, 1990. x+165 pp.

\bibitem{DOG} F. Dahmani, V. Guirardel, D. Osin, {\it Hyperbolically embedded subgroups and rotating families in groups acting on hyperbolic spaces}, Memoirs Amer. Math. Soc., v. {\bf 245} (2017), no. 1156, v+152 pp.





\bibitem{FJ} F.T. Farrell and L.E. Jones, {\it The lower algebraic K-theory of virtually infinite cyclic groups}, K-theory, {\bf 9}, no. 1 (1995), 13-30.


\bibitem{FW} N.J. Fine and H.S. Wilf, {\it Uniqueness theorem for periodic functions},
PAMS, {\bf 16} (1965), 109--114.



\bibitem{Gromov} M. Gromov, {\it Hyperbolic groups}, Essays in Group Theory, MSRI Series, Vol. {\bf 8},
(S.M. Gersten, ed.), Springer, 1987, 75-263.

\bibitem{Groves_1} D. Groves and M. Hull, {\it Homomorphisms to acylindrically hyperbolic groups I:
Equationally noetherian groups and families}, 2017, arXiv: 1704.03491v1

\bibitem{Groves_2} D. Groves and H. Wilton, {\it Enumerating limit groups}, Groups, Geometry and Dynamics,
{\bf 3}, no. 3 (2009), 389-399.

\bibitem{Guba} V. Guba, {\it Equivalence of infinite systems of equations in free groups and semigroups to finite subsystems}, Mat. Zametki, {\bf 40}, no. 3 (1986), 321-324.

\bibitem{Gupta_Rom} Ch.K. Gupta and N.S. Romanovskii, {\it The property of being equational Noetherian of some
solvuble groups}, Algebra and Logic, {\bf 46}, no. 1 (2007), 28-36.

\bibitem{Ham} U. Hamenst{\"a}dt, {\it Bounded cohomology and isometry groups of hyperbolic spaces},
J. Eur. Math. Soc., {\bf 10}, no. 2 (2008), 315-349.

\bibitem{Higman_Scott} G. Higman and E. Scott, {\it Existentially closed groups}, London Math. Soc. Monogr.,
Ser. 3, Clarendon Press, Oxford, 1988.

\bibitem{Hodges} W. Hodges, {\it Model Theory}, Cambridge Univedrsity Press, Cambridge, 1993.




\bibitem{HO} M. Hull, D. Osin, {\it Induced quasi-cocycles on groups with hyperbolically embedded subgroups},
Algebraic $\&$ Geometric Topology, {\bf 13}, no. 5 (2013), 2635-2665.

\bibitem{Ivanov} S.V. Ivanov, {\it On certain elements of free groups}, Communications in Algebra, {\bf 204},    (1998), 394-405.


\bibitem{JA} Eric Jaligot and Abderezak Ould Houcine, {\it Existentially closed CSA-groups},
J. Algebra, {\bf 280} (2004), 772-796.

\bibitem{JS} E. Jaligot and Z. Sela, {\it Makanin -- Razborov diagrams over free products}, Illinois J. Math.,
{\bf 54}, no. 1 (2010), 19-68.


\bibitem{KM_0} O. Kharlampovich and A. Myasnikov, {\it Irreducible affine varieties over a free group.
I. Irreducibility of quadratic equations and Nullstellensatz}, J. Algebra, {\bf 200}, Issue 2 (1998), 472-516.

\bibitem{KM_1} O. Kharlampovich and A. Myasnikov, {\it Elementary theory of free nonabelian groups},
J. Algebra, {\bf 302}, Issue 2 (2006), 451-552.




\bibitem{KM} A.A. Klyachko and A.M. Mazhuga, {\it Verbally closed virtually free subgroups}, Sb. Math.,
v. {\bf 209} (6) (2018), 850-856.

\bibitem{KMM} A.A. Klyachko, A.M. Mazhuga and V. Yu. Miroshnichenko, {\it Virtually free
finite-normal-subgroup-free groups are strongly verbally closed},
J. Algebra, {\bf 510} (2018), 319-330.


\bibitem{KT} Anton Klyachko and Andreas Thom, {\it New topological methods to solve equations over groups},
Algebr. Geom. Topol., {\bf 17} (1) (2017), 331-353.


\bibitem{Lee} D. Lee, {\it On certain $C$-test words for free groups}, J. Algebra, {\bf 247} (2002), 509-540.

\bibitem{Leinen} F. Leinen, {\it Existentially closed groups in specific classes}. In: {\it Finite and locally finite groups} (Istanbul, 1994), NATO Adv., Sci. Inst. Ser. C. Math. Phys. Sci. {\bf 471}, Kluwer Akademic Publishers, Dordrecht (1995), 285-326.


\bibitem{LS} R.C. Lyndon, P.E. Schupp, {\it Combinatorial group theory}, Springer Verlag, 1977.

\bibitem{Macintyre} A. Macintyre, {\it On algebraically closed groups}, Ann of Math., {\bf 96} (1972), 53-97.

\bibitem{Mazhuga_1} A.M. Mazhuga, {\it On free decompositions of verbally closed subgroups of free products of finite groups}, J. Group Theory, {\bf 20}, no. 5 (2017), 971-986.

\bibitem{Mazhuga_2} A.M. Mazhuga, {\it Strongly verbally closed groups}, J. Algebra, {\bf 493} (2018), 171-184.

\bibitem{Mazhuga_3} A.M. Mazhuga, {\it Free products of groups are strongly verbally closed}, ArXiv, 2018.
Available at https://arxiv.org/abs/1803.10634.


\bibitem{MO} K.V. Mikhajlovskii, A.Yu. Ol'shanskii, {\it Some constructions related to hyperbolic groups}. In
{\it Geometry and Cohomology in Group Theory} (Durham, 1994), London Math. Soc. Lecture Note Ser, {\bf 252},
Cambridge Univ. Press, Cambridge, 1998, pp. 263-290.

\bibitem{Myas_Remesl} A. Myasnikov, V. Remeslennikov, {\it Algebraic geometry over groups. II. Logical foundations}, J. Algebra, {\bf 234}, Issue 1 (2000), 225-276.

\bibitem{MR} A.G. Myasnikov, V. Roman'kov, {\it Verbally closed subgroups of free groups},
J. of Group Theory, {\bf 17}, no. 1 (2014), 29-40.

\bibitem{Neumann_0} B.H. Neumann, Adjunction of elements to groups, J. London Math. Soc., {\bf 18} (1943),
pp. 4-11.

\bibitem{Neumann_1} B.H. Neumann, {\it A note on algebraically closed groups},
J. London Math. Soc., {\bf 27} (1952), 227-242.

\bibitem{Neumann_2} B.H. Neumann, {\it The isomorphism problem for algebraically closed groups},
In.: Word Problems, Amsterdam: North Holland (1973), pp. 553-562.

\bibitem{Osin_0} D. Osin, {\it Relatively hyperbolic groups: Intrinsic geometry, algebraic properties, and algorithmic problems}, Memoirs Amer. Math. Soc., v. {\bf 179} (2006), no. 843.

\bibitem{Osin_1} D. Osin, {\it Acylindrically hyperbolic groups},
Trans. Amer. Math. Soc., v. {\bf 368} (2016), 851-888.

\bibitem{Osin_3} D. Osin, {\it Groups acting acylindrically on hyperbolic spaces}, ArXiv, 2017.
Available at
https://arxiv.org/pdf/1712.00814.pdf

\bibitem{Osin_4} Denis Osin, {\it Small cancellations over relatively hyperbolic groups and embedding theorems},
Annals of Math., {\bf 72} (2010), 1-39.

\bibitem{OT} Denis Osin and Andreas Thom, {\it Normal generation and $\ell^2$-Betti numbers of groups},
Mathematische Annalen, 2013, 1331-1347.

\bibitem{AOH} Abderezak Ould Houcine, {\it Limit groups of equationally noetherian groups.}
In Geometric group theory, Trends in Math., pages 103-119, Birkh{\"a}user, Basel, 2007.


\bibitem{Perin_2} Chlo\'{e} Perin, {\it Elementary embeddings in torsion-free hyperbolic groups}, Ann. Scient.
\'{E}c. Norm. Sup., {\bf 44} (4) (2011), 631-681.

\bibitem{Perin_3} Chlo\'{e} Perin, {\it Erratum to: Elementary embeddings in torsion-free hyperbolic groups},
Ann. Scient. \'{E}c. Norm. Sup., {\bf 46} (4) (2013), 851-856.


\bibitem{RW} C. Reinfeldt and R. Weidmann, {\it Makanin -- Razborov diagrams for hyperbolic groups}, 2014.

\bibitem{Roman'kov} Vitali\u{i} Roman'kov, {\it Equations over groups}, Groups Complexity Cryptology,
{\bf 4} (2) (2012), 191-239.

\bibitem{RK1} V.A. Roman'kov and N.G. Khisamiev, {\it Verbally and existentially closed subgroups of free
nilpotent groups},
Algebra and Logic, {\bf 52}, no. 4 (2013), 336-351.

\bibitem{RK2} V.A. Roman'kov and N.G. Khisamiev, {\it Existentially closed subgroups of free
nilpotent groups}, Algebra and Logic, {\bf 53}, no. 1 (2014), 29-38.

\bibitem{RKK} V.A. Roman'kov, N.G. Khisamiev, and A.A. Konyrkhanova
{\it Algebraically and verbally closed subgroups and retracts of finitely generated nilpotent groups},
Siberian Math. J., {\bf 58} (3) (2017), 536-545.

\bibitem{Rom_1} N.S. Romanovskii, {\it Equational Noetherianity of rigid solvable groups},
Algebra and Logic, {\bf 48}, no. 2 (2009), 147-160.

\bibitem{Scott} W.R. Scott, {\it Algebraically closed groups}, Proc. Amer. Math. Soc., {\bf 2} (1951), 118-121.

\bibitem{Sela_1} Z. Sela, {\it Diophantine geometry over groups. I. Makanin -- Razborov diagrams}, Publ. Math.
Inst. Hautes ${\rm \acute{E}}$tudes Sci., {\bf 93} (2001), 31-105.

\bibitem{Sela_6} Z. Sela, {\it Diophantine geometry over groups. VI. The elementary theory of free groups}, Geom. Funct. Anal., {\bf 16} (2006), 707-730.

\bibitem{Sela_2} Z. Sela, {\it Diophantine geometry over groups. VII:
The elementary theory of a hyperbolic group.}, Proc. Lond. Math. Soc., {\bf 99}, no. 1 (2009), 217-273.

\bibitem{Sela_3} Z. Sela, {\it Diophantine geometry over groups. X:
The elementary theory of free products.}, arXiv:1012.0044, 2010.

\bibitem{Sisto} A. Sisto, {\it Contracting elements and random walks}, J. Reine Angew. Math., {\bf 742} (2018), 79-114.



\end{thebibliography}
\end{document}

\section{Test words}

S. Ivanov introduced $C$-test words $w_n(x_1,\dots ,x_m)$ in free groups $F_n$ for every $m,n\geqslant 2$.

We fix $m\geqslant 2$ until the end of the paper.

\begin{defn}
{\rm A word $w(x_1,\dots ,x_m)$ is called a {\it $C$-test} word in a group $G$ if for every two tuples
$(g_1,\dots,g_m)$, $(g_1',\dots,g_m')$ of elements of $G$ the following holds:
If
$$
w_n(g_1,\dots ,g_m)= w_n(g_1',\dots ,g_m')\neq 1,
$$
then there exists $u\in G$ such that $u^{-1}g_iu=g_i'$, $i=1,\dots,m$.}
\end{defn}

We can interpret the $C$-test words as follows.

For every word $w(x_1,\dots,x_m)$, we consider the map
$$
\begin{array}{ll}
\Phi_w: & \underbrace{G\times \dots \times G}_{m}\rightarrow G,\\
& (g_1,\dots,g_m)\mapsto w(g_1,\dots,g_m).
\end{array}
$$

We say that $\Phi_w$ is {\it conjugacy-injective at} $(g_1,\dots ,g_m)$ if
for any nonconjugate tuple $(g_1',\dots ,g_m')$ the element $\Phi_w(g_1,\dots,g_m)$ is not conjugate
to $\Phi_w(g_1',\dots,g_m')$.

Then the above definition says that $\Phi_{w}$ is conjugacy injective at each point outside of $\Phi_{w}^{-1}(1)$.

The test words are complicated even in the case of free group. The idea of the following definition is to use a larger set of (maybe simpler)
words.

\begin{defn}
{\rm
A (possibly infinite) set of words $u_i=u_i(x_1,\dots ,x_m)$, $i\in I$, is called a {\it $C$-test set} in $G$ if for every tuple $g=(g_1,\dots,g_m)$, there exists $i\in I$ such that $\Phi_{u_i}$ is conjugacy injective at $g$.}
\end{defn}

\medskip

{\bf Problem 1.} Let $m\geqslant 2$.

a) Is the set of words $\{x_1^sx_2^s\dots x_m^s\,|\, s\in \mathbb{N}\}$ a test set for acylindrically hyperbolic groups with natural restrictions?

b) Is the set of words $\{x_1^{k_1}x_2^{k_2}\dots x_m^{k_m}\,|\, k_i\in \mathbb{N},\, i=1,\dots ,m \}$ a test set for acylindrically hyperbolic groups with natural restrictions?


{\bf Problem 2.}
Describe retracts in $\pi_1(S)$ for closed surfaces $S$. Are all of them geometric?

\section{General case}

Let ${N_1}=(k_1,l_1,m_1,k_2,l_2,m_2,s,p_1,q_1,r_1,p_2,q_2,r_2,t)\in \mathbb{N}^{14}$ be a 14-tuple of natural numbers. We set
$$
W_1(x_1,x_2,x_3,x_4;N_1)=
\Bigl(\bigl(x_1^{k_1}x_2^{l_1}\bigr)^{m_1}\bigl(x_3^{k_2}x_4^{l_2}\bigr)^{m_2}\Bigr)^{s}
\Bigl(\bigl(x_1^{p_1}x_2^{q_1}\bigr)^{r_1}\bigl(x_3^{p_2}x_4^{q_2}\bigr)^{r_2}\Bigr)^{t}.
$$

For a sequence $(N_i)_{i\in \mathbb{N}}$ of points in $\mathbb{N}^{14}$, we recursively define
words $W_i$ in variables $x_1,\dots,x_{4i}$:
$$
\begin{array}{ll}
& W_{i+1}(x_1,\dots ,x_{4i+4};N_1,\dots ,N_{i+1})\vspace*{2mm}\\
= & W_1\bigl(W_i(x_1,\dots ,x_{4i};N_1,\dots,N_{i}),x_{4i+2},x_{4i+3},x_{4i+4}; N_{i+1}\bigr).
\end{array}
$$


\section{A geometric piece}

We consider the equation $u^nv^m=a^nb^m$, where $u=r^{-1}a^kr$ and $v=s^{-1}b^ls$.
Let $ABCDEFGH$ be the geodesic 8-gon with $A=1$ and sides $[A,B], [B,C], \dots , [H, A]$ whose
labels are consecutive syllables of the word $a^{-n}r^{-1}a^{nk}rs^{-1}b^{lm}sb^{-m}$.

\medskip

{\it Case 1.} Suppose that the following two conditions are satisfied:

1) no one $a$-component of $[B,C]\cup [D,E]$ is connected with
$[A,B]\cup [C,D]$.

2) no one $b$-component of $[E,F]\cup [G,H]$ is connected with
$[F,G]\cup [H,A]$.

\medskip

We fix a geodesic segment $[A,E]$. We should additionally assume that $[D,E]$ contains at least one edge.
Consider the 5-gon $ABCDE$.

Note that $\widehat{d}_{\lambda_1}([C,D])$ is large
(we assume that $a$ has infinite order and $n$ is large).
Therefore the $a$-component $[C,D]$ is connected to an $a$-component $[A_1,E_1]$ of $[A,E]$.

Consider the 4-gon $ABCA_1$.
Since $\widehat{d}_{\lambda_1}([A,B])$ is large, $[A,B]$ is connected to an $a$-component of $[A,A_1]$.
Then the first edge of $[A,A_1]$ (which is the first edge of $[A,E]$) is an $a$-edge.

Analogously, the first edge of $[A,E]$ is a $b$-edge. A contradiction.

\medskip

{\it Case 2.}

\section{General case}

Let $G=\langle f_1,\dots,f_n\,|\, r_1,\dots,r_k\rangle$ and $H=\langle h_1,\dots,h_m\rangle$.

Let $x_1,\dots ,x_s$ be variables, where $s\geqslant 2$, and let $n_1,\dots,n_s$ be integer numbers.
We define a word $W(x_1,\dots,x_s; n_1,\dots,n_s)$ recursively by setting $W(x_1,x_2; n_1,n_2)=x_1^{n_1}x_2^{n_2}$ and
$$
W(x_1,\dots,x_s; n_1,\dots,n_s) =\bigl(W(x_1,\dots,x_{s-1})\bigr)^{n_{s-1}}x_s^{n_s}.
$$

\section{Old}

Before we consider the general case, we illustrate the main idea in the special case
of Proposition~\ref{prop 2.1}. We use there Lemma~\ref{lem 2.2} which we prove later.

\begin{prop}\label{prop 2.1}
Suppose that $G$ is acylindrically hyperbolic and has a finite presentation with two relations:
$G=\langle g_1,\dots, g_n\,|\, R_1,R_2\rangle$. Suppose that $H=\langle a_1,a_2\rangle$ is a two-generated verbally closed subgroup of $G$, where $a_1,a_2$ are loxodromic elements in $G$.
Then $H$ is a retract of~$G$.
\end{prop}

\medskip

{\it Proof.} Let $F=F(x_1,\dots,x_n)$ be a free group of rank $n$.
Let $w_1,w_2,v_1,v_2$ be words in $F$ such that $w_i(g_1,\dots,g_n)=R_i$ and $v_i(g_1,\dots,g_n)=a_i$ for $i=1,2$.
Consider the following equation in variables $x_1,\dots,x_n$:
$$
\begin{array}{ll}
  &
\Bigl(\bigl(a_1^{k_1}a_2^{l_1}\bigr)^{m_1}\bigl(a_1^{k_2}a_2^{l_2}\bigr)^{m_2}\Bigr)^{s}
\Bigl(\bigl(a_1^{p_1}a_2^{q_1}\bigr)^{r_1}\bigl(a_1^{p_2}a_2^{q_2}\bigr)^{r_2}\Bigr)^{t}\vspace*{2mm}\\
= &
\Bigl(\bigl(v_1^{k_1}v_2^{l_1}\bigr)^{m_1}\bigl((v_1w_1)^{k_2}(v_2w_2)^{l_2}\bigr)^{m_2}\Bigr)^{s}
\Bigl(\bigl(v_1^{p_1}v_2^{q_1}\bigr)^{r_1}\bigl((v_1w_1)^{p_2}(v_2w_2)^{q_2}\bigr)^{r_2}\Bigr)^{t}.
\end{array}
$$

This equation has the solution $g_1, \dots ,g_n$ in $G$, hence it has a solution $h_1,\dots, h_n$ in $H$.
We set $V_i=v_i(h_1,\dots,h_n)$, $W_i=w_i(h_1,\dots,h_n)$. Then we have

$$
\begin{array}{ll}
  &
\Bigl(\bigl(a_1^{k_1}a_2^{l_1}\bigr)^{m_1}\bigl(a_1^{k_2}a_2^{l_2}\bigr)^{m_2}\Bigr)^{s}
\Bigl(\bigl(a_1^{p_1}a_2^{q_1}\bigr)^{r_1}\bigl(a_1^{p_2}a_2^{q_2}\bigr)^{r_2}\Bigr)^{t}\vspace*{2mm}\\
= &
\Bigl(\bigl(V_1^{k_1}V_2^{l_1}\bigr)^{m_1}\bigl((V_1W_1)^{k_2}(V_2W_2)^{l_2}\bigr)^{m_2}\Bigr)^{s}
\Bigl(\bigl(V_1^{p_1}V_2^{q_1}\bigr)^{r_1}\bigl((V_1W_1)^{p_2}(V_2W_2)^{q_2}\bigr)^{r_2}\Bigr)^{t}.
\end{array}
$$

Let $U$ be the left side of this equation. Observe that $U\in H$. By Lemma~\ref{lem 2.2}, there exists $\alpha\in \mathbb{Z}$ such that
$$V_1=V_1W_1=a_1^{U^{\alpha}}\hspace*{2mm}{\text{\rm and}}\hspace*{2mm} V_2=V_2W_2=a_2^{U^{\alpha}}.\eqno{(1)}$$

In particular, $W_1=W_2=1$ and therefore there is a homomorphism $\varphi: G\rightarrow H$, sending $g_i$ to $h_i$, $i=1,\dots,n$. Moreover, the homomorphism $\widehat{U^{\alpha}}\circ \varphi:G\rightarrow H$ is a retraction,
since
$$
\begin{array}{ll}
\widehat{(U^{\alpha})^{-1}}\circ \varphi(a_i)& =\widehat{(U^{\alpha})^{-1}}\circ \varphi(v_i(g_1,\dots,g_n))\vspace*{2mm}\\
& = \widehat{(U^{\alpha})^{-1}}(v_i(h_1,\dots,h_n))=\widehat{(U^{\alpha})^{-1}}(V_i)\overset{(1)}{=}a_i.
\end{array}
$$
\hfill $\Box$

\begin{lem}\label{lem 2.2}
Let $a$ and $b$ be independent loxodromic elements.

$$
\begin{array}{ll}
  &
\Bigl(\bigl(a^{k_1}b^{l_1}\bigr)^{m_1}\bigl(a^{k_2}b^{l_2}\bigr)^{m_2}\Bigr)^{s}
\Bigl(\bigl(a^{p_1}b^{q_1}\bigr)^{r_1}\bigl(a^{p_2}b^{q_2}\bigr)^{r_2}\Bigr)^{t}\vspace*{2mm}\\
= &
\Bigl(\bigl(x_1^{k_1}x_2^{l_1}\bigr)^{m_1}\bigl(x_3^{k_2}x_4^{l_2}\bigr)^{m_2}\Bigr)^{s}
\Bigl(\bigl(x_1^{p_1}x_2^{q_1}\bigr)^{r_1}\bigl(x_3^{p_2}x_4^{q_2}\bigr)^{r_2}\Bigr)^{t}.
\end{array}
$$

Let $U$ be the left side of this equation. Then there is an integer number $\alpha$ such that
$x_1=x_3=a^{U^{\alpha}}$ and $x_2=x_4=b^{U^{\alpha}}$.

\end{lem}

\medskip

{\it Proof.}
By Lemma~\ref{lem 1.1}, there exists $\alpha\in \mathbb{Z}$ such that
$$
\Bigl(\bigl(a^{k_1}b^{l_1}\bigr)^{m_1}\bigl(a^{k_2}b^{l_2}\bigr)^{m_2}\Bigr)^{U^{\alpha}}=
\bigl(x_1^{k_1}x_2^{l_1}\bigr)^{m_1}\bigl(x_3^{k_2}x_4^{l_2}\bigr)^{m_2},\eqno{(2)}$$
$$
\Bigl(\bigl(a^{p_1}b^{q_1}\bigr)^{r_1}\bigl(a^{p_2}b^{q_2}\bigr)^{r_2}\Bigr)^{U^{\alpha}}=
\bigl(x_1^{p_1}x_2^{q_1}\bigr)^{r_1}\bigl(x_3^{p_2}x_4^{q_2}\bigr)^{r_2}.\eqno{(3)}
$$

We set
$$
V_1=\bigl(a^{k_1}b^{l_1}\bigr)^{m_1}\bigl(a^{k_2}b^{l_2}\bigr)^{m_2}\hspace*{2mm}{\text{\rm and}}\hspace*{2mm} V_2=\bigl(a^{p_1}b^{q_1}\bigr)^{r_1}\bigl(a^{p_2}b^{q_2}\bigr)^{r_2}.
$$

Applying Lemma~\ref{lem 1.1}
to (2) and (3), we get for some $\alpha_1,\alpha_2\in \mathbb{Z}$

$$
\bigl(a^{k_1}b^{l_1}\bigr)^{V_1^{\alpha_1}U^{\alpha}}=x_1^{k_1}x_2^{l_1},\hspace*{5mm}
\bigl(a^{k_2}b^{l_2}\bigr)^{V_1^{\alpha_1}U^{\alpha}}=x_3^{k_2}x_4^{l_2},\eqno{(4)}
$$

$$
\bigl(a^{p_1}b^{q_1}\bigr)^{V_2^{\alpha_2}U^{\alpha}}=x_1^{p_1}x_2^{q_1},\hspace*{5mm}
\bigl(a^{p_2}b^{q_2}\bigr)^{V_2^{\alpha_2}U^{\alpha}}=x_3^{p_2}x_4^{q_2}.\eqno{(5)}
$$

Finally, we set
$$
W_1=a^{k_1}b^{l_1},\hspace*{2mm} W_2=a^{k_2}b^{l_2},\hspace*{2mm}
W_3=a^{p_1}b^{q_1},\hspace*{2mm} W_4=a^{p_2}b^{q_2}.
$$
Applying Lemma~\ref{lem 1.1} to (4) and (5), we get for some $\beta_1,\beta_2,\beta_3,\beta_4\in \mathbb{Z}$

$$
a^{W_1^{\beta_1}V_1^{\alpha_1}U^{\alpha}}=x_1,\hspace*{2mm} b^{W_1^{\beta_1}V_1^{\alpha_1}U^{\alpha}}=x_2,\hspace*{2mm}
a^{W_2^{\beta_2}V_1^{\alpha_1}U^{\alpha}}=x_3,\hspace*{2mm} b^{W_2^{\beta_2}V_1^{\alpha_1}U^{\alpha}}=x_4,\hspace*{2mm} \eqno{(6)}
$$

$$
a^{W_3^{\beta_3}V_2^{\alpha_2}U^{\alpha}}=x_1,\hspace*{2mm} b^{W_3^{\beta_3}V_2^{\alpha_2}U^{\alpha}}=x_2,\hspace*{2mm}
a^{W_4^{\beta_4}V_2^{\alpha_2}U^{\alpha}}=x_3,\hspace*{2mm} b^{W_4^{\beta_4}V_2^{\alpha_2}U^{\alpha}}=x_4,\hspace*{2mm} \eqno{(7)}
$$

From the first equations in (6) and (7), we deduce that $W_1^{\beta_1}V_1^{\alpha_1}V_2^{-\alpha_2}W_2^{-\beta_3}$
centralizes $b$, hence $\alpha_1=\alpha_2=\beta_1=\beta_3=0$. Analogously, from the last equations in (6) and (7), we deduce $\alpha_1=\alpha_2=\beta_2=\beta_4=0$.

Then the first and the third equations in (6) imply
$x_1=x_3=a^{U^{\alpha}}$, and the second and the fourth imply $x_2=x_4=b^{U^{\alpha}}$.
\hfill $\Box$

/////////////////////////////////////////////

\begin{lem}\label{conclusion} Suppose that $a,b$ are
non-commensurable loxodromic with respect to $\mathcal{H}$ and $E_G(a)=\langle a\rangle$, $E_G(b)=\langle b\rangle$. Then there exists $N_0$ such that for all $n,m\geqslant N_0$, $n\neq m$ the following holds:

If $a^nb^m=c^nd^m$, then $c$ is conjugate to $a$ and $d$ is conjugate to $b$.
\end{lem}

\medskip

{\it Proof.} If at least one of $c,d$ is loxodromic with respect to $\mathcal{H}'$, then we get a contradiction.
Suppose that $c,d$ are elliptic with respect to $\mathcal{H}'$. Then each of $c,d$ is commensurable with
one of the elements, $a$ or $b$.

{\bf Case 1.} Suppose that they are commensurable with the same element of $\{a,b\}$, say $a$. Then we apply $\mathcal{Q}_b$ and get a contradiction.

\medskip

{\bf Case 2.} Suppose that $c$ is commensurable with $b$, and $d$ is commensurable with $a$.
Applying $\mathcal{Q}_a$ and $\mathcal{Q}_b$ and using assumption on elementary subgroups,
we conclude that $m$ is a divisor of $n$ and $n$ is a divisor of $m$.
A contradiction.

\medskip

{\bf Case 3.} Suppose that $c$ is commensurable with $a$, and $d$ is commensurable with $b$.
Applying $\mathcal{Q}_a$ and $\mathcal{Q}_b$, we conclude that $c$ is conjugate to $a$  and $d$ to $b$.
\hfill $\Box$

\medskip

\begin{lem}\label{lem 3.6}
Let $a$ and $b$ as above. Let $n,m$ be sufficiently large different positive integers. Suppose that
$c,d\in G$ have the following properties:

\begin{enumerate}
\item[(1)] $a^nb^m=c^nd^m$,

\item[(2)] $c$ and $d$ are elliptic with respect to $X\cup \mathcal{H}'$,
where $\mathcal{H}'=\mathcal{H}\cup H_a\cup H_b$.
\end{enumerate}

\noindent
Then $c$ is conjugate to $a$ and $d$ is conjugate to $b$.
\end{lem}

\medskip

{\it Proof.} 
We use Corollary~\ref{quasimorphisms} about extensions of quasimorphisms.
Recall that $a,b$ are loxodromic with respect to $\mathcal{H}$. Let $Ell$ be the set of elliptic elements with respect to $\mathcal{H}$.

Let $\mathcal{Q}_a:G\rightarrow \mathbb{R}$ be a homogenous quasimorphism, extending $q(a)=1$, $q(b)=0$, $q(Ell)=0$.
Let $\mathcal{Q}_b:G\rightarrow \mathbb{R}$ be a quasimorphism, extending $q(b)=1$, $q(a)=0$, $q(Ell)=0$.

\medskip

{\it Case 1.} Suppose that $c,d$  are both conjugate into $H_{\lambda}$, where $\lambda\neq a,b$,
i.e. $c,d$ are elliptic with respect to $\mathcal{H}$.
Then we use $\mathcal{Q}_a$ and get a contradiction for large $n$.

\medskip

{\it Case 2.} Suppose that one of $c$ or $d$ is conjugate into $\mathcal{H}_a\cup \mathcal{H}_b$
and the other one is conjugate into $H_{\lambda}$, where $\lambda\neq a,b$.
Then we use again $\mathcal{Q}_a$, or $\mathcal{Q}_b$ and get a contradiction for large $n$.

\medskip

{\it Case 3.} Suppose that both $c$ and $d$ are conjugate into $\mathcal{H}_a\cup \mathcal{H}_b$.
If they both conjugate into $\mathcal{H}_a$ or both conjugate into $\mathcal{H}_b$,
we use again $\mathcal{Q}_a$, or $\mathcal{Q}_b$.

Thus, two subcases remain:

\begin{enumerate}

\item[(a)] $c$ is conjugate to $a^s$ for some $s$ and
 $d$ is conjugate to $b^t$ for some $t$.

\item[(b)] $c$ is conjugate to $b^s$ for some $s$ and
 $d$ is conjugate to $a^t$ for some $t$.

\medskip

Applying $\mathcal{Q}_a$ and $\mathcal{Q}_b$, we get $s=t=1$ in both subcases.
Then the subcase $(a)$ is the desired statement.

\end{enumerate}

////////////////////////////////////

Suppose that $w(a,b)c^n=e$ for some elliptic element $e$. Then there exist $x,y\in G$
such that $e=x^{-1}yx$ and $|y|_{X\cup \mathcal{H}'}\leqslant 8\delta+1$.
Using a bit of hyperbolic geometry, we can find $q,d\in G$ such that $e=q^{-1}dq$, $|d|\leqslant 20\delta+4$, and $|e|\approx_1 2|q|$.

We fix $p,c_1\in G$ such that

-- the element $c_1$ is shortest in the conjugacy class of $c$,

-- the element $p$ is shortest with the property $c^n=p^{-1}c_1^np$,

\noindent

\vspace*{-20mm}
\hspace*{-8mm}
\includegraphics[scale=0.65]{Bild_4.eps}

\vspace*{-110mm}

\begin{center}
Fig. 3
\end{center}

\bigskip


First we consider the case where $n$ is even. For large $n=n(k)$, we have $e\neq 1$, and hence $d\neq 1$.

{\it Claim.} Let $M$ be the middle point of the geodesic $[BD]$. Then $d(E,M)\leqslant k+2\delta+2$
and $d(F,M)\leqslant k+2\delta+2$.

\medskip

{\it Proof.} \hfill $\Box$

By Lemma~\ref{blue}, $BC_1C_2\dots C_nD$ is a $(\varkappa,\epsilon)$-quasigeodesic.
Then $|MC_{\frac{n}{2}}|\leqslant \theta$ for some $\theta=\theta(\delta,\varkappa,\epsilon)$.
We have $|EC_{\frac{n}{2}}|\leqslant \varkappa_1$ for some constant $\varkappa_1$.

The ends of the quasigeodesic $C_{\frac{n}{2}}\dots C_nD$ and the geodesic $ED$ are $\varkappa_1$-close.
By Lemma~\ref{Hausdorff}, each is at Hausdorff distance at most $\tau_1$ from each other.

Therefore, for each point $C_{\frac{n}{2}+k}$, there exists a point $X_k\in [ED]$ such that
$|X_kC_{\frac{n}{2}+k}|\leqslant \tau_1$. Analogously, there  exists a point $Y_k\in [FG]$ such that
$|Y_kC_{\frac{n}{2}-k}|\leqslant \tau_2$. Let $\tau$ be the maximum of $\tau_1,\tau_2$.

We have $||EX_k|-|c_1^k||\leqslant 2\tau$ and $||FY_k|-|c_1^k||\leqslant 2\tau$.
Let $x_k$ be the label of $EX_k$, and let $y_k$ be the label of $FY_k$.
Then $x_k=y_kz_k$ for some $z_k\in G$ with $|z_k|\leqslant 4\tau$.

\medskip

$\bullet$ Reading the label of the closed path $EC_{\frac{n}{2}}C_{\frac{n}{2}+k}X_kE$, we have
$\alpha c_1^k \beta_k x_k^{-1}=1$ for short $\alpha, \beta_k$.

$\bullet$ Reading the label of the closed path $FC_{\frac{n}{2}}C_{\frac{n}{2}-k}Y_kF$, we have
$\gamma c_1^{-k}\delta_ky_k^{-1}=1$ for short $\gamma,\delta_k$.

We have seen above that $x_k=y_kz_k$ for short $z_k$. Combining, we have $c_1^k(\gamma^{-1}\alpha)c_1^k=\delta_kz_k\beta_k^{-1}$. By Lemma~\ref{red}, we have $\gamma^{-1}\alpha\in E_G(c_1)$ and $E_G(c_1)$ contains an infinite dihedral subgroup. But $\gamma^{-1}\alpha=\gamma^{-1}d^{-n}\gamma$.


The case of even $n$ is analogous.\hfill $\Box$

///////////////////////////

\begin{lem}\label{elliptic-ell} Suppose that $a,b$ are
non-commensurable loxodromic with respect to $\mathcal{H}$.
Then there exists $N_0$ such that for all $n,m\geqslant N_0$ the following holds:

Suppose that $a^nb^m=c^nd^m$, where $c,d$ are elliptic with respect to $\mathcal{H}'=\mathcal{H}\cup H_a\cup H_b$.
Then either $c$ is commensurable with $a$ and $d$ is commensurable with $b$, or
$c$ is commensurable with $b$ and $d$ is commensurable with $a$.
\end{lem}

\medskip

{\it Proof.}
We use Corollary~\ref{quasimorphisms} about extensions of quasimorphisms.
Let $Ell$ be the set of elliptic elements with respect to $\mathcal{H}$.

Let $\mathcal{Q}_a:G\rightarrow \mathbb{R}$ be a homogenous quasimorphism, extending $q(a)=1$, $q(b)=0$, $q(Ell)=0$.
Let $\mathcal{Q}_b:G\rightarrow \mathbb{R}$ be a quasimorphism, extending $q(b)=1$, $q(a)=0$, $q(Ell)=0$.

\medskip

Consider four cases:

{\bf Case 1.} Suppose that $c,d$ are elliptic with respect to $\mathcal{H}$. Then apply $\mathcal{Q}_a$
and get a contradiction.

\medskip

{\bf Case 2.} Suppose that $c$ is elliptic and $d$ is loxodromic with respect to $\mathcal{H}$. Then, by Lemma~\ref{elliptic-lox},
$d$ is commensurable with $a$ or with $b$. If $d$ is commensurable with $a$, then apply $\mathcal{Q}_b$ and get a contradiction.
If $d$ is commensurable with $b$, then apply $\mathcal{Q}_a$ and get a contradiction.

\medskip

{\bf Case 3.} Suppose that $d$ is elliptic and $c$ is loxodromic with respect to $\mathcal{H}$. This case is analogous to Case 3.

\medskip

{\bf Case 4.} Suppose that $c,d$ are loxodromic with respect to $\mathcal{H}$. Then, by Lemma~\ref{elliptic-lox},
$c,d$ are commensurable with elements of $\{a,b\}$. If they are commensurable with $a$, then we apply $\mathcal{Q}_b$,
and get a contradiction. If they are commensurable with $b$, then we apply $\mathcal{Q}_a$,
and get a contradiction.
Thus, $c$ and $d$ are commensurable with different elements of $\{a,b\}$. \hfill $\Box$

\subsection{Relatively elliptic elements}

\begin{lem}\label{Ell}
Suppose that $\{E_G(a_1),\dots ,E_G(a_k)\}\hookrightarrow_h (G,X)$.\marginpar{\tiny infinite cyclic}
Then there exists a constant $C>0$ such that the following holds:

For any two elliptic (with respect to $X'=X\cup \mathcal{E}$) elements $e_1,e_2\in G$,\break
for any expression of their product $e_1e_2$ as a word $w$ in the alphabet $X'$,
and for any $i\in \{1,\dots ,k\}$
such that $e_1$ and $e_2$ are not conjugate to a power of $a_i$,
we have $$|w|_{X'}\geqslant C\log_{a_i}(w),$$
where $\log_{a_i}:F(X\bigsqcup \mathcal{E})\rightarrow \mathbb{Z}$ is the homomorphism
which extends the following map on the alphabet $X\bigsqcup \mathcal{E}$:

$$
\log_{a_i}(y)=
\begin{cases}
\ell, & {\text{\rm if}}\hspace*{2mm} y\in X\hspace*{2mm}{\text{\rm and}}\hspace*{2mm}y
=a_i^{\ell}\hspace*{2mm}{\text{\rm in}}\hspace*{2mm} G,\\
\ell, & {\text{\rm if}}\hspace*{2mm} y\in \mathcal{E}\hspace*{2mm}{\text{\rm and}}\hspace*{2mm}y
=a_i^{\ell}\hspace*{2mm}{\text{\rm in}}\hspace*{2mm} G,\\
0, & {\text{\rm otherwise.}}
\end{cases}
$$
\end{lem}

\medskip

{\it Proof.} Denote $I=\{1,\dots,k\}$ and let $i\in I$ be an element such that $e_1$ and $e_2$ are not conjugate to a power of $a_i$.
By Corollary~\ref{quasimorphisms}, there exists a homogeneous quasi-morphism $q_i:G\rightarrow \mathbb{R}$
such that $q_i(a_i)=1$,  $q_i(a_j)=0$ for $j\in I\setminus \{i\}$, and $q(e)=0$ for
all $e\in {\text{\rm Ell}}(G,X\cup \mathcal{E})\setminus [\mathcal{E}]$. In particular, $q_i(e_1)=q_i(e_2)=0$.
We consider the identity $e_1e_2=w$ in $G$, where $w$ as in Lemma.
Then we have $$|q_i(w)|\leqslant q_i(e_1)+q_i(e_2)+ D(q_i)=D(q_i)\eqno{(6.1)}$$
On the other hand (using the definition of $q$ in the proof of Osin-Hull), we have
$$
|q_i(w)- \log_{a_i}(w)|\leqslant (|w|_{X'}-1)D(q).\eqno{(6.2)}
$$
It follows from (6.1) and (6.2) that $$|w|_{X'}\geqslant \frac{1}{D(q_i)}\log_{a_i}(w).$$
\hfill $\Box$

///////////////

\subsection{Special loxodromic elements}

\begin{defn}
{\rm Let $\Gamma(G,X)$ be acylindric and hyperbolic and acylindric.
An element $g\in G$ is called {\it special loxodromic} with respect to $X$ if it is loxodromic with respect to
$X$ and $E_G(g)=\langle g\rangle$.}
\end{defn}

\begin{lem}
Suppose that $H$ is an acylindrically hyperbolic group which does not normalize a nontrivial finite subgroup.
We fix the hyperbolic structure $E\hookrightarrow_h (H,X)$.
Let $a$ and $b$ be two special loxodromic and non-commensurable elements in $H$ with respect to $(X,\mathcal{E})$. Then, for each $m,k\in \mathbb{N}$, there exist $(k+1)$ different numbers $\alpha,\beta_1,\dots,\beta_k\in m\mathbb{N}$ such that the elements
$a^\alpha b^{\beta_i}$, $i=1,\dots ,k$ are special loxodromic and pairwise non-commensurable.
\end{lem}

\medskip

{\it Proof.} Using notation from~\cite{Osin_2}, we have $K(H)=1$.  Since $a$ and $b$ are non-commensurable and $E_H(a)=\langle a\rangle$, $E_H(b)=\langle b\rangle$, we have $E_H(a)\cap E_H(b)=1$,
$E_H(a)=E_H(a)^{+}$, and $E_H(b)=E_H(b)^{+}$.
This is a start point of the proof of~\cite[Lemma 6.18]{Osin_2} (where we set $h_1:=a^m$ and $h_2:=b^m$). According to the first two paragraphs of this proof, there exist the desired numbers $\alpha,\beta_1,\dots ,\beta_k$.\hfill $\Box$

\medskip

\begin{cor}\label{add_properties}
Suppose that $H$ is an acylindrically hyperbolic group which does not normalize a nontrivial finite subgroup.
We fix the hyperbolic structure $E\hookrightarrow_h (H,X)$.
Let $\mathcal{F}$ be a finite subset of $H$. Suppose that $a,b$ are special loxodromic and non-commensurable
elements in $H$ with respect to $(X,\mathcal{E})$. Then, for each $\ell\in \mathbb{N}$, there exist $n,m\in \ell\mathbb{N}$ such that $n\neq m$ and the element $a^nb^m$ is special loxodromic and non-commensurable with any element of $\mathcal{F}$.
\end{cor}

\section{Test words}

\begin{defn}
{\rm Let $H$ be a group, let $\mathcal{F}$ be a finite subset of $H$ and let $a_1,\dots,a_k$ be elements of $H$.
A word $W_k=W_k(x_1,\dots ,x_k)$ is called an
{\it $(a_1,\dots,a_k;\mathcal{F})$-test word} if the following two conditions hold:

\begin{enumerate}
\item[1)] For every solution $(b_1,\dots,b_k)$ of the equation
$$
W_k(a_1,\dots ,a_k)=W_k(x_1,\dots ,x_k)
$$
in $G$, there exists an integer number $\alpha\in \mathbb{Z}$ such that
$b_i=a_i^{U^{\alpha}}$ for $i=1,\dots, k$, where $U=W_k(a_1,\dots ,a_k)$.

\medskip

\item[2)] The element $W_k(a_1,\dots ,a_k)$ is special loxodromic and non-commens\-urable with elements of $\mathcal{F}$.
\end{enumerate}
}
\end{defn}
We skip $\mathcal{F}$ in the above notation if $\mathcal{F}=\emptyset$.

\medskip

\begin{rmk}\label{test_special_words_1}
{\rm Suppose that $H$ is acylindrically hyperbolic group, which does not normalise a nontrivial finite subgroup.
By Lemma~\ref{lem 2.2}, if $\mathcal{F}$ is a finite subset of $H$ and if $a_1,a_2,a_3$ are pairwise non-commensurable special loxodromic elements in $H$, then
there exist integer numbers $k_1$, $l_1$, $m_1$, $k_2$, $l_2$, $m_2$, $s$, $p$, $q$, $t$ such that
$$
W_3=\Bigl(\bigl(x_1^{k_1}x_3^{l_1}\bigr)^{m_1}\bigl( x_2^{k_2}x_3^{l_2}\bigr)^{m_2}\Bigr)^{s}
(x_2^p(x_3y_3)^q)^t
$$
is an $(a_1,a_2,a_3,1;\mathcal{F})$-test word.
}
\end{rmk}

The aim of this section is to prove the following general proposition.
In Section~8 we use this proposition in the case $\mathcal{F}=\emptyset$. The reason, why
we formulate this proposition in the general case (i.e., for any finite $\mathcal{F}$) is that we cannot prove it in the special case without this generalization.

\begin{prop}\label{test_special_words}
Let $H$ be an acylindrically hyperbolic group, which does not normalise a nontrivial finite subgroup.
For any finite subset $\mathcal{F}$, for any $k\geqslant 3$ and any pairwise non-commensurable special loxodromic elements
$a_1,\dots,a_{k}$ in $H$, there exists a word $W_k(x_1,\dots,x_k,y_3,\dots,y_k)$
which is an $(a_1,\dots,a_k,\underbrace{1,\dots,1}_{k-2};\mathcal{F})$-test word.
\end{prop}

This proposition follows by induction from Remark~\ref{test_special_words_1} (basis of induction for $k=3$)
and Lemma~\ref{test_special_words_2} (inductive step).

\medskip

\begin{lem}\label{test_special_words_2} Let $a_1,\dots,a_{k+1}$ be pairwise non-commensurable special loxodromic elements in $H$ and let
$\mathcal{F}$ be a finite subset of $H$.
Suppose that\break
$W_k(x_1,\dots,x_k,y_3,\dots,y_k)$ is an
$(a_1,\dots,a_k,\underbrace{1,\dots,1}_{k-2};\{a_k,a_{k+1}\})$-test word.
Then there exist integer numbers $k_1$, $l_1$, $m_1$, $k_2$, $l_2$, $m_2$, $s$, $p$, $q$, $t$
such that
$$
W_{k+1}(x_1,\dots,x_{k+1},y_3,\dots,y_{k+1})=\Bigl(\bigl(W_k^{k_1}x_{k+1}^{l_1}\bigr)^{m_1}\bigl( x_k^{k_2}x_{k+1}^{l_2}\bigr)^{m_2}\Bigr)^{s}
(x_k^p(x_{k+1}y_{k+1})^q)^t
$$
is an $(a_1,\dots,a_{k+1},\underbrace{1,\dots ,1}_{k-1};\mathcal{F})$-test word.
\end{lem}

\medskip

{\it Proof.} Denote $A=W_k(a_1,\dots,a_k,1,\dots,1)$.
Since the elements  $A$, $a_k$, $a_{k+1}$ are special loxodromic and pairwise non-commensurable, they satisfy the assumption of Lemma~\ref{lem 2.2}. By this lemma, there
exist integer numbers $k_1$, $l_1$, $m_1$, $k_2$, $l_2$, $m_2$, $s$, $p$, $q$, $t$
such that the word
$$
\Bigl(\bigl(X^{k_1}x_{k+1}^{l_1}\bigr)^{m_1}\bigl( x_k^{k_2}x_{k+1}^{l_2}\bigr)^{m_2}\Bigr)^{s}
(x_k^p(x_{k+1}y_{k+1})^q)^t\eqno{(7.1)}
$$
in variables $(X,x_k,x_{k+1},y_{k+1})$ is an $(A,a_k,a_{k+1},1;\mathcal{F})$-test word.

Now we prove that the word $W_{k+1}$ defined in Lemma~\ref{propos} for these exponents
is the required test word.
Suppose that for some elements $b_1,\dots,b_{k+1},c_3,\dots,c_{k+1}$ in $H$ we have
$$
\begin{array}{ll}
& \Bigl(\bigl(W_k^{k_1}(a_1,\dots,a_k,1,\dots ,1)a_{k+1}^{l_1}\bigr)^{m_1}\bigl( a_k^{k_2}a_{k+1}^{l_2}\bigr)^{m_2}\Bigr)^{s}
(a_k^pa_{k+1}^q)^t\vspace*{2mm}\\
= &
\Bigl(\bigl(W_k^{k_1}(b_1,\dots,b_k,c_3,\dots ,c_k)b_{k+1}^{l_1}\bigr)^{m_1}\bigl( b_k^{k_2}b_{k+1}^{l_2}\bigr)^{m_2}\Bigr)^{s}
(b_k^p(b_{k+1}c_{k+1})^q)^t.
\end{array}
$$
and let $U$ be the left side of this equation.

Denote $B:=W_k(b_1,\dots,b_k,c_3,\dots ,c_k)$ and write this equation shorter:
$$
\begin{array}{ll}
& \Bigl(A^{k_1}a_{k+1}^{l_1}\bigr)^{m_1}\bigl( a_k^{k_2}a_{k+1}^{l_2}\bigr)^{m_2}\Bigr)^{s}
(a_k^pa_{k+1}^q)^t\vspace*{2mm}\\
= &
\Bigl(B^{k_1}b_{k+1}^{l_1}\bigr)^{m_1}\bigl( b_k^{k_2}b_{k+1}^{l_2}\bigr)^{m_2}\Bigr)^{s}
(b_k^p(b_{k+1}c_{k+1})^q)^t.
\end{array}
$$

Since the word in (7.1) is an  $(A,a_k,a_{k+1},1;\mathcal{F})$-test word, we conclude:

1) there exists $\alpha\in \mathbb{Z}$ such that
$$
B =A^{U^{\alpha}},
\eqno{(7.2)}
$$

$$
b_k =a_k^{U^{\alpha}},\hspace*{3mm}
b_{k+1}=a_{k+1}^{U^{\alpha}},\hspace*{3mm}\\
c_{k+1}=1;
\eqno{(7.3)}
$$

2) The element $U$ is special loxodromic and non-commensurable with any element of $\mathcal{F}$.

From (7.2) we deduce
$$
W_k\bigl(b_1^{U^{-\alpha}},\dots,b_k^{U^{-\alpha}},c_3^{U^{-\alpha}},\dots ,c_k^{U^{-\alpha}}\bigr) =W_k(a_1,\dots,a_k,1,\dots ,1).
$$
Recall that $A$ is the right side of this equation. Since $W_k$ is an $(a_1,\dots,a_k,\underbrace{1,\dots,1}_{k-2})$-test word, there exists $\beta\in \mathbb{Z}$
such that

$$
\begin{array}{lll}
b_1^{U^{-\alpha}} & \hspace*{-1mm}=\hspace*{-1mm} & a_1^{A^{\beta}},\vspace*{2mm}\\
\vdots & & \vspace*{2mm}\\
b_k^{U^{-\alpha}} & \hspace*{-1mm}=\hspace*{-1mm} & a_k^{A^{\beta}},
\end{array}
\eqno{(7.4)}
$$
and
$$
c_3=\dots =c_k=1. \eqno{(7.5)}
$$
From the first equation in (7.3) and the last equation in (7.4), we deduce that $a_k^{A^{\beta}}=a_k$.
Since $A$ is non-commensurable with $a_k$, we have $\beta=0$.
Then (7.3)-(7.5) complete the proof.\hfill $\Box$

\section{Proof of the main theorem}

Let $G$ be acylindrically hyperbolic with respect to a generating set $X$.
An element $g\in G$ is called {\it special} if it is loxodromic and $E_G(g)=\langle g\rangle$,
Equivalently, $g\in G$ is called {\it special} if $E_G(g)=\langle g\rangle$ and $E_G(g)\hookrightarrow_h (G,Y)$ for some $X\subset Y\subset G$.

\medskip

The following lemma in the cases where $H$ is a relatively hyperbolic group was proved in\marginpar{\tiny $H$ and $G$}
\cite[Lemma 3.2]{OT}.  Its proof can be easily adopted to the general case, where $H$ is
acylindrically hyperbolic.

\begin{lem}\label{special-special}
Let $H$ be an acylindrically hyperbolic group with respect to a generating set $X$ and let $h\in H$ be
a special element. Then
for every $a\notin E_H(h)$, there exists a positive integer $n_0$ such that for every $n>n_0$ the element $g=ah^n$ is special loxodromic.
\end{lem}

\medskip

\begin{lem}
Suppose that $H$ is an $l$-generated (for some $l\geqslant 2$) acylindrically hyperbolic group
without nontrivial finite normal subgroups.
Then $H$ can be generated by a $(l+1)$ pairwise non-commensurable special loxodromic elements.
\end{lem}

\medskip

{\it Proof.}
Recall that in~\cite{Osin_2}, $K(H)$ denotes the maximal nontrivial normal finite subgroup of $H$. In our situation we have $K(H)=1$.
Then, by~\cite[Lemma 6.18]{Osin_2}, there exists a generating subset $Z\subseteq H$ and a loxodromic element $h\in H$ with respect to $Z$ such
that $E_G(h)=\langle h\rangle\times K(H)=\langle h\rangle$. In other words, $h$ is a special loxodromic element of $H$ with respect to $Z$.

Let $H=\langle b_1,\dots, b_l\rangle$ and let $a_1,\dots ,a_k$ be all elements of the set $\{b_1,\dots, b_l\}$ which are not powers of $h$.
Then $H=\langle h, a_1,\dots, a_k\rangle$.
By~\cite[Corollary 6.12]{Osin_2},
for each $i=1,\dots,k$, there exists a subset $A_i\subseteq \{a_ih^m\,|\, m\in \mathbb{Z}\}$
consisting of $i$ pairwise non-commensurable loxodromic elements with respect to $(Y\cup E_G(h))$.\marginpar{\tiny $Y$ is the original set of gen}
Using Lemma~\ref{special-special}, we may additionally assume that all elements of all $A_i$ are special.\marginpar{\tiny special}
Let $A_1=\{a_1h^{n_1}\}$.
Since $|A_2|=2$, there exists $a_2h^{n_2}\in A_2$ which is non-commensurable with $a_1h^{n_1}$.
Since $|A_3|=3$, there exists $a_3h^{n_3}\in A_3$ which is non-commensurable with $a_1h^{n_1}$ and $a_1h^{n_2}$.
By induction, we can find $n_1,\dots ,n_k$ such that the elements $a_1h^{n_1},\dots ,a_kh^{n_k}$
are pairwise non-commensurable.

All of them are non-commensurable with $h$ since they are loxodromic with respect to $Y\cup E_G(h)$, and
$h$ is elliptic since $|h^n|_{Y\cup E_G(h)}=1$ for any $n\in \mathbb{Z}\setminus \{0\}$.
All of them, including $h$, are loxodromic with respect to $Y$. These elements are chosen to be special.

At the end we have $k+1$ generators. If $k=1$, then $H=\langle h,a_1h^{k_1}\rangle$ and
we can restart the proof with the generators $h,a_1h^{k_1}, a_1$.
\hfill $\Box$


\medskip

\begin{thm}\label{prop 4.1}
Suppose that $G$ is a finitely presented group and $H$ is a finitely generated acylindrically hyperbolic subgroup of $G$ such that $H$ does not normalise a nontrivial finite subgroup of $H$.
If $H$ is verbally closed in $G$, then $H$ is a retract of~$G$.
\end{thm}

\medskip

{\it Proof.} Let $G=\langle g_1,\dots,g_n\,|\, R_1,\dots,R_m\rangle$ and $H=\langle a_1,\dots, a_k\rangle$.
We may assume that $k\geqslant 2$. Then every finite set of independent loxodromic elements of $H$ is contained in an infinite set of independent loxodromic elements of $H$.
Therefore we may enlarge $k$ if necessary and assume that $k=m+2$.

Let $v_i$ and $u_j$ be words in variables $x_1,\dots,x_n$ such that  $a_i=v_i(g_1,\dots,g_n)$ for $i=1,\dots,k$ and $R_j=u_j(g_1,\dots,g_n)$  for $j=1,\dots,m$.
Consider the following equation, where $W_k$ is an $(a_1,\dots,a_k,\underbrace{1,\dots,1}_{m})$-test word:

$$
\begin{array}{ll}
& W_k(a_1,\dots,a_k,1,\dots,1)\vspace*{2mm}\\
= & W_k(v_1(x_1,\dots,x_n),\dots,v_k(x_1,\dots,x_n),
u_1(x_1,\dots,x_n),\dots,u_m(x_1,\dots,x_n)).
\end{array}
$$

This equation has the solution $g_1, \dots ,g_n$ in $G$, hence it has a solution $h_1,\dots, h_n$ in $H$.
We set $V_i=v_i(h_1,\dots,h_n)$, $U_j=u_j(h_1,\dots,h_n)$. Then we have

$$
 W_k(a_1,\dots,a_k,1,\dots,1)= W_k(V_1,\dots,V_k,U_1,\dots,U_m).
$$

Let $U$ be the left side of this equation. Observe that $U\in H$. By Lemma~\ref{lem 2.2}, there exists $\alpha\in \mathbb{Z}$ such that

$$V_i=a_i^{U^{\alpha}},\hspace*{2mm} i=1,\dots,k,\eqno{(8.1)}$$
$$
U_j=1,\hspace*{2mm} j=1,\dots,m.\eqno{(8.2)}
$$

Because of (8.2), there is a homomorphism $\varphi: G\rightarrow H$, sending $g_i$ to $h_i$, $i=1,\dots,n$. Moreover, the homomorphism $\widehat{(U^{\alpha})^{-1}}\circ \varphi:G\rightarrow H$ is a retraction,
since
$$
\begin{array}{ll}
\widehat{(U^{\alpha})^{-1}}\circ \varphi(a_i)& =\widehat{(U^{\alpha})^{-1}}\circ \varphi(v_i(g_1,\dots,g_n))\vspace*{2mm}\\
& = \widehat{(U^{\alpha})^{-1}}(v_i(h_1,\dots,h_n))=\widehat{(U^{\alpha})^{-1}}(V_i)\overset{(8.1)}{=}a_i.
\end{array}
$$
\hfill $\Box$

\section{A very special case:\\ Infinite cyclic verbally closed subgroups}

\begin{prop} {\rm (see~\cite[Lemma 1]{Mazhuga})}
Let $G$ be a group, $K$ a verbal subgroup of $G$, and $\varphi:G\rightarrow G/K$ be the canonical epimorphism.
Suppose that $H$ is a subgroup of $G$. If $H$ is verbally closed in $G$, then $\varphi(H)$ is verbally closed in $\varphi(G)$.
\end{prop}

\begin{prop}
Let $G$ be a group with finitely generated abelianization $G/[G,G]$. Every torsion-free verbally closed abelian subgroup
$H$ of $G$ is a retract of $G$. (Of course, every such retraction goes through $G/[G,G]$.)
\end{prop}